\magnification=\magstep1
\input amstex
\documentstyle{amsppt}
\loadbold

\def\R{{\bold R}}
\def\C{{\bold C}}
\def\N{{\bold N}}
\def\Z{{\bold Z}}
\def\P{{\bold P}}
\def\dis{\mathop{\roman{dis}}\nolimits}
\def\id{\mathop{\roman{id}}\nolimits}
\def\dim{\mathop{\roman{dim}}\nolimits}
\def\codim{\mathop{\roman{codim}}\nolimits}
\def\Int{\mathop{\roman{Int}}\nolimits}
\def\supp{\mathop{\roman{supp}}\nolimits}
\def\Im{\mathop{\roman{Im}}\nolimits}
\def\Sing{\mathop{\roman{Sing}}\nolimits}
\def\const{\mathop{\roman{const}}\nolimits}
\def\Reg{\mathop{\roman{Reg}}\nolimits}

\def\Ker{\mathop{\roman{Ker}}\nolimits}
\def\graph{\mathop{\roman{graph}}\nolimits}

\topmatter
\title Analytic equivalence of normal crossing functions on a real analytic manifold\endtitle 
\rightheadtext{Analytic equivalence of normal crossing functions}
\author Goulwen Fichou and Masahiro Shiota\endauthor 
\address 
UFR Math\'ematiques de Rennes, Campus de Beaulieu, avenue du G\'en\'eral Leclerc, 35042 RENNES C\'edex, France and \endgraf
Graduate School of Mathematics, Nagoya University, Chikusa, Nagoya, 
464-8602, Japan\endaddress
\abstract
By Hironaka Desingularization Theorem, any real analytic function has only normal crossing singularities after a modification. We focus on the analytic equivalence of such functions with only normal crossing singularities. We prove that for such functions $C^{\infty}$ right equivalence implies analytic equivalence. We prove moreover that the cardinality of the set of equivalence classes is zero or countable.\endabstract 
\endtopmatter
\document

%%%%%%%%%%%%%%%%%%%%%%%%%%%%%%%%%%%%%%%%%%%%%%%%%%%%%%%%%%%%%%%%%%%%%%%%%%%%%%%%%%%%%%%%%%%%%%%%%%%
%%%%%%%%%%%%%%%%%%%%%%%%%%%%%%%%%%%%%%%%%%%%%%%%%%%%%%%%%%%%%%%%%%%%%%%%%%%%%%%%%%%%%%%%%%%%%%%%%%%
%%%%%%%%%%%%%%%%%%%%%%%%%%%%%%%%%%%%%%%%%%%%%%%%%%%%%%%%%%%%%%%%%%%%%%%%%%%%%%%%%%%%%%%%%%%%%%%%%%%

\head 
1. Introduction
\endhead

The classification of real analytic functions is a difficult but
fascinating topic in singularity theory. In this paper, we put our
interest on real analytic functions with only
normal crossing singularities. This case is of fundamental importance
since any analytic function becomes one with only normal 
crossing singularities after a finite sequence of blowings-up along
smooth center by Hironaka Desingularization Theorem [Hi]. 
Our goal is to establish the
cardinality of the set of equivalence classes of analytic functions with only
normal crossing singularities under analytic equivalence (theorem 3.2). 

Our first main result is theorem 3.1,(1) which asserts that
$C^{\infty}$ right equivalent real analytic functions with only normal
crossing singularities are automatically analytically right
equivalent. Its proof consists in a careful use of Cartan Theorems A
and B and Oka Theorem in order to use integration along analytic
vector fields to produce analytic isomorphisms. Theorem 3.1,(1) is a
crucial result in order to deal with cardinality issues, in particular
in view to make a reduction to the case of real analytic functions
with semialgebraic graph, called Nash functions.

The second main result (theorem 3.2) establishes the cardinality of
the set of equivalence classes of real analytic (respectively Nash)
functions with only normal crossing singularities on a compact
analytic manifold (resp. on a non-necessarily compact Nash manifold)
with respect to the analytic (resp. Nash) equivalence. To prove that
this cardinality is zero or countable, we first reduce the study to
the Nash case by theorem 3.1,(1), then from the non compact to the
compact case via Nash sheaf theory, a Nash version of Hironaka
Desingularization Theorem and a finer analysis of the normal crossing
property on a Nash manifold with corners. Finally Hardt triviality
[Ha], Artin-Mazur Theorem (see [S$_2$]) and Nash Approximation Theorems
[S$_2$], [C-R-S$_1$]
enable to achieve the proof. Note that along the way, we establish (as theorem 3.1,(3)) a $C^2$ plus semialgebraic version of theorem 3.1,(1), namely semialgebraically $C^{2}$ right equivalent Nash functions with only normal crossing singularities on a Nash manifold are Nash right equivalent (see also theorem 3.1,(2) for a $C^2$ version).

The paper is organized as follows. In section one, we recall some
definitions that are fundamental in the paper, in particular the
notion of normal crossing in the case of manifolds with corners. We
devote the second section to some preliminaries about real analytic
and Nash sheaf theory, that will be crucial tools for the proof of the
main theorems, and also a quick overview on the different topologies we will
consider on spaces of maps. Third section is dedicated to theorem
3.1,(1) and its proof, and the statement of theorem 3.1,(2) and
3.1,(3), the proof of which we postpone to section five. 
Actually, even though the statements are very similar, we need to prepare in
section four some materials for it. We prove in particular as lemma 4.6 that a normal crossing
Nash subset of a non-compact Nash manifold is trivial at infinity, and
we compactify in proposition
4.9 a Nash function with only normal crossing singularities. We
finally prove theorems 3.1,(2) and
3.1,(3) together with theorem 3.2 in the last section.

In this paper a manifold means a manifold without boundary, analytic
manifolds and maps mean real analytic ones unless otherwise specified,
and id stands for the identity map.

\subhead
1.1. Analytic functions with only normal crossing singularities
\endsubhead

\definition{Definition 1.1} Let $M$ be an analytic manifold. An analytic function with {\it only normal crossing singularities} at a point $x$ of $M$ is 
a function whose germ at $x$ is of the form $\pm x^\alpha
(=\pm\prod_{i=1}^n x_i^{\alpha_i})$ up to an additive constant, for some local analytic coordinate system $(x_1,...,x_n)$ at 
$x$ and some $\alpha=(\alpha_1,...,\alpha_n)\not=0\in\N^n$. 
If the function has only normal crossing singularities everywhere, we
say that the function has {\it only normal 
crossing singularities}. 
\enddefinition

An analytic subset of an analytic manifold is called {\it normal crossing} if it is the zero set of an 
analytic function with only normal crossing singularities. 
This analytic function is called {\it defined by} the analytic set. 
It is not unique. 
However, the sheaf of $\Cal O$-ideals {\it defined by} the analytic set is naturally defined and unique. 
We can naturally stratify a normal crossing analytic subset $X$ into analytic manifolds $X_i$ of dimension 
$i$. 
We call $\{X_i\}$ the {\it canonical stratification} of $X$. 

\subhead
1.2. Case of Nash manifolds
\endsubhead

\definition{Definition 1.2} 
A {\it semialgebraic} set is a subset of a Euclidean space which is described by finitely many equalities 
and inequalities of polynomial functions. 
A {\it Nash manifold} is a $C^\omega$ submanifold of a Euclidean space which is semialgebraic. 
A {\it Nash function} on a Nash manifold is a $C^\omega$ function with semialgebraic graph. 
A {\it Nash subset} is the zero set of a Nash function on a Nash manifold. 
(We call a germ {\it on} but not {\it at} $X$ in $M$ to distinguish the case where $X$ is a set from the 
case of a point.) 
\enddefinition

We define Nash functions with {\it only normal 
crossing singularities}, {\it
normal crossing} Nash subsets of a Nash manifold and the {\it canonical stratification} of a normal
crossing Nash subset similarly to the analytic case. 

For elementary properties of Nash manifolds and Nash
functions, we refer to [S$_2$].
As a general flavor, note that Nash
functions carry more structure than analytic or semialgebraic ones,
and therefore it is useful to dispose of approximation
results. In this paper, we will make an intensive use of the two
classical approximation theorems by Nash functions, which are quite
different in nature. The first one,
that we will refer to as Nash Approximation Theorem I, concerns the
approximation of semialgebraic $C^r$ maps by Nash maps (see [S$_2$]).
The topology we use in that case is the {\it semialgebraic $C^r$
topology} on spaces of semialgebraic $C^r$ maps (see subsection 2.3 for
an overview about topologies on spaces of maps). Note for instance that, in that topology, a semialgebraic 
$C^1$ map between semialgebraic $C^1$ manifolds close to a semialgebraic $C^1$ diffeomorphism is a 
diffeomorphism. 

\proclaim{Theorem}
(Nash Approximation Theorem I, [S$_2$]) Any semialgebraic $C^r$ map between
Nash manifolds can be approximated in the {\it semialgebraic $C^r$
topology} by a Nash map.
\endproclaim

The other one, say Nash Approximation Theorem
II, is a global version of Artin Approximation Theorem on a compact Nash
manifold.

\proclaim{Theorem}
(Nash Approximation Theorem II, [C-R-S$_1$]) Given a Nash function $F$
on $M_1\times M_2$ for a compact Nash manifold $M_1$ and a Nash
manifold $M_2$, and an analytic map $f:M_1 \to M_2$ with $F(x,f(x))=0$
for $x\in M_1$, then there exists a Nash approximation $\widetilde f:M_1 \to M_2$
of $f$ in the {\it$C^{\infty}$ topology} such that $F(x,\widetilde f(x))=0$ for $x\in M_1$.
\endproclaim

\subhead
1.3. Manifolds with corners
\endsubhead

Manifolds with corners appear naturally in the study of functions with
only normal crossing singularities. A manifold with corners is locally given by charts diffeomorphic to
$[0,\infty)^k\times \R^{n-k}$. In this paper we will consider analytic
manifold with corners as well as Nash ones. We refer to [K-S] for basics about
manifolds with corners.

The definition of the canonical stratification for manifolds can be
naturally extended to the boundary of an analytic manifold with
corners. However, concerning the notions of singularity and
normal crossings, we really need to adapt the definitions.

\definition{Definition 1.3} 
Let $f$ be an analytic function on analytic manifold with corners $M$. 
We say $f$ is {\it singular} at a point $x_0$ of $\partial M$ if the restriction of $f$ to the stratum 
of the canonical stratification of $\partial M$ containing $x_0$ is
singular at $x_0$.
\enddefinition

Note in particular that with such a definition, $f$ is singular at
points of the stratum of dimension 0 of the canonical stratification
of $\partial M$. This remark will be of importance when dealing with
proofs by induction.

To define a function with {\it only normal crossing singularities} on
a manifold with corners $M$, we need to extend $M$ beyond the
corners. More precisely, we can construct an analytic manifold $M'$ which contains $M$ and is of the same dimension by extending 
a locally finite system of analytic local coordinate neighborhoods of $M$. We call $M'$ an {\it analytic manifold extension} of $M$. 
In the same way, shrinking $M'$ if necessary we obtain a normal crossing analytic subset $X$ of $M'$ 
such that $\Int M$ is a union of some connected components of $M'-X$, and $f$ is extended to an analytic 
function $f'$ on $M'$. 

\definition{Definition 1.4} 
We say that $f$ has {\it only normal crossing singularities} if $f|_{\Int M}$ does so and if the germ of $(f-f(x_0))\phi$ at each point $x_0$ of $X$ has only normal 
crossing singularities, for $\phi$ an analytic function on $M'$ defined by $X$.

\enddefinition

Now we can define, similarly to the case without corners, a {\it normal crossing} analytic subset of $M$ and a {\it normal crossing} 
sheaf of $\Cal O$-ideals on $M$.

In the Nash case, we define analogously a {\it Nash manifold extension} of a Nash
manifold with corners, a Nash function with {\it only normal crossing
singularities} on a Nash manifold with corners, a {\it normal
crossing} Nash subset of $M$ and a {\it normal crossing} 
sheaf of $\Cal N$-ideals on $M$.
\par

\head
2. Preliminaries
\endhead
We dedicate this section to some remainder on real analytic sheaf
theory, and prove similar statements in the Nash case that will be of
importance in next sections. We finish with an overview of the
different topologies on spaces of functions we will make use in that
paper, in order to explain the major differences between them.

\subhead
2.1. Real analytic sheaves
\endsubhead

In this subsection, we deal with the real analytic case of Cartan Theorems A and B, and Oka Theorem.

Let $\Cal O$ and $\Cal N$ denote, respectively, the sheaves of analytic and Nash 
function germs on an analytic and Nash manifold and let $N(M)$ denote the ring of Nash functions on a Nash manifold 
$M$. 
We write $\Cal O_M$ and $\Cal N_M$ when we emphasize the domain $M$. 
Let $f_x$, $X_x$, $v_x$ and $\Cal M
_x$ denote the germs of $f$ and $X$ at a point $x$ of $M$, the tangent vector assigned to $x$ by $v$ 
and the stalk of $\Cal M$ at $x$ for a function $f$ on an analytic (Nash) manifold $M$, a subset $X$ of $M$, a vector field $v$ on $M$ 
and for a sheaf of $\Cal O$- ($\Cal N$-) modules $\Cal M$ on $M$, respectively. 
For a compact semialgebraic subset $X$ of a Nash manifold $M$, let $\Cal N(X)$ denote the germs of Nash 
functions on $X$ in $M$, with the topology of the inductive limit
space of the topological spaces $N(U)$ endowed with the compact-open $C^\infty$ topology, where $U$ runs through the family of open semialgebraic neighborhoods 
of $X$ in $M$. 
In the same way, we define $\Cal O(X)$ for a compact semianalytic subset $X$ of an analytic manifold $M$. 
Here a {\it semianalytic} subset is a subset whose germ at each point of $M$ is described by finitely many 
equalities and inequalities of analytic function germs. 
\proclaim{Theorem 2.1}
(Cartan Theorem A) 
Let $\Cal M$ be a coherent sheaf of $\Cal O$-modules on an analytic manifold $M$. 
Then for any $x\in M$, the germ $\Cal M_x$ is equal to $H^0(M,\Cal M)\Cal O_x$. 
\endproclaim
See [G-R] for Cartan Theorems A and B in the complex case and [Ca] for the real case. Next corollary will 
be useful in this paper. 
It deals with the case where the number of local generators is uniformly bounded. 

\proclaim{Corollary 2.2}
In theorem 2.1, assume that $\Cal M_x$ is generated by a uniform
number of elements for any $x$ in $M$. 
Then $H^0(M,\Cal M)$ is finitely generated as a $H^0(M,\Cal O)$-module. 
\endproclaim

The corollary is proved in [Co] in the complex case. The real case follows from a complexification of $\Cal M$ as in [Ca].

\proclaim{Theorem 2.3}
(Cartan Theorem B)
Let $\Cal M$ be a coherent sheaf of $\Cal O$-modules on an analytic manifold $M$. 
Then $H^1(M,\Cal M)$ is equal to zero.
\endproclaim
\proclaim{Corollary 2.4}
Let $M$ be an analytic manifold and $X\subset M$ be a global analytic set---the zero set of an analytic 
function. 
Let $\Cal I$ be a coherent sheaf of $\Cal O$-ideals on $M$ such that any element of $\Cal I$ vanishes on 
$X$. 
Then any $f\in H^0(M,\Cal O/\Cal I)$ can be extended to some $F\in C^\omega(M)$, i.e., $f$ is the image 
of $F$ under the natural map $H^0(M,\Cal O)\to H^0(M,\Cal O/\Cal I)$. \par
If $X$ is normal crossing, we 
can choose $\Cal I$ to be the function germs vanishing on $X$. 
Then $H^0(M,\Cal O/\Cal I)$ consists of functions on $X$ whose germs at each point of $X$ are extensible 
to analytic function germs on $M$. 
\endproclaim
Corollary 2.4 follows from theorem 2.3 by considering the exact sequence $0\to\Cal I\to\Cal O\to\Cal I/\Cal O$. 
\proclaim{Theorem 2.5}
(Oka Theorem)
Let $\Cal M_1$ and $\Cal M_2$ be coherent sheaves of $\Cal O$-modules on an analytic manifold $M$, and 
$h:\Cal M_1\to\Cal M_2$ be an $\Cal O$-homomorphism. 
Then $\Ker h$ is a coherent sheaf of $\Cal O$-modules. 
\endproclaim
See [G-R] in the complex case. 
The real case follows from complexification [Ca] of $M,\Cal M_1,\Cal M_2$ and $h$. \par

\subhead
2.2. Nash sheaves
\endsubhead

In this subsection $M$ stands for a Nash manifold. 
A sheaf of $\Cal N$-modules $\Cal M$ on $M$ is called {\it finite} if for some finite open 
semialgebraic covering $\{U_i\}$ of $M$ and for each $i$ there exists an exact sequence $\Cal N^{m_i}|
_{U_i}\longrightarrow\Cal N^{n_i}|_{U_i}\longrightarrow\Cal M|_{U_i}\longrightarrow 0$ of 
$\Cal N$-homomorphisms, with $m_i,n_i\in\N$. 
Non-finite examples are the sheaf of $\Cal N$-ideals $\Cal I$ on $\R$ of germs vanishing on $\Z$ and 
$\Cal N/\Cal I$. \par
\proclaim{Theorem 2.6}
(Nash case of Oka Theorem) 
Let $h$ be an $\Cal N$-homomorphism between finite sheaves of $\Cal N$-modules on a Nash manifold. 
Then $\Ker h$ is finite. 
\endproclaim
\demo
{Proof}
Let $h:\Cal M_1\to\Cal M_2$ be such a homomorphism on a Nash manifold $M$. 
There exists a finite open semialgebraic covering $\{U_i\}$ of $M$
such that $\Cal M_j|_{U_i}$, for$j=1,2
$, satisfy the condition of exact sequence in the definition of a finite sheaf. Therefore it suffices to 
prove the theorem on each $U_i$. 
Now we may assume that $\Cal M_j$, for $j=1,2$, are generated by global cross-sections $\alpha_1,...,\alpha_{n_1}
$ and $\beta_1,...,\beta_{n_2}$, respectively, and there are Nash maps $\gamma_1,...,\gamma_{n_3}\in 
N(M)^{n_2}$ which are generators of the kernel of the surjective $\Cal N$-homomorphism $p:\Cal N^{n_
2}\supset\Cal N_x^{n_2}\ni(\phi_1,...,\phi_{n_2})\to\sum_{i=1}^{n_2}\phi_i\beta_{ix}\in\Cal M_{2x}
\subset\Cal M_2,\ x\in M$. 
Let $\overline\alpha_1,...,\overline\alpha_{n_1}$ denote the images of $\alpha_1,...,\alpha_{n_1}$ 
in $H^0(M,\Cal M_2)$ under the homomorphism $h_*:H^0(M,\Cal
M_1)\to H^0(M,\Cal M_2)$ induced by $h$. 
\par 
We prove the theorem by induction on $n_2$. For $n_2=1$, there exist $\hat\alpha_1,...,\hat\alpha_{n_1}\in H^0(M,\Cal N)$ such that $p_*(\hat\alpha_i)=
\overline\alpha_i$, for $ i=1,...,n_1$ because the application $p_*:H^0(M,\Cal N)\to H^0(M,\Cal M
_2)$ is surjective by theorem 2.8 for $\Cal M_1=\Cal N$ ([C-R-S$_1$] and [C-S$_3$] ). 
Let $\delta_1,...,\delta_{n_4}\in N(M)^{n_1}$ be generators of the kernel of the surjective homomorphism 
$\Cal N^{n_1}\supset\Cal N_x^{n_1}\ni(\phi_1,...,\phi_{n_1})\to\sum_{i=1}^{n_1}\phi_i\alpha_{ix}\in
\Cal M_{1x}\subset\Cal M_x,\ x\in M$ 
(we choose the above $\{U_i\}$ so that $\delta_1,...,\delta_{n_4}$ exist). 
Multiplying $\alpha_i,\overline\alpha_i,\hat\alpha_i,\ \gamma_i$ 
and $\delta_i$ by a small positive Nash function, we can assume by the \L ojasiewicz inequality that the Nash maps $\hat\alpha_i,\ \gamma_i$ 
and $\delta_i$ are bounded. 
Then by Proposition VI.2.8 in [S$_2$] we can regard $M$ as the interior of a compact 
Nash manifold possibly with corners $\tilde M$ and the maps as the restrictions to $M$ of Nash maps 
$\tilde{\hat\alpha}_i,\ \tilde\gamma_i$ and $\tilde\delta_i$ on $\tilde M$. 
Replace $\Cal M_1$ and $\Cal M_2$ by the sheaves of $\Cal N$-modules
on $\tilde M$ given by $\Cal N^{n_1}/(
\tilde\delta_1,...,\tilde\delta_{n_4})\Cal N^{n_1}$ and $\Cal N/(\tilde\gamma_1,...,\tilde\gamma_{n_
3})\Cal N$, respectively, and replace $h:\Cal M_1\to\Cal M_2$ with the $\Cal N$-homomorphism $\tilde h:
\tilde M_1\to\tilde M_2$ defined by 
$$
\tilde h(0,...,0,\overset i\to1,0,...,0)=\tilde{\hat\alpha}_i\quad\text{mod}\ (\tilde\gamma_1,...,
\tilde\gamma_{n_3})\Cal N,\ i=1,...,n_1. 
$$
Then it suffices to see that $\Ker\tilde h$ is finite. 
Hence we assume from the beginning that $M$ is a compact Nash manifold possibly with corners. 
Then $\Ker h$ is isomorphic to $\Cal N\otimes_{N(M)}\Ker h_*$ by Theorem 5.2 in [C-R-S$_1$]. 
Hence $\Ker h$ is finite. \par
Let $n_2>1$ and assume that the theorem holds for $n_2-1$. 
Let $\Cal M_0$ denote the sheaf of $\Cal N$-ideals with $\Cal M_{0x}=\{0\}$. 
Set $\Cal M_3=\Cal M_2/p(\Cal N\times\Cal M_0\times\cdots\times\Cal M_0)$ and let $h_3:\Cal M_1\to
\Cal M_3$ denote the composite of $h$ with the projection from $\Cal
M_2$ to $\Cal M_3$. 
Then $\Cal M_3$ is generated by the images $\overline\beta_2,...,\overline\beta_{n_2}$ of $\beta_2,...,
\beta_{n_2}$, and $\gamma'_1,...,\gamma'_{n_3}\in N(M)^{n_2-1}$ are generators of the kernel of 
the $\Cal N$-homomorphism 
$$\Cal N^{n_2-1}\supset\Cal N_x^{n_2-1}\ni(\phi_1,...,\phi_{n_2-1})\to
\sum_{i=1}^{n_2-1}\phi_i\overline\beta_{i+1x}\in\Cal M_{3x}\subset\Cal
M_3,\ x\in M$$ 
where $\gamma
_i=(\gamma_{i,1},...,\gamma_{i,n_2})=(\gamma_{i,1},\gamma'_i)$, for $i=1,...,n_3$. 
Hence $\Cal M_3$ is finite, and by induction hypothesis $\Ker h_3$ is finite. 
Consider $h|_{\Ker h_3}:\Ker h_3\to\Cal M_2$. 
The image is contained in $p(\Cal N\times\Cal M_0\times\cdots\times\Cal M_0)$ which is isomorphic to 
$(\Ker p\cup\Cal N\times\Cal M_0\times\cdots\times\Cal M_0)/\Ker p$ and then to $\Cal N\times\Cal M_0
\times\cdots\times\Cal M_0/(\Ker p\cap\Cal N\times\Cal M_0\times\cdots\Cal M_0)$. 
Hence we can regard $h|_{\Ker h_3}$ as an $\Cal N$-homomorphism from $\Ker h_3$ to $\Cal N\times\Cal M_0
\times\cdots\times\Cal M_0/(\Ker p\cap\Cal N\times\Cal M_0\times\cdots\Cal M_0)$. 
In order to achieve the proof, we need to prove that $\Ker p\cap\Cal
N\times\Cal M_0\times\cdots\Cal M_0$ is finite. Define a sheaf of $\Cal N$-submodules $\Cal M$ of $\Cal N^{n_3}$ on $M$ by 
$$
\Cal M_x=\{(\phi_1,...,\phi_{n_3})\in\Cal N_x^{n_3}:\sum_{i=1}^{n_3}\phi_i\gamma_{i,jx}=0,\ j=2,...,n_
2\}. 
$$
Then it suffices to see that $\Cal M$ is finite because $\Ker p\cap\Cal N\times\Cal M_0\times\cdots\times
\Cal M_0$ is the image of $\Cal M$ under the $\Cal N$-homomorphism $:\Cal N^{n_3}\supset\Cal N_x^{n_3}
\ni(\phi_1,...,\phi_{n_3})\to(\sum_{i=1}^{n_3}\phi_i\gamma_{i,1x},0,...,\allowmathbreak0)\in\Cal N_x
\times\{0\}\times\cdots\times\{0\}\subset\Cal N\times\Cal M_0\times\cdots\times\Cal M_0,\ x\in M$. 
On the other hand, if we define an $\Cal N$-homomorphism $r:\Cal N^{n_3}\to\Cal N^{n_2-1}$ by $r(\phi
_1,...,\phi_{n_3})=(\sum_{i=1}^{n_3}\phi_i\gamma_{i,2x},...,\allowmathbreak\sum_{i=1}^{n_3}\phi_i
\gamma_{i,n_2x})$ for $(\phi_1,...,\phi_{n_3})\in\Cal N_x^{n_3},\ x\in M$, then $\Ker r=\Cal M$. 
As in the case of $n_2=1$ we reduce the problem to the case where $\gamma_{i,j}$ are bounded and then 
$M$ is a compact Nash manifold possibly with corners. 
Then $\Ker r$ is finite by Theorem 5.2 in [C-R-S$_1$]. \par
Thus $\Ker p\cap\Cal N\times\Cal M_0\times\cdots\times\Cal M_0$ is finite. 
We can regard it as a sheaf of $\Cal N$-ideals. 
Hence by the result in case of $n_2=1$, $\Ker(h|_{\Ker h_3})=\Ker h$ is finite. 
\qed
\enddemo
The following two theorems do not hold for general sheaves of $\Cal N$-modules, [Hu], [B-C-R] and VI.2.10 
in [S$_2$]. 
However, our case is sufficient for the applications we have in mind
in this paper. 

\proclaim{Theorem 2.7}
(Nash case of Cartan Theorem A) 
Let $\Cal M$ be a finite sheaf of $\Cal N$-submodules of $\Cal N^n$ on a Nash manifold $M$ for $n>0\in
\N$. 
Then $\Cal M$ is finitely generated by its global cross-sections. 
\endproclaim
\demo{Proof}
We assume that $n>1$ and proceed by induction on $n$. 
Let $p:\Cal N^n\to\Cal N^{n-1}$ denote the projection forgetting the first factor, and set $\Cal M_1=
\Ker p|_\Cal M$ and $\Cal M_2=\Im p|_\Cal M$. 
Then the sequence $0\longrightarrow\Cal M_1\overset p_1\to\longrightarrow\Cal M\overset p_2\to
\longrightarrow\Cal M_2\longrightarrow0$ is exact, we can regard $\Cal M_1$ as a sheaf of 
$\Cal N$-ideals, which is finite by theorem 2.6, and $\Cal M_2$ is clearly a finite sheaf of 
$\Cal N$-submodules of $\Cal N^{n-1}$. 
By induction hypothesis we have global generators $h_1,...,h_l$ of $\Cal M_1$ and $g_1,...,g_k$ of 
$\Cal M_2$. 
Then it suffices to find $f_1,...,f_k\in H^0(M,\Cal M)$ such that $p_{2*}(f_i)=g_i,\ i=1,...,k$ 
because $f_1,...,f_k,h_1,...,h_l$ are generators of $\Cal M$. \par
Fix $i$. 
Since $H^0(M,\Cal M)\subset N(M)^n$ and $H^0(M,\Cal M_2)\subset N(M)^{n-1}$, setting $g_i=(g_{i,2},
...,g_{i,n})$ we construct $g_{i,1}\in N(M)$ such that $(g_{i,1},...,g_{i,n})\in H^0(M,\Cal M)$. 
For each $x\in M$, the set $\Phi_x=\{\phi\in\Cal N_x:(\phi,g_{ix})\in\Cal M_x\}$ is a residue class of 
$\Cal N_x$ modulo $\Cal M_{1x}$, and the correspondence $\Phi:x\to\Phi_x$ is a global cross-section 
of $\Cal N/\Cal M_1$. Actually, it suffices to check it on each member of a finite open semialgebraic covering of $M$, we assume that 
$\Cal M$ is generated by global cross-sections $\alpha_1=(\alpha_{1,1},...,\allowmathbreak\alpha_{1,n}),...,\alpha_
{k'}=(\alpha_{k',1},...,\alpha_{k',n})\in N(M)^n$. 
Then $\alpha'_1=(\alpha_{1,2},...,\alpha_{1,n}),...,\alpha'_{k'}=(\alpha_{k',2},...,\alpha_{k',n})$ 
are also generators of $\Cal M_2$. 
Let $\Cal M_3$ denote the kernel of the $\Cal N$-homomorphism $\Cal N^{k'+1}\supset\Cal N_x^{k'+1}\ni(
\phi_1,...,\phi_{k'+1})\to\sum_{j=1}^{k'}\phi_j\alpha'_{jx}-\phi_{k'+1}g_{ix}\in\Cal N_x^{n-1}\subset
\Cal N^{n-1},\ x\in M$. 
Then $\Cal M_3$ is finite by theorem 2.6, and each stalk $\Cal M_{3x}$ contains a germ of the form 
$(\phi_1,...,\phi_{k'},1)$. 
Hence refining the covering if necessary, we assume that $\Cal M_3$ is generated by a finite number of 
global cross-sections. 
Then we have $\beta_1,...,\beta_{k'}\in N(M)$ such that $g_i=\sum_{j=1}^{k'}\beta_j\alpha'_j$. 
It follows $\Phi=\sum_{j=1}^{k'}\beta_j\alpha_{j,1}$ mod $\Cal M_1$. 
Thus $\Phi$ is a global cross-section. \par 
Apply the next theorem to the projection $\Cal N\to\Cal N/\Cal M_1$ and $\Phi$. 
Then there exists $g_{i,1}\in N(M)$ such that $g_{i,1x}=\Phi_x$ mod $\Cal M_{1x}$ for $x\in M$ and 
hence $(g_{i,1},...,g_{i,n})\in H^0(M,\Cal M)$. 
\qed
\enddemo
\proclaim{Theorem 2.8}
(Nash case of Cartan Theorem B) 
Let $h:\Cal M_1\to\Cal M_2$ be a surjective $\Cal N$-homomorphism between finite sheaves of $\Cal N
$-modules on a Nash manifold $M$. 
Assume that $\Cal M_1$ is finitely generated by its global cross-sections. 
Then the induced map $h_*:H^0(M,\Cal M_1)\to H^0(M,\Cal M_2)$ is surjective. 
\endproclaim
\demo{Proof}
We can assume that $\Cal M_1=\Cal N^n$ for some $n>0\in\N$ because there exist global generators $g_1,...,
g_n$ of $\Cal M_1$ and then we can replace $h$ with the surjective homomorphism $\Cal N^n\supset\Cal N
_x^n\ni(\phi_1,...,\phi_n)\to h(\sum_{i=1}^n\phi_ig_{ix})\in\Cal M_{2x}\subset\Cal M_2,\ x\in M$. 
Set $\Cal M=\Ker h$. 
Then by theorem 2.6, $\Cal M$ is a finite sheaf of $\Cal N$-submodules of $\Cal N^n$, and $h:\Cal N^n
\to\Cal M_2$ coincides with the projection $p:\Cal N^n\to \Cal N^n/\Cal M$. 
Hence we consider $p$ in place of $h$. 
Assume that $n>1$ and the theorem holds for smaller $n$. \par
Let $f\in H^0(M,\Cal N^n/\Cal M)$. 
We need to find $g\in H^0(M,\Cal N^n)=N(M)^n$ such that $p_*(g)=f$. 
Let $\Cal M_0$ denote the sheaf of $\Cal N$-ideals with $\Cal M_{0x}=\{0\}$ for $x\in M$. 
Then the homomorphism $\Cal M_0\times\Cal N^{n-1}\to\Cal N^n/(\Cal M+\Cal N\times\Cal M_0\times\cdots\times
\Cal M_0)$ is surjective and we can regard it as the projection $\Cal N^{n-1}\to\Cal N^{n-1}/\Cal L$ for 
some finite sheaf of $\Cal N$-submodules $\Cal L$ of $\Cal N^{n-1}$. 
Hence by induction hypothesis there exists $(0,g_2,...,g_n)\in H^0(M,\Cal M_0\times\Cal N^{n-1})$ 
whose image in $H^0(M,\Cal N^n/(\Cal M+\Cal N\times\Cal M_0\times\cdots\times\Cal M_0))$ coincides 
with the image of $f$ there. 
Replace $f$ with the difference of $f$ and the image of $(0,g_2,...,g_n)$ in $H^0(M,\Cal N^n/\Cal M)$. 
Then we can assume from the beginning that $f\in H^0(M,(\Cal M+\Cal N\times\Cal M_0\times\cdots\times
\Cal M_0)/\Cal M)$. 
Hence we regard $f$ as a global cross-section of $\Cal N\times\Cal M_0\times\cdots\times\Cal M_0/(
\Cal M\cap\Cal N\times\Cal M_0\times\cdots\times\Cal M_0)$ since $(\Cal M+\Cal N\times\Cal M_0\times
\cdots\times\Cal M_0)/\Cal M$ is naturally isomorphic to $\Cal N\times\Cal M_0\times\cdots\times
\Cal M_0/(\Cal M\cap\Cal N\times\Cal M_0\times\cdots\times\Cal M_0)$. 
It was shown in the proof of theorem 2.6 that $\Cal M\cap\Cal N\times\Cal M_0\times\cdots\times\Cal M_0$ 
is finite. 
Hence $f$ is the image of some global cross-section $g$ of $\Cal N\times\Cal M_0\times\cdots\times
\Cal M_0$ under the projection $H^0(M,\Cal N\times\Cal M_0\times\cdots\times\Cal M_0)\to H^0(M,
\Cal N\times\Cal M_0\times\cdots\times\Cal M_0/(\Cal M\cap\Cal N\times\Cal M_0\times\cdots\times
\Cal M_0))$ because this is the case of $\Cal M_1=\Cal N$ in the theorem. 
Then $p_*(g)=f$. 
\qed
\enddemo
Let $X$ be a Nash subset of $\R^n$ and $f_1,...,f_k$ be generators of the ideal of $N(\R^n)$ of 
functions vanishing on $X$. 
Let $\Sing X$ denote the subset of $X$ where the Jacobian matrix rank of $f_1,...,f_k$ is smaller than 
$\codim X$. 
Let a complexification $X^\C$ of $X$ in $\C^n$ be defined to be the common zero set of some 
complexifications $f_1^\C,...,f_k^\C$ of $f_1,...,f_k$. 
Then by Lemma 1.9 and Theorem 1.10 in [C-R-S$_2$] and theorem 2.7, we obtain the next remark. 
\remark{Remark}
$\Sing X$ is the smallest Nash subset of $X$ whose complement is a
Nash manifold. But it does not coincide in general with points in $X$
where the germ of $X$ is not a Nash manifold germ of $\dim
X$. Moreover $\Sing X$ is also equal to
$X\cap\Sing X^\C$, where $\Sing X^\C$ denotes the $C^\omega$ singular point set of $X^\C$. 
\endremark
\medskip
We deduce from [Hi] a Nash version of Hironaka Desingularization
Theorem that will be useful in our context. 
\proclaim{Theorem 2.9}
(Nash case of Main Theorem I of [Hi]) 
Let $X$ be a Nash subset of $\R^n$. 
Then there exists a finite sequence of blowings-up $X_r\overset\pi_r\to\longrightarrow\cdots\overset\pi_1
\to\longrightarrow X_0=X$ along smooth Nash centers $C_i\subset X_i,\ i=1,...,r-1$, such that $X_r$ is 
smooth and $C_i\subset\Sing X_i$. 
\endproclaim

\demo{Proof}
Since $N(\R^n)$ is a Noetherian ring ([E] and [Ri]), we have generators $f_1,...,f_k$ of the ideal of $N(\R^n)$ of functions vanishing on $X$. 
Set $F=(f_1,...,f_k)$, which is a Nash map from $\R^n$ to $\R^k$, and $Y=\graph F$. 
Let $Y^Z$ denote the Zariski closure of $Y$ in $\R^n\times\R^k$ and let $\widetilde{Y^Z}\subset\R^n\times
\R^k\times\R^{n'}$ be an algebraic set such that the restriction $p$
to $\widetilde{Y^Z}$ of the projection
$\R^n\times\R^k\times\R^{n'}\to\R^n\times\R^k$ is the 
normalization of $Y^Z$ 
(we simply call $\widetilde{Y^Z}$ the {\it normalization} of $Y^Z$). 
Then by Artin-Mazur Theorem (see Theorem I.5.1 in [S$_2$]) there exists a connected component $L$ of 
$\widetilde{Y^Z}$ consisting of only regular points such that $p(L)=Y$ and $p|_L:L\to Y$ is a Nash 
diffeomorphism. 
Let $q_1:\widetilde{Y^Z}\to\R^n$ and $q_2:\widetilde{Y^Z}\to\R^k$ denote the restrictions to $\widetilde
{Y^Z}$ of the projections $\R^n\times\R^k\times\R^{n'}\to\R^n$ and $\R^n\times\R^k\times\R^{n'}\to\R^k$, 
respectively. 
Then $q_1|_L$ is a Nash diffeomorphism onto $\R^n$, the set
$q_2^{-1}(0)$ is algebraic, the equality $(q_1|_L)^{-1}(X)
=(q_2|_L)^{-1}(0)$ holds, and $(q_1|_L)^{-1}(\Sing X)$ is equal to the intersection of $L$ with the algebraic singular 
point set of $q_2^{-1}(0)$. 
Indeed, $(q_1|_L)^{-1}(\Sing X)$ is contained in the above intersection because
$(q_1|_L)
^{-1}(\Sing X)$ is the smallest Nash subset of $(q_1|_L)^{-1}(X)\,(=(q_2|_L)^{-1}(0))$ whose complement 
is a Nash manifold (by the remark before theorem 2.9), and the converse inclusion follows from the equality $q_2=F\circ q_1$ on $L$. 
Hence we can replace $X$ by $L\cap q_2^{-1}(0)$---the union of some connected components of $q_2^{-1}(0)$. 
By Main Theorem I there exists a finite sequence of blowings-up $\tilde X_r\overset\tilde\pi_r\to
\longrightarrow\cdots\overset\tilde\pi_1\to\longrightarrow\tilde X_0=q_2^{-1}(0)$ along smooth algebraic 
centers $\tilde C_i\subset\tilde X_i$, for $i=0,...,r-1$, such that $\tilde X_r$ is smooth and $\tilde C_i
\subset\Sing\tilde X_i$. 
Then $\tilde X_r\cap(\tilde\pi_1\circ\cdots\circ\tilde\pi_r)^{-1}(L)\to\cdots\to\tilde X_0\cap L$ 
fulfills the requirements. 
\qed\enddemo
A sheaf of $\Cal N$-($\Cal O$-)ideals on a Nash (analytic) manifold $M$ is called {\it normal crossing} 
if there exists a local Nash (analytic) coordinate system $(x_1,...,x_n)$ of $M$ at each point such that 
the stalk of the sheaf is generated by $\prod_{i=1}^nx_i^{\alpha_i}$ for some $(\alpha_1,...,\alpha_n)\in
\N^n$. 

\proclaim{Theorem 2.10}
(Nash case of Main Theorem II of [Hi]) 
Let $M$ be a Nash manifold and let $\Cal I_1$ and $\Cal I_2$ be finite sheaves of non-zero $\Cal N$-ideals 
on $M$. 
Assume that $\Cal I_2$ is normal crossing. 
Then there exists a finite sequence of blowings-up $M_r\overset\pi_r\to\longrightarrow\cdots\overset\pi_1
\to\longrightarrow M_0=M$ along smooth Nash centers $C_i\subset M_i$,
for $i=1,...,r-1$, such that $(\pi_1\circ
\cdots\circ\pi_r)^{-1}\Cal I_1\Cal I_2\Cal N_{M_r}$ is normal crossing, each $C_i$ is normal crossing with 
$(\pi_1\circ\cdots\circ\pi_i)^{-1}(\supp\Cal N_M/\Cal I_2)\cup\cup_{j=1}^i(\pi_j\circ\cdots\circ\pi_i)^{-1
}(C_{j-1})$ and $\pi_1\circ\cdots\circ\pi_i(C_i)$ is contained in the
subset of $M$ consisting of $x$ such that even $\Cal I_{1x}$ is not generated 
by any power of one regular function germ or $\Cal I_{1x}+\Cal I_{2x}\not=\Cal N_x$. 
\endproclaim
Note that $(\pi_1\circ\cdots\circ\pi_i)^{-1}\Cal I_2\Cal N_{M_i}$, for
$i=1,...,r$, are normal crossing. 
\demo{Proof}
Let $f_1,...,f_{k'}\in N(\R^n)$ and $f_{k'+1},...,f_k\in N(\R^n)$ be global generators of $\Cal I_1$ and 
$\Cal I_2$ (theorem 2.7), respectively, and define $F,Y,Y^Z,\widetilde{Y^Z},L,q_1:\widetilde{Y^Z}\to\R^n$ 
and $q_2:\widetilde{Y^Z}\to\R^k$ as in the last proof. 
Let $W$ be the subset of $\widetilde{Y^Z}$ consisting of points where $f_{k'+1}\circ q_1,...,f_k\circ q_1$ do not 
generate a normal crossing sheaf of $\Cal N$-ideals. 
Consider the algebraic $\R$-scheme of the topological underlying space $\widetilde{Y^Z}-\Sing\widetilde{Y
^Z}-W$, and let $\Cal J_1$ and $\Cal J_2$ denote the sheaf of ideals of the scheme generated by $f_1
\circ q_1,...,f_{k'}\circ q_1$ and by $f_{k'+1}\circ q_1,...,f_k\circ q_1$, respectively. 
Then we can replace $M,\Cal I_1$ and $\Cal I_2$ with the scheme, $\Cal J_1$ and $\Cal J_2$. 
Hence the theorem follows from Main Theorem II. 
\qed\enddemo

\remark{Remark}
Note that main Theorems I and II of [Hi] state some additional
conditions that are automatically satisfied in the Nash case.
\endremark
%%%%%%%%%%%%%%%%%%%%%%%%%%%%%%%%%%%%%%%%%%%%%%%%%%%%%%%%%%%%%%%%%%%%%%%%%%%%%%%%%%%%%%%%%%
%%%%%%%%%%%%%%%%%%%%%%%%%%%%%%%%%%%%%%%%%%%%%%%%%%%%%%%%%%%%%%%%%%%%%%%%%%%%%%%%%%%%%%%%%%
\subhead
2.3. Topologies on function spaces
\endsubhead

Let $M$ be a $C^\infty$ manifold. 
We use three kinds of topologies on $C^\infty(M)$ as a topological
linear space.

The first is the classical compact-open $C^r$ topology,
$r=0,...,\infty$, for which $C^\infty(M)$ is a Fr\'echet space if $r=\infty$.

The second is the {\it Whitney $C^r$ topology}, $r=0,...,\infty$. Even
if it is well-known, we
recall its definition because we will define the third topology below by comparison with it. 
If $M=\R^n$, then a system of open neighborhoods of the zero function
in $C^\infty(M)$ is given by
$$
U_{r',g_\alpha}=\{f\in C^\infty(\R^n):|D^\alpha f(x)|<g_\alpha(x),\,\alpha\in\N^n,\,|\alpha|\le r\}
$$
where $r'$ runs in $\{m\in\N:m\le r\}$ and $g_\alpha$ runs in $C^\infty(\R^n)$ with $g_\alpha>0$ everywhere for each $\alpha\in\N^n$ with $|\alpha|\le r'$. 
If $M$ is an open subset of $\R^n$, we define the topology on $C^\infty(M)$ in the same way. 
In general, embed $M$ in some $\R^n$ and let $p:V\to M$ be the orthogonal projection of a tubular neighborhood of $M$ in $\R^n$. 
Then $p$ induces an injective linear map $C^\infty(M)\ni f\to f\circ
p\in C^\infty(V)$ whose image is closed in $C^\infty(V)$ in the
Whitney $C^r$ topology. Hence $C^\infty(M)$ inherits a topology as a
closed subspace of $C^\infty(V)$. We call it the Whitney $C^\infty$ topology. \par
  The {\it strong Whitney $C^\infty$ topology} is the third topology which we will consider. 
Assume first that $M=\R^n$, and let $g_\alpha$ be a
positive-valued $C^\infty$ function on $\R^n$ and $K_\alpha$ be a
compact subset of $\R^n$ for each $\alpha\in\N^n$ such that
$\{\R^n-K_\alpha\}$ is locally finite. Set $g=(g_\alpha)_\alpha$ and $K=(K_\alpha)_\alpha$. 
Then a system of open neighborhoods of the zero function in $C^\infty(\R^n)$ is given by the family of sets 
$$
U_{g,\alpha}=\{f\in C^\infty(\R^n):|D^\alpha f(x)|<g_\alpha(x)\ \text{for}\ x\in\R^n-K_\alpha\ \text{for}\ \alpha\in\N^n\}
$$
for all $g$ and $K$. 
We define the strong Whitney $C^\infty$ topology on a general
manifold $M$ in the same way as in the case of the Whitney $C^\infty$ topology. 
Moreover, we shall need to consider $C^\infty$ functions on an
analytic set and the strong Whitney $C^\infty$ topology on the
space. To this aim, we use another equivalent definition of the topology. 
Let $\{M_l\}$ be a family of compact $C^\infty$ submanifolds of $M$ possibly with boundary such 
that $\{\Int M_l\}$ is a locally finite covering of $M$. 
Regard $C^\infty(M)$ as a subset of $\prod_l C^\infty(M_l)$ by the injective map $C^\infty(M)\ni 
f\to\prod_l f|_{M_l}\in\prod_l C^\infty(M_l)$. 
Then the family of sets $C^\infty(M)\cap\prod_l O_l$ is the system of open sets of $C^\infty(M)$, where $ O_l$ are open subsets of $C^\infty(M_l)$ in the $C^\infty$ topology. 
Note that the product topology of $\prod_l C^\infty(M_l)$ induces the compact-open $C^\infty$ topology. \par
For an analytic manifold $M$, we endow $C^\omega(M)$ with the three topologies in the same way, and we extend naturally the definition of the topologies to the spaces of $C^\infty$ or $C^\omega$ maps between $C^\infty$ or $C^\omega$ manifolds. 
\remark{Remark 2.11} 
The first three remarks explain essential differences between the three topologies. \par
(1) The compact-open $C^\infty$ topology, the Whitney $C^\infty$ topology and the strong Whitney $C^\infty$ topology coincide if $M$ is compact. \par
(2) The strong Whitney $C^\infty$ topology is stronger than the Whitney $C^\infty$ topology if $M$ is not compact. \par
(3) $C^\infty(M)$ is not a Fr\'echet space in the Whitney $C^r$
topology nor the strong Whitney $C^\infty$ topology if $M$ is not
compact. Indeed, it is even not metrizable. \par
The following remarks will be useful in the sequel. \par
(4) Whitney Approximation Theorem---any $C^\infty$ function on an analytic manifold is approximated by a $C^\omega$ function---holds also in any of these topologies (see [W]). \par
Finally, an advantage of the strong Whitney $C^\infty$ topology is
that we can reduce many global problems to local problems using
partition of unity. \par
(5) Let $\{\phi_i\}$ be a partition of unity of class $C^\infty$ on $M$. 
Then for a neighborhood $U$ of 0 in $C^\infty(M)$ in the strong Whitney $C^\infty$ topology there exists another $V$ such that if $f\in V$ then $\phi_if\in U$ for all $i$ and conversely if $\phi_if\in V$ for all $i$ then $f\in U$. 
\endremark
\medskip
Let $M$ be a Nash manifold. 
We give a topology on $N(M)$, called the {\it semialgebraic $C^r$
topology}, $r=0,...,\infty$, so that a system of open neighborhoods of
0 in $N(M)$ is given by the family $U_{r',g_\alpha}$ defined in the
above definition of the Whitney $C^r$ topology, where $g_\alpha$ runs
here in $N(M)$ only. 
If $r=\infty$, we call it the {\it Nash topology}. 
For $r<\infty$, let $N^r(M)$ denote the space of semialgebraic $C^r$
functions on $M$. We define {\it semialgebraic $C^{r'}$ topology} on $N^r(M)$ for $r'\le r$ in the same way. 
We do not need the analog on $N(M)$ of the strong Whitney $C^\infty$
topology. 
When $M$ is not compact, it is the discrete topology by Proposition VI.2.8, [S$_2$]
and next remark. 
A {\it partition of unity of class semialgebraic} $C^r,\ r\in\N$, on
$M$ is a {\bf finite} family of non-negative semialgebraic $C^r$ functions on $M$ whose sum equals 1. 
\remark{Remark 2.11,(5)$'$} 
Let $r'\le r\in\N$, and let $\{\phi_i\}$ be a partition of unity of class semialgebraic $C^{r'}$ on $M$. 
Then for a neighborhood $U$ of 0 in $N^r(M)$ in the $C^{r'}$ topology there exists a neighborhood $V$ of 0 in $N^r(M)$ such that if $f\in V$ then $\phi_if\in U$ for all $i$ and conversely if $\phi_if\in V$ for all $i$ then $f\in U$. \par
The reason is that $\{\phi_i\}$ is a finite family and the map $N^r(M)\ni f\to\phi_if\in N^r(M)$ is continuous for each $i$ by lemma II.1.6, [S$_2$], which states that $N^r(M)$ and $N(M)$ are topological rings in the semialgebraic $C^{r'}$ topology. 

\endremark
\medskip
We need also the following lemma many times. 
\proclaim{Lemma 2.12} Let $M$ be an analytic manifold and $\xi_1,...,\xi_l$ be analytic functions on $M$. 
Then the maps $\Xi^\infty:C^\infty(M)^l\ni(h_1,...,h_l)\to\sum_{i=1}^l\xi_ih_i\in\sum_{i=1}^l\xi_iC^\infty(M)$ and $\Xi^\omega:C^\omega(M)^l\ni(h_1,...,h_l)\to\sum_{i=1}^
l\xi_ih_i\in\sum_{i=1}^l\xi_iC^\omega(M)$ are open in both the compact-open $C^\infty$ topology and the strong Whitney $C^\infty$ topology. 
\endproclaim
Note that in the case $l=1$, the lemma is much easier to prove because
the involved maps are injective. Moreover, the lemma does not necessarily hold in the Whitney $C^\infty$ topology. 
This is one reason why we need to have recourse to the strong Whitney
$C^\infty$ topology in the paper. 
\demo{Proof}
Consider $\Xi^\infty$ in the compact-open $C^\infty$ topology. 
It is well-known that the ideal of $C^\infty(M)$ generated by a finite
number of analytic functions is closed in $C^\infty(M)$ in any of the
$C^\infty$ topologies (which follows from Theorems III.4.9 and VI.1.1$'$, [Ml]). 
In particular $\sum_{i=1}^l\xi_iC^\infty(M)$ is a Fr\'echet space in
the compact-open $C^\infty$ topology, and $\Xi^\infty$ is open by the open mapping theorem on Fr\'echet spaces. 
Note that the above proof is still valid in the case of an analytic manifold with corners. \par
  Consider $\Xi^\infty$ in the strong Whitney $C^\infty$ topology. 
Let $M_j$ be compact $C^\omega$ submanifolds of $M$ with boundary such that $\{\Int M_j\}$ is a locally finite covering of $M$. 
Let $\{\phi_j\}$ be a partition of unity of class $C^\infty$ subordinate to $\{\Int M_j\}$. 
As shown above, the map $C^\infty(M_j)^l\ni(h_1,...,h_l)\to\sum_{i=1}^l\xi_i|_{M_j}h_i\in\sum_{i=1}^l\xi_i|_{M_j} C^\infty(M_j)$ is open for each $j$. 
Hence for each $h=(h_1,...,h_l)\in C^\infty(M)^l$ and $g\in\sum_{i=1}^l\xi_iC^\infty(M)$ sufficiently close to $\sum_{i=1}^l\xi_ih_i$ in the strong Whitney $C^\infty$ topology there exist $g_j=(g_{1,j},...,g_{l,j})\in C^\infty(M_j)^l$ close to $h|_{M_j}$ such that $\sum_{i=1}^l\xi_i|_{M_j}g_{i,j}=g|_{M_j}$ for any $j$. 
Then $\sum_j\phi_jg_j$ is a well-defined $C^\infty$ map from $M$ to $\R^l$ and close to $h$ by remark 2.11,(5), and $\sum_{i=1}^l\xi_i\sum_j\phi_jg_{i,j}=\sum_j\phi_jg=g$. 
Thus $\Xi^\infty$ is also open in the strong Whitney $C^\infty$ topology. \par
  We finally consider $\Xi^\omega$ only in the strong Whitney
$C^\infty$ topology (the proof is similar, and even easier, in the compact-open $C^\infty$ topology). 
Let $(h_1,...,h_l)\in C^\omega(M)^l$ such that $\sum_{i=1}^l\xi_ih_i$ is small. 
Then, by openness of $\Xi^\infty$, there exists small $(h'_1,...,h'_l)\in C^\infty(M)^l$ such that 
$\sum_{i=1}^l\xi_ih'_i=\sum_{i=1}^l\xi_ih_i$ and hence $(h_1-h'_1,...,h_l-h'_l)\in\Ker\Xi^\infty$. 
Therefore, it suffices to see that $\Ker\Xi^\omega$ is dense in $\Ker\Xi^\infty$. 
Let $H=(h_1,...,h_l)\in\Ker\Xi^\infty$. 
We want to approximate $H$ by an element of $\Ker\Xi^\omega$. \par
Let $\Cal J$ denote the kernel of the homomorphism $:\Cal O^l\supset\Cal O^l_a\ni(\phi_1,...,
\phi_l)\to\sum_{i=1}^l\xi_{ia}\phi_i\in\Cal O_a\subset\Cal O,\ a\in M$, which is a coherent sheaf 
of $\Cal O$-submodules of $\Cal O^l$ by theorem 2.5. 
Let $M^\C$ and $\Cal J^\C$ be Stein and coherent complexifications of $M$ and $\Cal J$ which are 
complex conjugation preserving. 
Let $\{U_i\}$ be a locally finite open covering of $M^\C$ such that each $\overline U_i$ is 
compact. Let $H_{i,j}=(h_{1,i,j},...,h_{l,i,j})$, for $j=1,...,n_i,\,
i=1,2,...$, be global 
cross-sections of $\Cal J$ such that $H_{i,j}|_M$ are real valued and $H_{i,1},...,
H_{i,n_i}$ generate $\Cal J^\C$ on $U_i$ (theorem 2.1) for each $i$. 
Then $H_{i,1}|_{U_i\cap M},...,H_{i,n_i}|_{U_i\cap M}$ are generators
of $\Ker\Xi^\infty|_{U_i\cap M}$. Actually, by Theorem VI,1.1$'$ in [Ml] it is equivalent to prove that $\Ker\Xi_a=\Cal F_a\Ker\Xi^\omega_a,\ a\in U_i$, 
where $\Cal F_a$ is the completion of $\Cal O_a$ in the $\frak p$-adic 
topology and the homomorphisms $\Xi^\omega_a:\Cal O^l_a\to\Cal O_a$ and $\Xi_a:\Cal F^l_a\to
\Cal F_a$ are naturally defined. 
However, this condition is the flatness of $\Cal F_a$ over $\Cal O_a$, 
which is well-known (see [Ml]). 
Thus $H_{i,1}|_{U_i\cap M},...,H_{i,n_i}|_{U_i\cap M}$ generate $\Ker\Xi^\infty|_{U_i\cap M}$. 
Let $\{\rho_i\}$ be a partition of unity of class $C^\infty$ subordinate to $\{U_i\cap M\}$. 
Then $\rho_iH\in\Ker\Xi^\infty|_{U_i\cap M}$ and we have $C^\infty$
functions $\chi_{i,j}$ on $M$, for $ j=1,...,n_i,\ i=1,2,...$, such that $\supp\chi_{i,j}\subset U_i\cap M$ and $\rho_iH=\sum_{j=1}^{n
_i}\chi_{i,j}H_{i,j}|_M$. 
By remark 2.11,(4) we can approximate $\chi_{i,j}$ by analytic functions $\chi'_{i,j}$. 
Moreover, as in [W] we can approximate so that each $\chi'
_{i,j}$ can be complexified to a complex analytic function $\chi^{\prime\C}_{i,j}$ on $M^\C$ and 
$\sum_{i,j}|\chi^{\prime\C}_{i,j}H_{i,j}|$ is locally uniformly bounded. 
Then $\sum_{i,j}\chi^{\prime\C}_{i,j}H_{i,j}$ is a complex analytic map from $M^\C$ to $\C^l$, 
and its restriction to $M$ is both an approximation of $H$ and an element of $\Ker\Xi^\omega$. 
Thus $\Xi^\omega$ is open. 
\qed
\enddemo
Lemma 2.12 holds in the Nash case for a compact Nash manifold. 
However, we do not know whether lemma 2.12 still holds for a
non-compact Nash manifold. Consequently, we have recourse many times in this paper to compactification arguments that require much care to deal with. 
%%%%%%%%%%%%%%%%%%%%%%%%%%%%%%%%%%%%%%%%%%%%%%%%%%%%%%%%%%%%%%%%%%%%%%%%%%%%%%%%%%%%%%%%%%%%%%%%%%%
%%%%%%%%%%%%%%%%%%%%%%%%%%%%%%%%%%%%%%%%%%%%%%%%%%%%%%%%%%%%%%%%%%%%%%%%%%%%%%%%%%%%%%%%%%%%%%%%%%%
%%%%%%%%%%%%%%%%%%%%%%%%%%%%%%%%%%%%%%%%%%%%%%%%%%%%%%%%%%%%%%%%%%%%%%%%%%%%%%%%%%%%%%%%%%%%%%%%%%%

\head
3. Equivalence of normal crossing functions
\endhead

%%%%%%%%%%%%%%%%%%%%%%%%%%%%%%%%%%%%%%%%%%%%%%%%%%%%%%%%%%%%%%%%%%%%%%%%%%%%%%%%%%%%%%%%%%%%%%%%%%%
%%%%%%%%%%%%%%%%%%%%%%%%%%%%%%%%%%%%%%%%%%%%%%%%%%%%%%%%%%%%%%%%%%%%%%%%%%%%%%%%%%%%%%%%%%%%%%%%%%%

\subhead
3.1. On $C^{\infty}$ equivalence of analytic functions with only normal crossing singularities
\endsubhead

Let us compare $C^\omega$ and $C^\infty$ right equivalences of two analytic functions on an analytic 
manifold. The $C^\infty$ right equivalence is easier to check. 
The
$C^\omega$ right equivalence implies the latter. 
However the converse is not necessarily true. 
We will show that this is the case for analytic functions with only normal crossing singularities, and apply 
the fact to the proof of the main theorem 3.2. \par

The main theorem of this section is

\proclaim{Theorem 3.1 }
(1) Let $M$ be a $C^\omega$ manifold and $f,g\in C^\omega(M)$. 
Assume that $f$ and $g$ admit only normal crossing singularities. 
If $f$ is $C^\infty$ right equivalent to $g$, then $f$ is $C^\omega$ right equivalent to $g$. \par
(2) If $C^\infty$ functions $f$ and $g$ on a $C^\infty$ manifold $M$ admit only normal crossing 
singularities and are proper and $C^2$ right equivalent, then they are $C^\infty$ right equivalent. 
\par
(3) If $f$ and $g$ are semialgebraically $C^2$ right equivalent Nash functions on a Nash manifold 
$M$ with only normal crossing singularities, then they are Nash right equivalent. \par
The case where $M$ has corners also holds. 
\endproclaim
\remark{Remark}
(i) The germ case is also of interest. 
Let $M,\,f$ and $g$ be the same as in above (1). 
Let $\phi$ be a $C^\infty$ diffeomorphism of $M$ such that $f=g\circ\phi$. 
Set $X=\Sing f$ and 
$Y=\Sing g$, and let $\{X_i\}$ and $\{Y_i\}$ be the irreducible analytic components of $X$ and $Y$, 
respectively. 
Let $A$ and $B$ be the unions of some intersections of some $X_i$ and $Y_i$, respectively. 
Assume that $\phi(A)=B$. 
Then we can choose a $C^\omega$ diffeomorphism $\pi$ so that $f=g\circ\pi$ and $\pi(A)=B$. 
Consequently, theorem 3.1,(1) holds for the germs of $f$ on $A$ and $g$ on $B$. 
Similar statements for (2) and (3) hold. \par
(ii) In the Nash case, $C^\infty$ right equivalence does not imply Nash right equivalence. Indeed, let $N$ be a compact contractible Nash manifold with non-simply connected boundary of 
dimension $n>3$ (e.g., see [Mz]). 
Set $M=(\Int N)\times(0,\,1)$ and let $f:M\to(0,\,1)$ denote the projection. 
Then $M$ and $f$ are of class Nash, and $M$ is Nash diffeomorphic to
$\R^{n+1}$. Actually, smooth the corners of $N\times[0,\,1]$. 
Then $N\times[0,\,1]$ is a compact contractible Nash manifold with simply connected boundary 
of dimension strictly more than four. 
Hence by the positive answers to Poincar\'e conjecture and Sch\"onflies problem (Brown-Mazur 
Theorem) $N\times[0,\,1]$ is $C^\infty$ diffeomorphic to an $(n+1)$-ball. 
Hence by Theorem VI.2.2 in [S$_2$] $M$ is Nash diffeomorphic to an open $(n+1)$-ball. 
Let $g:M\to\R$ be a Nash function which is Nash right equivalent to the projection $\R^n
\times(0,\,1)\to(0,\,1)$. 
Then $f$ and $g$ are $C^\omega$ right equivalent since $\Int N$ is $C^\omega$ diffeomorphic 
to $\R^n$, but they are not Nash equivalent because $\Int N$ and 
$\R^n$ are not Nash diffeomorphic, by Theorem VI 2.2 in [S$_2$]. 
\endremark
For the proof of part (2) and (3) of the theorem, we need to prepare
some material in next part. Therefore we postpone their proof to the last part of the paper.

\demo{Proof of theorem 3.1,(1)}
In this proof we apply the strong Whitney $C^\infty$ topology unless otherwise specified. 
The idea of the proof is taken from [S$_1$]. 
Let us consider the case without corners. 
The proof is divided into three steps. 
Denote by $X$ and $Y$ the extended critical sets of $f$ and $g$, that is $X=f^{-1}(f(\Sing f))$ 
and $Y=g^{-1}(g(\Sing g))$. 
Note that $X$ and $Y$ are not necessarily analytic sets. 
Let $M$ be analytic and closed in the ambient Euclidean space $\R^N$, and consider the Riemannian metric 
on $M$ induced from that of $\R^N$. Set $n=\dim M$. \par
\enddemo
\example{Step 1}
Assume that $X$ is an analytic set. 
Let $\phi$ denote a $C^\infty$ diffeomorphism of $M$ such that $f=g\circ\phi$. 
Then there exists a $C^\omega$ diffeomorphism $\pi$ of $M$ arbitrarily close to $\phi$ such that 
$\pi(X)=Y$. 
\endexample
{\it Proof of step 1.} 
Let $\{X_i:i=0,...,n-1\}$ and $\{Y_i:i=1,...,n-1\}$ be the canonical stratifications of $X$ and $Y$ 
respectively, and put $X_n=M-X$ and $Y_n=M-Y$. \par
Before beginning the proof, we give some definitions and a remark. 
Fix $X_i$. 
Let $\{M_l\}$ be a family of compact $C^\infty$ manifolds of dimension $i$ possibly with boundary 
such that $\{\Int M_l\}$ is a locally finite covering of $\cup_{j=0}^iX_j$. 
A function on $\cup_{j=0}^iX_j$ is called {\it of class} $C^\infty$ if its restriction to each 
$M_l$ is of class $C^\infty$. 
Thus $C^\infty(\cup_{j=0}^iX_j)$ is a subset of $\prod_lC^\infty(M_l)$. 
We give to $C^\infty(\cup_{j=0}^iX_j)$ the product topology of the $C^\infty$ topology on each 
$C^\infty(M_l)$, i.e. the compact-open $C^\infty$ topology. 
We give also the strong Whitney $C^\infty$ topology on $C^\infty(\cup_{j=0}^i X_j)$ in the same way. 
Then lemma 2.12 holds for the map $C^\infty(\cup_{j=0}^iX_j)\ni f\to h|_{\cup_{j=0}^iX_j}f\in C^\infty(\cup_{j=0}^iX_j)$ for an analytic function $h$ on $M$, which is proved in the same way. 
We will use this generalized version of the lemma below. \par
Let $X'$ be a normal-crossing $C^\omega$ subset of $M$ contained in $X$. 
Assume that the sheaf of $\Cal O$-ideals on $M$ defined by $X'$ is generated by a single $C^\omega$ 
function $\xi$ on $M$. 
Let $V$ denote the subspace of $C^\infty(\cup_{j=0}^iX_j)$ consisting of functions which vanish 
on $X'$. 
Then $V=\xi C^\infty(\cup_{j=0}^iX_j)$ by Theorem VI,3.10 in [Ml]. 
We will use this remark later in this proof. \par
Now we begin the proof. 
By induction, for some $i\in\N$, assume that we have constructed a $C^\infty$ diffeomorphism $\pi_{i-1}$ 
of $M$ close to $\phi$ such that $\pi_{i-1}|_{\cup_{j=0}^{i-1}X_j}$ is of class $C^\omega$ (in the 
sense that $\pi_{i-1}|_{\cup_{j=0}^{i-1}X_j}\in\prod_lC^\omega(M_l))$ and $\pi_{i-1}(X_j)=Y_j$ for 
any $j$. 
Let $\Cal M$ denote the sheaf of $\Cal O$-ideals on $M$ defined by $\cup_{j=0}^{i-1}X_j$, which 
is coherent because $X$ is normal-crossing. 
Then $\pi_{i-1}|_{\cup_{j=0}^{i-1}X_j}\in H^0(M,\Cal O/\Cal M)^N$ for the following reason. 
As the problem is local, we can assume that $M=\R^n$ and $X=\{(x_1,...,x_n)\in\R^n:x_1\cdots x_{n'}=0\}$ 
for some $n'\le n\in\N$. 
Moreover, we suppose that $i=n$ because if for each irreducible analytic component $E$ of $\cup_{j=0}^i
X_j$ we can extend $\pi_{i-1}|_{E\cap\cup_{j=0}^{i-1}X_j}$ to an analytic map on $E$ then the 
extensions for all $E$ define an analytic map from $\cup_{j=0}^iX_j$
to $\R^N$, and hence it suffices to work on each $E$ in place of $\R^n$. 
Then what we see is that an analytic function on $X=\{(x_1,...,x_n)\in\R^n:x_1\cdots x_{n'}=0\}$ is an element of $H^0(M,\Cal O/\Cal M)$, namely, can be extensible to an analytic function on $\R^n$ (corollary 2.4), where $\Cal M$ is defined by $X$. 
We proceed by induction on $n'$. 
Since the statement is clear if $n'=0$ or $n'=1$, assume that $n'>1$
and that the restriction of $f$ to $\{x_1=0\}$ is extensible to an analytic function $\tilde f$ on $\R^n$. 
Then $f-\tilde f|_X$ vanishes on $\{x_1=0\}$ and hence is divisible by $x_1$. 
Apply the induction hypothesis to $(f-\tilde f)/x_1$ and $\{x_2\cdots x_{n'}=0\}$ and let $\Tilde{\Tilde f}$ be an extension of $(f-\tilde f)/x_1$. 
Then $\tilde f+x_1\Tilde{\Tilde f}$ is the required extension of $f$. \par
  Consider any $C^\omega$ extension $\alpha:M\to\R^N$ of $\pi_{i-1}|_{\cup_{j=0}^{i-1}X_j}$ (corollary 
2.4). 
Here we can choose $\alpha$ to be sufficiently close to $\pi_{i-1}$ and so that $\Im\alpha
\subset M$ for the following reason. 
Let $\gamma_1,...,\gamma_k\in C^\omega(M)$ be generators of $\Cal M$ (corollary 2.2). 
Then there exist $\delta_1,...,\delta_k\in C^\infty(M,\R^N)$ such that $\pi_{i-1}-\alpha=
\sum_{j=1}^k\gamma_j\delta_j$ by the above remark. 
Let $\tilde\delta_j$ be $C^\omega$ approximations of $\delta_j$. 
Replace $\alpha$ with the composite of $\alpha+\sum\gamma_j\tilde\delta_j$ and the orthogonal projection of a neighborhood of $M$ in $\R^N$ to $M$. 
Then $\alpha$ satisfies the requirements. 
Let $p_j:U_j\to Y_j$ be the orthogonal projection of a tubular neighborhood of $Y_j$ in $\R^N$. 
Here $U_j$ is described as $\cup_{y\in Y_j}\{x\in\R^N:\ |x-y|<\epsilon(y),\,(x-y)\perp T_yY_j\}$ 
for some positive $C^0$ function $\epsilon$ on $Y_j$ where $T_yY_j$ denotes the tangent space of $Y_
j$ at $y$, and we can choose $\epsilon$ so large that $\epsilon(y)\ge\epsilon_0\dis(y,\cup_{k=0}^{j-
1}Y_k)$ locally at each point of $\cup_{k=0}^{j-1}Y_k$ for some positive number $\epsilon_0$ because 
$Y$ is normal crossing. 
Then $\alpha$ can be so close to $\pi_{i-1}$ that $\alpha(X_i)\subset U_i$ since $\alpha=\pi_
{i-1}$ on $\cup_{j=0}^{i-1}X_j$ and hence $d\alpha_x v=d\pi_{i-1 x}v$ for any $x\in\cup_{j=0}^{i-1}
X_j$ and for any tangent vector $v$ at $x$ tangent to $\cup_{j=0}^{i-1}X_j$. 
Define $\pi_i$ on $\cup_{j=0}^iX_j$ to be $\alpha$ on $\cup_{j=0}^{i-1}X_j$ and $p_i\circ\alpha$ 
on $X_i$. 
Note that $\pi_i$ is a $C^\omega$ map from $\cup_{j=0}^i X_j$ to $\cup_{j=0}^i Y_j\subset\R^N$ 
and close to $\pi_{i-1}|_{\cup_{j=0}^iX_j}$. Actually $\pi_{i-1}|_{\cup_{j=0}^iX_j}=\overline{p_i}
\circ\pi_{i-1}|_{\cup_{j=0}^iX_j}$ and moreover since $\alpha$ is
close to $\pi_{i-1}$, then $\overline{p_i}\circ\alpha$ on $\cup_{j=0}^iX_j$ $(=\pi_i)$ is close to $\overline{p_i}\circ\pi
_{i-1}$ on $\cup_{j=0}^iX_j$, where $\overline{p_i}:\overline{U_i}\to\cup_{j=0}^iY_j$ is the natural 
extension of $p_i$. 
We need to extend $\pi_i$ to a $C^\infty$ diffeomorphism of $M$ which is close to $\pi_{i-1}$ and 
carries each $X_j$ to $Y_j$. 
Compare $\pi_i\circ\pi_{i-1}^{-1}|_{\cup_{j=0}^iY_j}$ and the identity map of $\cup_{j=0}^iY_j$. 
Then they are close each other and what we have to prove is the following statement: let $\tau$ be a $C^\infty$ map between $\cup_{j=0}^iY_j$ close to $\id$. 
Then we can extend $\tau$ to a $C^\infty$ diffeomorphism of $M$ which is close to $\id$ and 
carries each $Y_j$ to $Y_j,\, j=i+1,...,n$. \par
By the second induction, it suffices to extend $\tau$ to a $C^\infty$ map between $\cup_{j=0}^
{i+1}Y_j$ close to $\id$. 
We reduce the problem to a trivial case. 
First it is enough to extend $\tau$ to a $C^\infty$ map from $\cup_{j=0}^{i+1}Y_j$ to $\R^N$ 
close to $\id$ by virtue of $p_{i+1}:U_{i+1}\to Y_{i+1}$ as above. 
Secondly, if we replace $\tau$ with $\tau-\id|_{\cup_{j=0}^iY_j}$ then the problem is that 
for a $C^\infty$ map $\tau:\cup_{j=0}^iY_j\to\R^N$ close to the zero map we can extend $\tau$ to 
a $C^\infty$ map from $\cup_{j=0}^{i+1}Y_j$ to $\R^N$ close to 0. 
Thirdly, we can assume that $\tau$ is a function. 
Hence we can use a partition of unity of class $C^\infty$ and the problem becomes local (remark 2.11,(5)). 
So we assume that $M=\R^n$, $Y$ is the union of some irreducible analytic components of $\{y_1\cdots 
y_n=0\}$ and $\tau\in C^\infty(\cup_{j=0}^iY_j)$ is close to 0 and vanishes on $\{y\in\cup_{j=0}^i
Y_j\subset\R^n:\ |y|\ge1\}$, where $(y_1,...,y_n)\in\R^n$. 
Let $\xi$ be a $C^\infty$ function on $M$ such that $\xi=1$ on $\{y\in\R^n:\ |y|\le 1\}$ and 
$\xi=0$ on $\{|y|\ge 2\}$. 
If $n=0$ or 1 we have nothing to do. 
Hence by the third induction on $n$, assume that we have a $C^\infty$ extension $\tau_1(y_2,...,y_n)$ 
of $\tau|_{\{y_1=0\}\cap\cup_{j=0}^iY_j}$ to $\{y_1=0\}\cap\cup_{j=0}^{i+1}Y_j$ close to 0. 
Regard $\tau_1(y_2,...,y_n)$ as a $C^\infty$ function on $\cup_{j=0}^{i+1}Y_j$, which is possible 
because $\cup_{j=0}^{i+1}Y_j$ is contained in the product of $\R$ and the image of $\{y_1=0\}\cap\cup_{j=0}^
{i+1}Y_j$ under the projection $\R\times\R^{n-1}\to\R^{n-1}$. 
Replace $\tau$ with $\tau-\tau_1\xi$, which vanishes on $\{y_1=0\text{ or }|y|\ge 2\}\cap\cup_{j=0
}^iY_j$. 
Next consider $(\tau-\tau_1\xi)/y_1|_{\{y_2=0\}\cap\cup_{j=0}^iY_j}$ and apply the generalized lemma 2.12, the above remark and the 
same arguments as above. 
Then we reduce the problem to the case where $\tau=0$ on $\{y_1y_2=0\text{ or }|y|\ge3\}\cap\cup_
{j=0}^iY_j$ and by the fourth induction to the case where $\tau=0$ on $\cup_{j=0}^iY_j$. 
Thus step 1 is proved. 
\example{Step 2}
Assume that $X=Y$, $X$ is an analytic set and there exists a $C^\infty$ diffeomorphism $\phi$ of $M$ 
such that $f=g\circ\phi$. 
Then there exists a $C^\omega$ diffeomorphism $\pi$ of $M$ close to $\phi$ such that $f=g\circ\pi$. 
\endexample
{\it Proof of step 2.} 
By step 1 there exists a $C^\omega$ diffeomorphism $\phi'$ of $M$ arbitrarily close to $\phi$ such that $\phi'(X)=X$. 
Then $f\circ\phi^{\prime-1}=g\circ\phi\circ\phi^{\prime-1}$, $f\circ\phi^{\prime-1}$ is analytic and $\phi\circ\phi^{\prime-1}$ is close to id. 
Hence we assume in step 2 that $g$ is fixed and $f$ and $\phi$ can be chosen so that $\phi$ and $f-g$ are arbitrarily close to id and 0 respectively. 
We construct $\pi$ by integrating along a well-chosen vector field on $M$. 
There exist analytic vector fields $w_1,...,w_N$ on $M$ which span the tangent space at each point 
of $M$, e.g., $w_{i x}=dp_x\frac{\partial\ \ }{\partial x_i},\, x\in M$, where $(x_1,...,x_N)\in\R^N$ 
and $p$ is the orthogonal projection of a tubular neighborhood of $M$ in $\R^N$. 
Consider a vector field $v=\frac{\partial\ }{\partial t}+\sum_{i=1}^N a_i w_i$ on $M\times[0,\,1]$ 
where $a_i\in C^\omega(M\times[0,\,1])$ for $i\in\{1,...,N\}$. 
Put $F(x,t)=(1-t)f(x)+t g(x)$ for $(x,t)\in M\times[0,\,1]$. \par
Assume that we have found such $a_i,\,i=1,...,N$, that $v(F)=0$ and $|\sum_{i=1}^N a_i w_i|$ is small. 
Then $F$ is constant along integral curves of $v$, therefore, the flow of $v$ furnishes an analytic 
diffeomorphism $\pi$ so that $f=g\circ\pi$. \par
Therefore, what we have to do is to construct the relevant $a_i,\,i\in\{1,...,N\}$. 
First look at the local case. 
We will show that there exist a compact neighborhood $U$ of each point of $M$ and $a_i\in C^\omega(U\times[
0,\,1])$, $i=1,...,N$, such that $v(F)=0$ on $U\times[0,\,1]$ and $|\sum_{i=1}^N a_i w_i|$ is small. 
If the point is in $X$, we can write $U=\{x\in\R^n:\ |x|\le 1\},\ g(x)-c=\prod_{i=1}^n x_i^{n_i}$ and 
$f(x)-c=\lambda(x)(g(x)-c)$ for $x\in U$, where $c\in\R$, $\lambda$ is a $C^\omega$ function on $U$ 
and close to 1 by lemma 2.12, and at least one of $n_i$'s, say $n_1$, is non-zero. 
Assume that $c=0$ without loss of generality. 
Then there exists $v$ of the form $\frac{\partial\ }{\partial t}+b_1\frac{\partial\ \ }{\partial 
x_1},\ b_1\in C^\omega(U\times[0,\,1])$, which satisfies $v(F)=0$. 
Actually
$$
\gather
(\frac{\partial\ }{\partial t}+b_1\frac{\partial\ \ }{\partial x_1})F(x,t)=\qquad\qquad\qquad\qquad
\qquad\qquad\qquad\qquad\qquad\qquad\qquad\\
\qquad(1-\lambda)g(x)+b_1(x,t)\big( n_1(t+(1-t)\lambda(x))\frac{g(x)}{x_1}+(1-t)g(x)\frac{\partial
\lambda\ }{\partial x_1}(x)\big) =0,\\
b_1(x,t)=\frac{-(1-\lambda(x))x_1}{n_1(t+(1-t)\lambda(x))+(1-t)x_1\frac{\partial\lambda\ }
{\partial x_1}(x)},
\endgather
$$
which is an analytic function in $U\times[0,\,1]$ and close to 0. 
Shrink $U$ if necessary. 
Then for some $0<i_1<\cdots<i_n\le N$, the vector fields $w_{i_1},...,w_{i_n}$ span the tangent space 
there, and $b_1\frac{\partial\ \ }{\partial x_1}$ is described uniquely by $\sum_{j=1}^n a_{i_j}w_{i_j}$ for 
some $C^\omega$ functions $a_{i_j}$. 
Hence $a_{i_j},\ j=1,...,n$, and $a_i=0,\ i\not\in\{i_1,...,i_n\}$, fulfill the requirements. \par
Next consider the situation at a point outside of $X$. 
Note that the values of $f$ and $g$ at the point may be different, and hence the above arguments do not work. 
We can choose its local coordinate system so that $U=\{x\in\R^n:\ |x|\le 1\}$ and $\frac{
\partial g\ }{\partial x_1}=1$ on $U$. 
Then 
$$
\gather
(\frac{\partial\ }{\partial t}+b\frac{\partial\ \ }{\partial x_1})F(x,t)=-f+g+b_1((1-t)\frac{\partial 
f\ }{\partial x_1}+t\frac{\partial g\ }{\partial x_1})=0,\\
b_1(x,t)=\frac{f-g}{(1-t)\frac{\partial f\ }{\partial x_1}+t\frac{\partial g\ }{\partial x_1}},
\endgather
$$
and $(1-t)\frac{\partial f\ }{\partial x_1}+t\frac{\partial g\ }{\partial x_1}$ and $f-g$ are close 
to 1 and 0, respectively. 
Hence there exist $U$ and $a_i$ as before. \par
Consequently, using a partition of unity of class $C^\infty$ we obtain a $C^\infty$ vector 
field $v'=\frac{\partial\ }{\partial t}+\sum_{i=1}^N a'_i w_i$ on $M\times[0,\,1]$ such that $v(F)=0$ 
and $|\sum_{i=1}^N a'_i w_i|$ is small (remark 2.11,(5)). \par
Now, to construct the global analytic vector filed $v$ on $M\times[0,\,1]$ we use Cartan Theorems A 
and B. 
Consider the sheaf of relations $\Cal J$ on $M\times[0,\,1]$ defined by
$$
\multline
\Cal J=\cup_{(x,t)\in M\times[0,\,1]}\{(\beta,\alpha_1,\ldots,\alpha_N)\in\Cal O_{(x,t)}^{N+1}:\\ 
\beta(f_x-g_x)+\sum_{i=1}^N\alpha_i(w_i((1-t)f+tg))_{(x,t)}=0\}. 
\endmultline
$$
The sheaf $\Cal J$ is a coherent sheaf of $\Cal O$-modules by Oka Theorem 2.5. 
Later we will find $l\in \N$ and global cross-sections $(b_k,a_1^k,\ldots,a_N^k)\in H^0(M\times[0,\,1],
\Cal J)$, $k\in\{1,\ldots,l\}$, such that for any $(x,t)\in M\times[0,\,1]$, any $C^\omega$ vector field 
germ $\omega$ at $(x,t)$ in $M\times[0,\,1]$ with $\omega(F_{(x,t)})=0$ is of the form $\sum_{k=1}^l\xi_
k v_{k(x,t)}$ for some $C^\omega$ function germs $\xi_k$ at $(x,t)$ in $M\times[0,\,1]$, where $v_k=b_k
\frac{\partial\ }{\partial t}+\sum_{i=1}^N a_i^k w_i$. 
Assume the existence of such $l$ and $(b_k,a_1^k,\ldots,a_N^k)$. 
Then by the above method of construction of $v'$ and by a partition of unity of class $C^\infty$ 
there exist $C^\infty$ functions $\theta_k$ on $M$ such that $v'=\sum_{k=1}^l\theta_k v_k$. 
Approximate $\theta_k$ by $C^\omega$ functions $\tilde\theta_k$, and set $\tilde v=\sum_{k=1}^l
\tilde\theta_k v_k$. 
Then $\tilde v$ is a $C^\omega$ vector field close to $v'$ such that $F(\tilde v)=0$ and is described by 
$\tilde a_0\frac{\partial\ }{\partial t}+\sum_{i=1}^N\tilde a_i w_i,\ \tilde a_i\in C^\omega(M\times[0,
\,1])$, for the following reason. 
Let $\Cal I$ denote the coherent sheaf of $\Cal O$-submodules of the sheaf of $\Cal O$-modules 
of germs of $C^\omega$ vector fields on $M\times[0,\,1]$ defined by 
$$
\Cal I_{(x,t)}=\{\omega:\ \omega(F_{(x,t)})=0\}\quad\text{for}\ (x,t)\in M\times[0,\,1],
$$
and define an $\Cal O$-homomorphism $\delta:\Cal O^l\to\Cal I$ by 
$$
\delta(\gamma_1,...,\gamma_l)=\sum_{k=1}^l\gamma_k v_{k(x,t)}\quad\text{for}\ (\gamma_1,...,\gamma
_l)\in\Cal O_{(x,t)}^l,\ (x,t)\in M\times[0,\,1].
$$
Then $\delta$ is surjective, $H^0(M\times[0,\,1],\Cal I)$ is the set of all $C^\omega$ vector fields 
$w$ on $M\times[0,\,1]$ with $w(F)=0$, and hence by Cartan Theorem B the homomorphism $C^\omega(M
\times[0,\,1])^l\ni(d_1,...,d_l)\to\sum_{k=1}^l d_k v_k\in H^0(M\times[0,\,1],\Cal I)$ is surjective, 
i.e., $\tilde v$ is of the form $\sum_{k=1}^l d_k v_i$ for some $C^\omega$ functions $d_k$ on $M
\times[0,\,1]$. 
Therefore, we have $\tilde v=\tilde a_0\frac{\partial\ }{\partial t}+\sum_{i=1}^N\tilde a_i w_i$ for 
$\tilde a_0=\sum_{k=1}^l d_k$ and $\tilde a_i=\sum_{k=1}^l d_k a_i^k,\ i=1,...,N$. 
Here $\tilde a_0$ is unique and hence close to 1, and $|\sum_{i=1}^N\tilde a_i w_i|$ is small. 
Thus $v=\frac{\partial\ }{\partial t}+\sum_{i=1}^N(\tilde a_i/\tilde a_0)w_i$ is what we wanted. \par
It remains to find $(b_k,a_1^k,...,a_N^k),\,k=1,...,l$. 
That is equivalent to prove that $H^0(M\times[0,\,1],\Cal I)$ is finitely generated by Cartan Theorem B 
because the homomorphism $\Cal J_{(x,t)}\ni(\beta,\alpha_1,...,\alpha_N)\to\beta\frac{\partial\ }
{\partial t}+\sum_{i=1}^N\alpha_i w_{i(x,t)}\in\Cal I_{(x,t)}$ is surjective. 
Moreover, it suffices to see that each stalk $\Cal I_{(x,t)}$ is generated by a uniform number of elements 
by corollary 2.2. 
Note that $F$ is an analytic function with only normal crossing singularities. 
Hence we replace $F$ with $f$ to simplify notation. 
Let $\Cal K$ denote the sheaf of $\Cal O$-modules of $C^\omega$ vector field germs on $M$ such that 
$$
\Cal K_x=\{\omega:\omega(f_x)=0\}\quad\text{for}\ x\in M. 
$$
Then it suffices to choose $l\in\N$ so that for any $x_0\in M$, $\Cal K_{x_0}$ is generated by $l$ 
elements. 
Since the problem is local we can assume that $M=\R^n,\ x_0=0$ and $f(x)=\prod_{i=1}^k x_i^{n_i}$ with 
$n_1,...,n_k>0,\ 0<k\le n$. 
Write $\omega\in\Cal K_0$ as $\sum_{i=1}^n\alpha_i\frac{\partial\quad}{\partial x_i},\ \alpha_i\in
\Cal O_0$, and set $h(x)=\prod_{i=1}^kx_i$. 
Then $\omega(f_{x_0})=0$ means 
$$
\sum_{i=1}^k n_i\alpha_i f(x)/x_i=0,\quad\text{hence}\ \sum_{i=1}^k n_i\alpha_i h(x)/x_i=0. 
$$
Therefore, each $\alpha_i$ is divisible by $x_i$. 
Hence, setting $\alpha'_i=\alpha_i/x_i$ we obtain $\sum_{i=1}^kn_i\alpha'_i\allowmathbreak=0$. 
It is clear that $\{(\alpha'_1,...,\alpha'_n)\in\Cal O_0^n:\sum_{i=1}^k\alpha'_i=0\}$ is generated 
by $n-1$ elements, which proves step 2. \par
The proof of step 1 shows that any $C^\infty$ diffeomorphism of $M$ carrying $\Sing f$ to $\Sing g$ 
is approximated by an analytic diffeomorphism of $M$ with the same property. 
Hence it suffices to prove the next statement. \par
\example{Step 3}
Assume that $\Sing f=\Sing g$ and there exists a $C^\infty$ diffeomorphism $\phi$ of $M$ such that $f=g\circ\phi$. 
Then there exists a $C^\omega$ diffeomorphism $\pi$ of $M$ such that $f=g\circ\pi$. 
\endexample
{\it Proof of step 3.} 
As at the beginning of the proof of step 2 we fix $g$ and modify $f$ and $\phi$ so that $\phi$ and $f-g$ are sufficiently close to id and 0 respectively. 
Set $Z=\Sing f$ and let $Z_i,\ i=1,2,...,$ be connected components of $Z$. 
Let $U_i$ be disjoint small open neighborhoods of $Z_i$ in $M$ such that if $\phi(U_i)\cap U_{i'}
\not=\emptyset$ then $i=i'$. 
Then by steps 1 and 2 there exist $C^\omega$ diffeomorphisms $\pi_i:U_i\to\phi(U_i)$ close to 
$\phi|_{U_i}:U_i\to\phi(U_i)$ such that $f=g\circ\pi_i$ on $U_i$. 
Note that if we define a map between $M$ to be $\pi_i$ on each $U_i$ and $\phi$ elsewhere, then the map 
is a $C^\infty$ diffeomorphism by the definition of the strong Whitney $C^\infty$ topology. 
For $x_0\in M$, let $m(x_0)$ denote the multiplicity of $g-g(x_0)$ at $x_0$, i.e., $m(x_0)=|\alpha
|=\alpha_1+\cdots+\alpha_n$ for $\alpha=(\alpha_1,...,\alpha_n)\in\N^n$ such that $g(x)-g(x_0)$ is 
written as $\pm x^\alpha$ for some local coordinate system $(x_1,...,x_n)$ at $x_0$. 
There exists $h\in C^\omega(M)$ such that $h^{-1}(0)=Z$ and $h$ is $m(x)$-flat at each $x\in Z$ 
for the following reason. 
For each $i$, let $\{Z_{i,j}\}_j$ be the stratification of $Z_i$ by multiplicity number, and for 
each $Z_{i,j}$, consider the smallest analytic set in $U_i$ and hence in $M$ containing each 
connected component of $Z_{i,j}$. 
Then we have a locally finite decomposition of $Z_i$ into irreducible analytic sets $\{W_{i,j}\}_j$ 
in $M$ such that $m(x)$ is constant, say $m_{i,j}$, on each $W_{i,j}-\cup_{j'}\{W_{i,j'}:\dim W_{i,
j'}<\dim W_{i,j}\}$. 
By corollary 2.2 there exists $h_{i,j}\in C^\omega(M)$---e.g., the $m_{i,j}$th power of the square 
sum of a finite number of global generators of the sheaf of $\Cal O$-ideals defined by $W_{i,j}
$---such that $h_{i,j}^{-1}(0)=W_{i,j}$ and $h_{i,j}$ is $m_{i,j}$-flat at $W_{i,j}$, and then 
considering the sheaf of $\Cal O$-ideals $\prod_{i,j}h_{i,j}\Cal O$ we obtain $h$ in the same way. 
\par
We will reduce the problem to the case where $\pi_i-\id$ on $U_i$ and $f-g$ are divisible by $h$. 
Since $\supp\Cal O/h\Cal O=Z$, $\{\pi_i\}_i$ defines an element of $H^0(M,(\Cal O/h\Cal O)^N)$. 
Hence applying Cartan Theorem B to the exact sequence $0\to(h\Cal O)^N\to\Cal O^N\to(\Cal O/h\Cal 
O)^N\to 0$, we obtain $\pi'\in C^\omega(M)^N$ such that $\pi_i-\pi'\in hC^\omega(U_i)^N$ for each 
$i$. 
We need to modify $\pi'$ to be a diffeomorphism of $M$. 
Let $\xi$ be a $C^\infty$ function on $M$ such that $\xi=0$ outside of a small neighborhood of $Z$ 
and $\xi=1$ on a smaller one. 
Approximate $C^\infty$ maps $\sum_i\xi(\pi_i-\pi')/h$ and $(1-\xi)(\phi-\pi')/h$ from $M$ to 
$\R^N$ by $C^\omega$ maps $H_1$ and $H_2$, respectively. 
Then $h H_1+h H_2+\pi'$ is an analytic approximation of $\phi:M\to\R^N$ whose difference with 
$\pi_i$ on $U_i$ is divisible by $h$. 
Hence its composite $\pi''$ with the orthogonal projection of $M$ in $\R^N$ is an analytic 
approximation of $\phi:M\to M$ such that $\pi_i-\pi''$ is divisible by $h$ by the next fact. 
Given $\theta_1,\theta_2\in\R\langle\langle x_1,...,x_n\rangle\rangle^m,\ \eta\in\R\langle\langle x_1,...,x_n
\rangle\rangle$ and $\rho\in\R\langle\langle y_1,...,y_m\rangle\rangle$ with $\theta_1(0)=\theta
_2(0)=\eta(0)=\rho(0)=0$, then $\rho\circ\theta_1-\rho\circ(\theta_1+\eta\theta_2)$ is divisible by $\eta$ as an element of $\in\R\langle\langle 
x_1,...,x_n\rangle\rangle$. 
Now replace $\phi,\ \pi_i$ and $f$ with $\phi\circ\pi^{\prime\prime-1},\ \pi_i\circ\pi
^{\prime\prime-1}$ and $f\circ\pi^{\prime\prime-1}$, respectively. 
Then the equalities $f=g\circ\phi$ and $f=g\circ\pi_i$ continue to hold, and $\pi_i-\id$ and 
hence $f-g$ are divisible by $h$ and, moreover, by $h^{3+s}$ by the same way, where $s\in\N$ is 
such that $\alpha^s$ is contained in the ideal of $\R\langle\langle x_1,...,x_n\rangle\rangle$ 
generated by $\frac{\partial\psi_l}{\partial x_1},...,\frac{\partial\psi_l}{\partial x_n}$ 
for $\psi_l(x)=\prod_{j=1}^l x_j\in\R\langle\langle x_1,...,x_n\rangle\rangle,\, l\le 
n$, and for $\alpha\in\R\langle\langle x_1,...,x_n\rangle\rangle$ which vanishes on $\Sing\psi_l
$ (Hilbert Zero Point Theorem). 
Set $h_1=(f-g)/h^{3+s}\in C^\omega(M)$, which is close to 0 by lemma 2.12. \par
As in the proof of step 2, we define $C^\omega$ vector fields $w_i,\,i=1,...,N$, and a $C^\omega$ 
function $F$ on $M\times[0,\,1]$, and it suffices to find a $C^\omega$ vector field $v$ of the 
form $\frac{\partial\ }{\partial t}+\sum_{i=1}^N a_i w_i$ on $M\times[0,\,1]$ such that $v(F)=0$ 
and $|\sum_{i=1}^N a_i w_i|$ is bounded. 
Since $f=g+h^{3+s}h_1$, then $F=g+(1-t)h^{3+s}h_1$, and the equality $v(F)=0$ becomes 
$$
h^{3+s}h_1=\sum^N_{i=1}a_i(w_i g+(1-t)h^{2+s}h_{2,i})
$$
for some $C^\omega$ functions $h_{2,i}$ on $M$ close to 0. 
This is solvable locally. 
Indeed, for each $x_0\in M-Z$, at least one of $w_i g$, say $w_1 g$, does not vanish at $x_0$. 
Hence $a_1=h^{3+s}h_1/(w_1 g+(1-t)h^{2+s}h_{2,1}),\ a_2=\cdots=a_N=0$ is a solution on a 
neighborhood of $x_0$. 
Assume that $x_0\in Z$. 
Then choose an analytic local coordinate system $(x_1,...,x_n)$ at $x_0$ in $M$ so that $g(x)=
\prod_{i=1}^n x_i^{n_i}+\const$, where $\sum_{i=1}^n n_i=m(x_0)>1$. 
Here we can assume that $x_0=0,$ $\const=0$, $n_1,...,n_l>0$ and $n_{l+1}=\cdots=n_n=0$. 
Note that $m(0,...,0,x_{l+1},...,x_n)=m(0)$ for $(x_{l+1},...,x_n)\in\R^{n-l}$ near 0. 
What we prove is that for each $t_0\in[0,\,1]$, the ideal $I$ of $\Cal O_{(0,t_0)}$ generated by 
$\frac{\partial g_{(0,t_0)}}{\partial x_i\ \quad}+(1-t_{(0,t_0)})h_{(0,t_0)}^{2+s}h_{2,i(0,t_0)},\ 
i=1,...,l$, contains $h_{(0,t_0)}^{3+s}h_{1(0,t_0)}$. 
Let $J$ denote the ideal of $\Cal O_{(0,t_0)}$ generated by $\frac{\partial g_{(0,t_0)}}{\partial 
x_i\ \quad},\ i=1,...,l$. 
Then it suffices to see that $h^{1+s}_{(0,t_0)}\in J$ because if so, $J\supset I,\ J\ni h_{(0,t_0)}^{
3+s}h_{1(0,t_0)},\ J=I+\frak m J$ and hence by Nakayama lemma $J=I$, where $\frak m$ is the 
maximal ideal of $\Cal O_{(0,t_0)}$. 
Moreover, assuming $g(x)=x_1\cdots x_l$ we prove that $h_{(0,t_0)}^s\in J$, which is sufficient because 
$\frac{\partial g_{(0,t_0)}}{\partial x_i\quad\ }=n_i\prod_{j=1}^l x_j^{n_j-1}\frac{\partial x_1
\cdots x_l}{\partial x_i\quad\ }$ and $h_{(0,t_0)}$ is divisible by $\prod_{j=1}^l x_j^{n_j-1}$ 
by the definition of $h$. 
However, $h_{(0,t_0)}^s\in J$ is clear by the definition of $s$. 
Note that we can choose local $v=\frac{\partial\ }{\partial t}+\sum_{i=1}^N a_i w_i$ in any case so 
that $|\sum_{i=1}^N a_i w_i|$ is arbitrarily small. \par
We continue to proceed in the same way as in the proof of step 2. 
We obtain $C^\omega$ vector fields $v_k=b_k\frac{\partial\ }{\partial t}+\sum_{i=1}^N a_i^k w_i,
\ k=1,...,l$, by local existence and a $C^\infty$ vector field $v'=\frac{\partial\ }{\partial t}
+\sum_{i=1}^N a'_i w_i$ such that $v_k(F)=v'(F)=0,\ |\sum_{i=1}^N a'_i w_i|$ is small and $v'$ 
is of the form $\sum_{k=1}^l\theta_k w_i$. 
After then we approximate $\theta_k$ by $C^\omega$ functions $\tilde\theta_k$, and $v=\sum_{k=1}
^l\tilde\theta v_k$ fulfills the requirements. 
Thus we complete the proof of (1) in the case of without corners. \par
The case with corners is proved in the same way. 
\qed

%%%%%%%%%%%%%%%%%%%%%%%%%%%%%%%%%%%%%%%%%%%%%%%%%%%%%%%%%%%%%%%%%%%%%%%%%%%%%%%%%%%%%%%%%%%%%%%%%%%
%%%%%%%%%%%%%%%%%%%%%%%%%%%%%%%%%%%%%%%%%%%%%%%%%%%%%%%%%%%%%%%%%%%%%%%%%%%%%%%%%%%%%%%%%%%%%%%%%%%

\subhead
3.2. Cardinality of the set of equivalence classes
\endsubhead\par
Our main theorem establishes the cardinality of analytic (respectively Nash) R-L equivalence classes of analytic (respectively Nash) 
functions on $M$ with only normal crossing singularities. 
\proclaim{Theorem 3.2} 
Let $M$ be a compact analytic (respectively, Nash) manifold of strictly positive dimension. 
Then the cardinality of analytic (resp., Nash) R-L equivalence classes of analytic (resp. Nash) 
functions on $M$ with only normal crossing singularities is 0 or countable. 
In the Nash case, the compactness of $M$ is not necessary, and if moreover $M$ is non-compact then the 
cardinality is countable. 
\endproclaim
The proof of theorem 3.2 runs as follows. 
We reduce the $C^\omega$ case to the Nash case by proposition 4.1 and then the non-compact Nash 
case to the compact Nash case by proposition 4.9. 
Lemmas 4.3 and 4.4 together with Nash Approximation Theorem II prove
the compact Nash case. We postpone the proof of theorem 3.2 to the
last part of the paper.

\remark{Remark}
(i) The case where the cardinality is zero may appear, e.g. $M=S^2,\P(2)$ 
(for the proof, see the arguments in (v) below in case $M=\R^2$). \par
(ii) In the theorem we do not need to fix $M$, namely, the cardinality of equivalence classes of 
analytic or Nash functions on all compact analytic manifolds or Nash manifolds, respectively, with 
only normal crossing singularities is also countable. 
Indeed, the cardinality is clearly infinite, and there are only a countable number of compact 
analytic manifolds and (not necessarily compact) Nash manifolds up to analytic diffeomorphism 
and Nash diffeomorphism, respectively, which will be clear in the proof of lemma 4.4. \par
(iii) Theorem 3.2 does not hold for analytic functions on a non-compact analytic manifold. 
To be precise, for a non-compact analytic manifold $M$, the cardinality of analytic R-L equivalence 
classes of (proper) analytic functions on $M$ with only normal crossing singularities is of the 
continuum (0 or of the continuum, respectively). We prove this fact below. \par
(iv) On any non-compact connected 
analytic (Nash, respectively,) manifold $M$, there exists a
non-singular analytic (Nash, respectively,) function. We give the construction
below. \par
(v) An example of non-compact $M$ where there is no proper analytic (Nash) function with only 
normal crossing singularities is $\R^2$. 
We see this by reduction to absurdity. 
Assume that there exists such an $f$. 
Note that each level of $f$ is a finite union of Jordan curves. 
Let $a_1\in\R$ be a point of $\Im f$ and $X_1\subset\R^2$ be a Jordan curve in 
$f^{-1}(a_1)$ that does not intersect with $f^{-1}(a_1)$
inside of $U_1$. 
Next choose $a_2\in f(U_1)$, a Jordan curve $X_2$ in $f^{-1}(a_2)\cap U_1$ and $U_2$ in the same 
way. 
If we continue these arguments, we arrive at a contradiction to the above note. \par
{\it Proof of (iii) for proper functions.} 
Assume that there exists a proper analytic function $f$ on a non-compact analytic manifold $M$ with only 
normal crossing singularities. 
Replacing $f$ with $\pi\circ f$ for some proper analytic function
$\pi$ on $\R$ if necessary, we can assume that $f(\Sing f)=\N$ because $f(\Sing f)$ has no accumulating points in $\R$. 
Define a map $\alpha_f:\N\to\N$ so that for each $n\in\N$, $f-n$ is $\alpha_f(n)$-flat at any point 
of $f^{-1}(n)\cap\Sing f$ and not $(\alpha_f(n)+1)$-flat at some point of $f^{-1}(n)\cap\Sing f$. 
If a proper analytic function $g$ with $g(\Sing g)=\N$ is $C^\omega$ R-L equivalent to $f$ then 
$\alpha_f=\alpha_g$. 
Consider all proper $C^\omega$ functions $\pi$ on $\R$ such that $\Sing\pi=\N$ and $\pi=\id$ on $\N$. 
Then the cardinality of $\{\alpha_{\pi\circ f}\}$ is of the continuum. 
Hence the cardinality of $C^\omega$ R-L equivalence classes of proper analytic functions on $M$ with 
only normal crossing singularities is of the continuum. \qed\par
{\it Proof of (iv).} 
Assume that $\dim M>1$. 
We use the idea of handle body decomposition by Morse functions (see [Mi]). 
Let $f$ be a non-negative proper $C^\infty$ function on a non-compact connected $C^\omega$ manifold 
$M$ with singularities of Morse type. 
Approximate $f$ and changing $\R$ by some $C^\omega$ diffeomorphism of $\R$, we assume that $f$ is of class 
$C^\omega$, that $f|_{\Sing f}$ is injective and $f(\Sing f)=2\N$. 
For each $k\in\N$, let $A_k$ be the union of $f^{-1}(k)\cap\Sing f$ with one point in each connected 
component of $f^{-1}(k)$ not containing points of $\Sing f$. 
Consider the 1-dimensional simplicial complex $K$ whose 0-skeleton $K^0$ is $\cup_{k\in\N}A_k$ and 
whose 1-skeleton $K^1$ consists of 1-simplexes $\overline{a b}$, for $a,b\in K^0$, such that $f(\{a,b\})=
\{k,k+1\}$ for some $k\in\N$ and there exists a connected component $C$ of $f^{-1}((k,\,k+1))$ with 
$\overline C\ni a,b$. 
Note that such a $C$ is unique because $f_{f^{-1}((2k',\,2k'+2))}:f^{-1}((2k',\,2k'+2))\to(2k',\,2k'+2)$ 
is a proper submersion for $k'\in\N$ and that conversely for each connected component $C$ of $f^{-1}((
k,\,k+1))$ there exist $a,b\in K^0$ such that $f(\{a,b\})=\{k,k+1\}$ and $\overline C\ni a,b$. 
In other words, we can identify $K^1$ with the set of connected components of $f^{-1}((k,\,k+1)):k\in\N$. 
Moreover, for $\overline{a b}\in K^1$ there exist an injective $C^\omega$ map $l_{a,b}:[0,\,1]\to M$ with 
$l_{a,b}(0)=a$, $l_{a,b}(1)=b$, $f\circ l_{a,b}(t)=f(a)\pm t$ and $\Im l_{a,b}=\Im l_{b,a}$. 
Here for $\overline{a b}\not=\overline{a' b'}$, then $\Im l_{a,b}\cap\Im l_{a',b'}=\{a\}$ if $a=a'$ or 
$a=b'$, or $\Im l_{a,b}\cap\Im l_{a',b'}=\{b\}$ if $b=a'$ or $b=b'$, and $\Im l_{a,b}\cap\Im l_{a',b'}=
\emptyset$ otherwise. 
Hence we can identify the underlying polyhedron $|K|$ with the subset $\cup_{\overline{a b}\in K^1}\Im 
l_{a,b}$ of $M$, i.e., $K$ is realized in $M$. 
Note also that there exists a unique $C^0$ retraction $r:M\to\cup_{\overline{a b}\in K^1}\Im l_{a,b}$ 
such that $f\circ r=f$. \par
We will see that each $a\in K^0$ is the end of some half-polygon in $|K|$, i.e., there exist distinct 
$a_0=a,a_1,a_2,...\in K^0$ such that $\overline{a_i a_{i+1}}\in K^1$
for $i\in\N$. 
Note that $a_i\to\infty$ (i.e., $f(a_i)\to\infty$) as $i\to\infty$. 
Since $M$ is non-compact and connected, there exists a proper $C^1$ map $l:[0,\,\infty)\to M$ such 
that $l(0)=a$. 
We can move $\Im l$ into $|K|$ by $r$ so that $\Im l$ is the underlying polyhedron of some subcomplex 
of $K$. 
If there is a 1-simplex $s$ in $K_l\overset{\text{def}}\to=K|_{\Im l}$
with an end $v$ not equal to $l(0)$ nor equal to another 1-simplex in $K_l$, then remove $s$ and $v$ from $K_l$, and repeat 
this operation as many times as possible. 
Then $K_l$ becomes a simplicial subcomplex of $K$ and $|K_l|$ is the union of a half-polygon and 
Jordan curves. 
Remove, moreover, some vertices except $l(0)$ and 1-simplexes so that $|K_l|$ is a half-polygon. 
Then we obtain an injective simplicial map $l:\tilde{\N}\to K$ with $l(0)=a$, where $\tilde{\N}=\N
\cup\{[i,\,i+1]:i\in\N\}$. 
Let $L_a$ denote all of such $l$, and let $l_a$ be such that $\min f\circ l_a=\max\{\min f\circ l:l
\in L_a\}$ and $\#(f\circ l_a)^{-1}(\min f\circ l_a)\le\#(f\circ l)^{-1}(\min f\circ l)$ 
for $l\in L_a$ with $\min f\circ l=\min f\circ l_a$. \par
Next we show that $\min f\circ l_a\to\infty$ as $a\to\infty$. 
Otherwise, there would exist distinct $a_1,a_2,...$ in $K^0$ such that
$\min f\circ l_{a_i}$ remains constant, 
say equal to $m$. 
Note that $a_i\to\infty$ as $i\to\infty$. 
Since $f^{-1}(m)$ is compact we have a subsequence of $a_1,a_2,...$ where $\Im l_{a_i}$ contain one 
point $b_0\in K^0$ with $f(b_0)=m$. 
Next, choose a subsequence so that $\Im l_{a_i}$ contain $\overline{b_0 b_1}\in K^0$ for some $b_1
\in K^0$ and $l_{a_i}(k_i+1)=b_0$ and $l_{a_i}(k_i)=b_1$ for some $k_i\in\N$. 
Repeating these arguments we obtain sequences $a_1,a_2,...$ and $b_0,b_1,...$ in $K^0$ such that 
$l_{a_i}(k_i+i)=b_0,...,l_{a_i}(k_i)=b_i$ for some $k_i\in\N,\,i=1,2,...$ 
Then $\cup_{i\in\N}\overline{b_i b_{i+1}}$ is a half-polygon. 
Fix $i$ so large that $f(b_j)>m,\ j=i,i+1,...$, and consider a polyhedron $l_{a_i}([0,\,k_i])\cup
\overline{b_i b_{i+1}}\cup\overline{b_{i+1}b_{i+2}}\cup\cdots$. Remove vertices and open 
1-simplexes from it, as in above construction of $l$, so that the polyhedron becomes a half-polygon 
starting from $a_i$. 
This half-polygon defines a new $l\in L_{a_i}$. 
Clearly $\min f\circ l\ge m=\min f\circ l_{a_i}$ for this $l$ by the definition of $l$. 
However, $\min f\circ l=\min f\circ l_{a_i}$ by the definition of $l_a$. 
Then the difference between this $l$ and $l_{a_i}$ is $\#(f\circ l)^{-1}(m)<\#(f\circ l_{a_i}
)^{-1}(m)$ since $f\circ l_{a_i}(k_i+1)=m,$ the inclusion $f\circ l_{a_i}([0,\,k_i])\supset f\circ l([0,\,k_i])$ holds
and since $f(\overline{b_ib_{i+1}}\cup\overline{b_{i+1}b_{i+2}}\cup\cdots)>m$, which contradicts the 
definition of $l_{a_i}$. 
Thus $\min f\circ l_a\to\infty$ as $a\to\infty$. \par
Note that $M$ is $C^\infty$ diffeomorphic to $M-\Im l_a$ for any $a\in K^0$ and the diffeomorphism can 
be chosen to be id outside of a small neighborhood of $\Im l_a$. 
Hence if $\Im l_a\cap\Im l_{a'}=\emptyset$ for any $a\not= a'\in\Sing f$, there exists a $C^\omega$ 
diffeomorphism $\pi:M\to M-\cup_{a\in\Sing f}\Im l_a$ such that $f\circ\pi$ is a non-singular analytic 
function on $M$. 
Consider the case where $\Im l_a\cap\Im l_{a'}\not=\emptyset$ for some $a\not= a'\in\Sing f$. 
Set $\{a_0,a_1,...\}=\Sing f$, set $X_0=Y_0=\Im l_{a_0}$ and $Z_0=\emptyset$. 
Let $i\in\N$. 
Assume by induction that we have defined subpolyhedra $X_i\supset Y_i\supset Z_i$ of $|K|$. 
If $X_i\cap\Im l_{a_{i+1}}=\emptyset$, set $X_{i+1}=X_i\cup\Im
l_{a_{i+1}}$, set $ Y_{i+1}=\Im l_{a_{i+1}}$ 
and $Z_{i+1}=\emptyset$. 
Otherwise, set $X_{i+1}=X_i\cup\Im l_{a_{i+1}}([0,\,k_{i+1}])$, define
$Z_{i+1}$ to be the closure of the unbounded 
connected component of the set of difference of the connected component of $X_i$ containing $l_{a_{i+1}}(k_
{i+1}))$ and of $l_{a_{i+1}}(k_{i+1})$, and set $Y_{i+1}=Z_{i+1}\cup\Im l_{a_{i+1}}([0,\,k_{i+1}])$, where 
$k_{i+1}=\min\{k\in\N:X_i\cap l_{a_{i+1}}([0,\,k])\not=\emptyset\}$. 
Then $X=\cup_{i\in\N}X_i$ is the underlying polyhedron of a subcomplex of $K$, and for each $i$ there 
exists a $C^\infty$ diffeomorphism $\pi_i:M-Z_i\to M-Y_i$ such that $\pi_i=\id$ on $X_i-Z_i$ and 
outside of a small neighborhood of $Y_i-Z_i$ in $M-Z_i$. 
Since $\min f\circ l_a\to\infty$ as $a\to\infty$, we see that $\cdots\circ\pi_1\circ\pi_0:M\to M$ 
is a well-defined $C^\infty$ diffeomorphism to $M-X$. 
Approximate it by a $C^\omega$ diffeomorphism $\pi:M\to M-X$. 
Then $f\circ\pi$ is the required non-singular analytic function on $M$. \par
Consider the case where $M$ is a non-compact connected Nash manifold. 
Then there exists a proper Nash function on $M$ with only
singularities of Morse type. Actually, by Theorem VI.2.1 in [S$_2$],
the manifold $M$ is Nash diffeomorphic to the interior of a compact Nash manifold with 
boundary $M'$, which is called a {\it compactification} of $M$. 
Then by using a partition of unity of class semialgebraic $C^2$ we obtain a nonnegative semialgebraic 
$C^2$ function $\phi$ on $M'$ with zero set $\partial M'$ and with only singularities of Morse type. 
Approximating the semialgebraic $C^2$ function $1/\phi$ on $M$ by a Nash function $\psi$ in the 
semialgebraic $C^2$ topology (Approximation Theorem I), we obtain the required function. 
Note that $\#\Sing\psi<\infty$ because $\Sing\psi$ is semialgebraic. 
Hence in the same way as in the analytic case, we can find a Nash function on $M$ without singularities by 
the following fact. \par 
Let $X$ be a 1-dimensional closed semialgebraic connected subset of $M$ which is a union of smooth 
curves $X_0,...,X_k$ such that any $X_i$ is closed in $M$, any $X_i$ and $X_j$ intersect transversally 
and for each $a\in X$ there exists one and only one path from $a$ to $\infty$ in $X$. 
Then $M$ and $M-X$ are Nash diffeomorphic. \par
We prove this fact as follows. 
Assume that $M=\Int M'$ for $M'$ as above. 
Then the closure $\overline X$ of $X$ in $M'$ intersects with $\partial M'$ at one point. 
By moving $X$ by a semialgebraic $C^1$ diffeomorphism of $M$ and then by a Nash diffeomorphism 
(Approximation Theorem I) we assume that $\overline X$ is smooth at $\overline X\cap\partial M'$ and 
$\overline X$ and $\partial M'$ intersect transversally. 
Let $\xi$ denote the function on $M'$ which measures distance from
$\overline X$. This function being 
semialgebraic, we approximate $\xi|_{M'-\overline X}$ by a positive Nash function $\tilde\xi$ on $M'-
\overline X$ so that $\tilde\xi(x)\to0$ as $x$ converges to a point of $\overline X$. 
Let $\epsilon>0$ be small enough. 
Then $\tilde\xi|_{\tilde\xi^{-1}((0,\,\epsilon])}:\tilde\xi^{-1}((0,\,\epsilon])\to(0,\,\epsilon]$ 
is a proper trivial Nash submersion by [Ha]. 
Hence $M'-\overline X-\tilde\xi^{-1}((0,\,\epsilon])$ and $M'-\overline X$ are semialgebraically 
$C^1$ diffeomorphic and, moreover, Nash diffeomorphic by Approximation Theorem I. 
On the other hand, $M'-\overline X-\tilde\xi^{-1}((0,\,\epsilon))$ is a compact Nash manifold with 
corners, and if we smooth the corners then $M'-\overline X-\tilde\xi^{-1}((0,\,\epsilon))$ is 
$C^\infty$ and hence Nash (Theorem VI.2.2 in [S$_2$]) diffeomorphic to $M'$ by the assumptions on $X$, 
which implies that $M-X-\tilde\xi^{-1}((0,\,\epsilon])$ and $M$ are Nash diffeomorphic. 
Therefore, $M-X=\Int(M'-\overline X)$ is Nash diffeomorphic to $M$. \qed\par
{\it Proof of (iii) for non-proper functions.}
As in the situation above, we have a non-bounded non-singular non-negative analytic function $f$ on a non-compact 
connected analytic manifold $M$. 
Let $\pi$ be a proper analytic function on $\R$ such that $\Sing\pi=\N$ and $\pi=\id$ on $\N$. 
Then $\pi\circ f(\Sing\pi\circ f)=\N$. 
Hence as in the case of proper functions, we see that the cardinality of $C^\omega$ R-L equivalence 
classes of analytic functions on $M$ with only normal crossing singularities is of the continuum. 
\qed
\endremark
\medskip

%%%%%%%%%%%%%%%%%%%%%%%%%%%%%%%%%%%%%%%%%%%%%%%%%%%%%%%%%%%%%%%%%%%%%%%%%%%%%%%%%%%%%%%%%%%
%%%%%%%%%%%%%%%%%%%%%%%%%%%%%%%%%%%%%%%%%%%%%%%%%%%%%%%%%%%%%%%%%%%%%%%%%%%%%%%%%%%%%%%%%%%
%%%%%%%%%%%%%%%%%%%%%%%%%%%%%%%%%%%%%%%%%%%%%%%%%%%%%%%%%%%%%%%%%%%%%%%%%%%%%%%%%%%%%%%%%%%
\head 
4. Reductions
\endhead
In order to prove theorems 3.1,(2) and 3.1,(3) and theorem 3.2, we
proceed to some reductions. Firstly, we reduce the analytic case to
the Nash one, secondly we reduce the non-compact Nash case to the
compact one.
\subhead
4.1. Reduction to the Nash case
\endsubhead

By the following proposition we reduce the $C^\omega$ case of theorem 3.2 to the Nash case. 

\proclaim{Proposition 4.1}
Let $M$ be a compact Nash manifold possibly with corners, and $f$ a $C^\omega$ function on $M$ with 
only normal crossing singularities. 
Then $f$ is $C^\omega$ right equivalent to some Nash function. 
\endproclaim
\remark{Remark}
If $M$ is a non-compact Nash manifold, proposition 4.1 does not hold. 
For example, consider $M=\R$ and $f(x)=\sin x$. 
\endremark
\demo{Proof of proposition 4.1}
Set $X=f^{-1}(f(\Sing f))$. 
Let $g:\tilde X\to M$ be a $C^\omega$ immersion of a compact $C^\omega$ manifold possibly with corners 
such that $\Im g=X$, $g|_{g^{-1}(\Reg X)}$ is injective and $g_x(\tilde X_x)$ 
is an analytic subset germ of $M_{g(x)}$ for each $x\in\tilde X$. 
Here we construct $\tilde X$ and $g$ locally and then paste them. 
For a connected component $C$ of $\tilde X$ there are two possible
cases to consider: either $g(C)\subset
\partial M$ or $g(C)\not\subset\partial M$. 
Assume that $g(C)\not\subset\partial M$ for any $C$. 
Then $g(\Int\tilde X)\subset\Int M$ and $g(\partial\tilde
X)\subset\partial M$, and moreover $X$ is a normal crossing 
analytic subset of $M$. 
Consider the family of all $C^\infty$ maps $g':\tilde X\to M$ with $g'(\Int\tilde X)\subset\Int M$ and 
$g'(\partial\tilde X)\subset\partial M$. 
Then 
\enddemo
\proclaim{Lemma 4.2}
Let $r\,(>0)\in\N\cup\{\infty\}$. 
Then $g$ is $C^\infty$ stable in family, in the sense that any such $C^\infty$ map $g':\tilde X\to M$ 
close to $g$ in the $C^r$ topology is $C^\infty$ R-L equivalent to $g$. 
\endproclaim
\remark{Remark}
The proof we produce below shows that lemma 4.2 holds even if $M$ is a
non-compact Nash ($C^\infty$) manifold possibly with corners, using the 
Whitney $C^r$ (strong Whitney $C^\infty$, respectively) topology. 
\endremark
\demo{Proof of lemma 4.2}
It suffices to find a $C^\infty$ diffeomorphism of $M$ which carries $\Im g$ to $\Im g'$. 
As usual, using a tubular neighborhood of $M$ in its ambient Euclidean space, the orthogonal 
projection to $M$ and a partition of unity of class $C^\infty$, we reduce the problem to the 
following local problem. \par
Assume that $M=\R^n\times[0,\,\infty)^m$ and $X=\{x_1\cdots x_l=0\}$,
for $l\le n$. Let $y_1=y_1(x)$ be 
a $C^\infty$ function on $M$ which is close to the function $x_1$ in the Whitney $C^r$ topology and 
coincides with $x_1$ outside of a neighborhood of 0. 
Then there exists a $C^\infty$ diffeomorphism $\pi$ of $M$ which is id outside of a 
neighborhood of 0 and close to $\id$ in the Whitney $C^r$ topology and carries $\{x_2\cdots x_l=0\}\cup
\{y_1(x)=0\}$ to $X$. \par
This is true since $\pi(x_1,...,x_{n+m})=(y_1(x),x_2,...,x_{n+m})$ satisfies the requirements. 
\qed
\enddemo
\demo{Continued proof of proposition 4.1}
Let $0\ll r'\ll r\in\N$.\par
Case without corners. 
Give a Nash manifold structure to $\tilde X$ (Theorem of Nash, see Theorem I.3.6 in [S$_2$]).
Let $g':\tilde X\to M$ be a Nash approximation of $g$ in the $C^r$ topology, e.g., the composite of a 
polynomial approximation of the map $g$ from $\tilde X$ to the ambient Euclidean space of $M$ with the 
orthogonal projection of a tubular neighborhood of the space, and set $X'=\Im g'$. 
Then by lemma 4.2 there exists a $C^\infty$ diffeomorphism $\pi$ of $M$ which carries $X$ to $X'$, 
and by the above proof $\pi$ can be arbitrarily close to the identity map in the $C^r$ topology. 
Let $t_1,...,t_l$ be the critical values of $f$. 
We assume that $t_i\,>0$. 
We want to construct a Nash function $f'$ on $M$ such that
$(f')^{-1}(f'(\Sing f'))=X'$ and $f'\circ
\pi=f$ on $X$ for some modified $\pi$ and, moreover, $f'\circ\pi-t_i$ and $f-t_i$ have the same 
multiplicity at each point of $f^{-1}(t_i)$ for each $i$. 
For each $t_i$, let $\Cal I_i$ denote the sheaf of $\Cal N$-ideals with zero set $\pi(f^{-1}(t_i))$ 
and having the same multiplicity as $f\circ\pi^{-1}-t_i$ at each point of $\pi(
f^{-1}(t_i))$. 
Such a sheaf exists because a non-singular
semialgebraic and analytic set germ is a 
non-singular Nash set germ. 
Then $\Cal I_i$ is generated by a finite number of global Nash functions (theorem 2.7). 
Let $\phi_i$ denote the square sum of the generators and define a Nash function $\psi_i$ on $M$ so 
that $\psi_i^2=\phi_i$ and $\psi_i$ has the
same sign as $f\circ\pi^{-1}-t_i$ everywhere. 
Note that $\psi_i^{-1}(0)=\pi(f^{-1}(t_i))$ and $\psi_i$ and $f\circ\pi^{-1}-t_i$ have the same 
multiplicity at each point of $\psi_i^{-1}(0)$. 
Set $\phi=\prod\phi_i$. 
We have a global cross-section of the sheaf of $\Cal N$-modules $\Cal N/\prod_i\Cal I^2_i$ whose value 
at each point $x$ of $\psi^{-1}_i(0)$ equals $\psi_{ix}+t_i$ mod $\Cal I^2_{ix}$. 
Apply theorem 2.8 to the homomorphism $\Cal N\to\Cal N/\prod_i\Cal I^2_i$ and the global cross-section. 
Then there exists a Nash function $\psi$ on $M$ such that $\psi-t_i$ and $f\circ\pi^{-1}-
t_i$ have the same sign at each point of a neighborhood of $\psi^{-1}_i(0)$ and the same multiplicity 
at each point of $\psi^{-1}_i(0)$ for each $i$. 
We need to modify $\psi$ so that $X'=\psi^{-1}(\psi(\Sing\psi))$. 
Let $f''$ be a $C^\infty$ function on $M$, constructed by a partition of 
unity of class $C^\infty$, such that $f''=\psi$ on a small neighborhood of $X'$ and 
$X'=(f'')^{-1}(f''(\Sing f''))$. 
Then $f''-\psi$ is of the form $\phi\xi$ for some $C^\infty$ function $\xi$ on $M$. 
Let $\tilde\xi$ be a strong Nash approximation of $\xi$ in the $C^\infty$ topology, and set $f'=\psi+\phi\tilde\xi$. 
Then $f'$ is a Nash function, $X'=(f')^{-1}(f'(\Sing f'))$ and $f'-t_i$ and $f\circ\pi^{-1}-t_i$ 
have the same multiplicity at each point of $\pi(f^{-1}(t_i))$ for each $i$. \par
By Theorem 3.1,(1) it suffices to see that $f$ is $C^\infty$ 
right equivalent to the function $h$ defined to be $f'\circ\pi$. 
Note that $h^{-1}(h(\Sing h))=X$, and $f-t_i$ and $h-t_i$ have the same multiplicity at each point of $f^{
-1}(t_i)$ for each $i$. 
Remember that $\pi$ is close to id in the $C^r$ topology. 
We can choose $f'$ so that $f$ and $h$ are close each other in the $C^r$ topology. 
Indeed, $f\circ\pi^{-1}-f'$ is of the form $\eta\prod_i\psi_i$ for some $C^\infty$ function $\eta$ 
on $M$. 
Replace $f'$ with $f'+\tilde\eta\prod_i\psi_i$ for a strong Nash approximation $\tilde\eta$ of $\eta$ in the $C^\infty$ topology. 
Then $f$ and $h$ are close. 
Hence we can reduce the problem, as usual, to the following local problem. \par
Let $M=\R^n,\ f(x)=x_1^{\alpha_1}\cdots x_l^{\alpha_l}$ and $h(x)=a(x)x_1^{\alpha_1}\cdots x_l^{
\alpha_l}$ for some $C^\infty$ function $a(x)$ on $M$ close to 1 in the Whitney $C^{r'}$ topology (lemma 2.12). 
Assume that $\alpha_1>0$. 
Then $f$ and $h$ are $C^\infty$ right equivalent through a $C^\infty$ diffeomorphism close to id in 
the Whitney $C^{r'}$ topology. \par
This is true since the $C^\infty$ diffeomorphism $\R^n\ni(x_1,...,x_n)\to(a^{1/\alpha_1}(x)x_1,x_2,...,
\allowmathbreak x_n)\in\R^n$ satisfies the requirements. 
Thus the case without corners is proved. \par
Case with corners. 
Let $M_1$ be a Nash manifold extension of $M$. 
We can assume that $M$ is the closure of the union of some connected components of $M_1-Y$ for a normal crossing 
Nash subset $Y$ of $M_1$. 
Let $U$ be an open semialgebraic neighborhood of $M$ in $M_1$ so small that $f$ is extensible to an 
analytic function $f_1$ on $U$ with only normal crossing singularities. 
Shrinking $U$ if necessary, we replace $X$ in the above proof with $X_1=f^{-1}_1(f_1(\Sing f_1))$, 
and we define a $C^\omega$ manifold $\tilde X_1$ and a $C^\omega$ immersion $g_1:\tilde X_1\to U$ in the same 
way. 
For each connected component $C$ of $\tilde X_1$ there are two
possible cases to consider: either $g_1(C)\subset 
Y$ or $g_1(C)\not\subset Y$. 
If $g_1(C)\subset Y$, then $g_1(C)$ is a Nash subset of $U$ with only normal crossing singularities and $C$ has an abstract Nash manifold structure such that $g_1|_C$ is a Nash 
diffeomorphism to $g_1(C)$ (see [S$_2$] for the definition of an abstract Nash manifold). 
Apply Artin-Mazur Theorem to $g_1(C)$. 
Then $C$ with this abstract Nash manifold structure is a Nash manifold. 
Set $g'_1=g_1$ on $C$. 
If $g_1(C)\not\subset Y$, give a Nash manifold structure to $C$, approximate $g_1|_C$ by a Nash 
immersion $g'_1|_C:C\to M_1$. 
In this way we define a Nash immersion $g'_1:\tilde X_1\to M_1$ and set $X'=\Im g'_1\cap M$. 
The rest proceeds in the same way as the case without corners. 
\qed
\enddemo
The following lemma is the $C^\omega$ or Nash version of lemma 4.2 and is used to prove theorems 3.1,(2), 3.1,(3) and lemma 
4.4. 
\proclaim{Lemma 4.3}
Let $r\,(>0)\in\N\cup\{\infty\}$. 
Let $M$ and $N$ be compact $C^\omega$ manifolds possibly with corners
such that $\dim M=1+\dim N$. Let $\phi:N\to M$ be a $C^\omega$ immersion such that $\phi(\Int N)\subset\Int M$, $\phi(\partial N)\subset
\partial M$, $\Im\phi$ is a normal crossing analytic subset of $M$ and the restriction of $\phi$ to 
$\phi^{-1}(\Reg\Im\phi)$ is injective. 
Then $\phi$ is $C^\omega$ stable in the family of $C^\omega$ maps from $N$ to $M$ carrying $\partial N$ 
to $\partial M$ in the same sense as in lemma 4.2. 
If $M,\ N$ and $\phi$ are of class Nash, then $\phi$ is Nash stable in the family of Nash maps with the same 
property as above. 
\endproclaim
\remark{Remark}
In the case of a non-compact $M$ and proper $\phi$, we see easily that the former half part of lemma 4.3 holds 
in the Whitney $C^r$ topology, $r\,(>0)\in\N\cup\{\infty\}$. 
We can prove the latter half in the non-compact case in the semialgebraic $C^r$ topology by reducing 
to the compact case by lemmas 4.5 and 4.6. 
\endremark
\demo{Proof of lemma 4.3}
Let $\psi$ be an analytic approximation of $\phi$ in family in the analytic case. 
Then by lemma 4.2 $\psi$ is $C^\infty$ R-L equivalent to $\phi$, namely, there exists a 
$C^\infty$ diffeomorphism $\pi$ of $M$ which carries $\Im\phi$ to $\Im\psi$. 
Note that we can choose $\pi$ to be close to id in the $C^r$ topology by the proof of lemma 4.2. 
Then by step 1 in the proof of theorem 3.1,(1) and its proof, we can choose an analytic $\pi$ even in the case 
with corners. 
The existence of an analytic diffeomorphism $\tau$ of $N$ with $\psi\circ\tau=\pi\circ\phi$ is clear because 
$\tau=\psi^{-1}\circ\pi\circ\phi$ on $\phi^{-1}(\Reg\Im\phi)$. 
Thus $\phi$ and $\psi$ are $C^\omega$ R-L equivalent. \par
Assume that $M,\ N,\ \phi$ and $\psi$ are of class Nash. 
It suffices to find a Nash diffeomorphism of $M$ which carries $\Im\phi$ to $\Im\psi$. 
Let $\pi$ be such a diffeomorphism of $M$ of class $C^\omega$. \par
Case without corners. 
Let $\Cal I_\phi$ and $\Cal I_\psi$ denote the sheaves of $\Cal N$-ideals on $M$ defined by $\Im
\phi$ and $\Im\psi$, respectively, and let $\{f_i\}$ and $\{g_j\}$ be a finite number of their 
respective global generators. 
Then $\{g_j\circ\pi\}$ is a set of global generators of the sheaf of $\Cal O$-ideals $\Cal I_\phi\Cal O$ 
on $M$. 
Hence in the same way as in step 2 of the proof of theorem 3.1,(1) we obtain $C^\omega$ functions $\alpha_{i,j
}$ and $\beta_{i,j}$ on $M$ such that 
$$
f_i=\sum_j\alpha_{i,j}\cdot(g_j\circ\pi)\quad\text{and}\quad g_j\circ\pi=\sum_i\beta_{i,j}f_i.
$$
Let $M\subset\R^n$ and let $h$ be a Nash function on $\R^n$ with zero set $M$. 
Extend $g_j$ to Nash functions on $\R^n$ and use the same notation $g_j$ (theorem 2.8). 
Consider the following equations of Nash functions in variables $(x,y,a_{i,j},b_{i,j})\in M\times\R^n\times\R^
{n'}\times\R^{n'}$, for some $n'$:
$$
h(y)=0,\ \ f_i(x)-\sum_ja_{i,j}g_j(y)=0,\ \ g_j(y)-\sum_ib_{i,j}f_i(x)=0.
$$
We have a $C^\omega$ solution $y=\pi(x),\ a_{i,j}=\alpha_{i,j}(x)$ and $b_{i,j}=\beta_{i,j}(x)$. 
Hence by Nash Approximation Theorem II there exists a Nash solution $y=\pi'(x),\ a_{i,j}=\alpha'_{i,j}(x)$ 
and $b_{i,j}=\beta'_{i,j}(x)$, which are close to $\pi,\ \alpha_{i,j}$ and $\beta_{i,j}$, 
respectively. 
Then 
$$
\Im\pi'=M,\ \ f_i=\sum_j\alpha'_{i,j}\cdot(g_j\circ\pi'),\ \ g_j\circ\pi'=\sum_i\beta_{i,j}f_i.
$$
Hence $\pi'$ is a Nash diffeomorphism of $M$ and carries $\Im\phi$ to $\Im\psi$. \par
Case with corners. 
We can assume that for Nash manifold extensions $M_1$ and $N_1$ of $M$ and $N$, respectively, $\phi,
\psi$ are extensible to Nash immersion $\phi_1$ and $\psi_1$
of $N_1$ into $M_1$ and $\pi$ to a 
$C^\omega$ embedding $\pi_1$ of a semialgebraic open neighborhood $U$ of $M$ in $M_1$ into $M_1$, 
and moreover there exist normal crossing Nash subsets $Y$ of $M_1$ and $Z$ of $N_1$ such 
that $M$ and $N$ are closures of the unions of some connected components of $M_1-Y$ and of $N_1-Z$, 
respectively, $\phi_1(Z)\subset Y$ and $\psi_1(Z)\subset Y$. 
Let $M_1$ be contained and closed in an open semialgebraic set $O$ in
$\R^n$, and $h_1$ be a Nash 
function on $O$ with zero set $M_1$. 
Take a small open semialgebraic neighborhood $V$ of $M$ in $M_1$ and shrink $M_1,N_1$ and $U$ so 
that $\pi(U)\subset V$ and $U\cap\Im\phi_1$ and $V\cap\Im\psi_1$ are normal crossing Nash subsets 
of $U$ and $V$, respectively. 
Then in the same way as above, we obtain a finite number of global generators $\{f_{1,i}\}$ and 
$\{g_{1,i}\}$ of the sheaves of $\Cal N$-ideals on $U$ and $V$ defined by $U\cap\Im\phi_1$ and 
$V\cap\Im\psi_1$, respectively, and analytic functions $\alpha_{1,i,j},\beta_{1,i,j}$ on $U$ such 
that 
$$
f_{1,i}=\sum_j\alpha_{1,i,j}\cdot(g_{1,j}\circ\pi_1)\quad\text{and}\quad g_{1,j}\circ\pi_1=\sum_i
\beta_{1,i,j}f_{1.i}\quad\text{on}\ U.
$$
We need to describe the condition $\pi(\partial M)=\partial M$, i.e., $\pi_1(U\cap Y)\subset Y$. 
Let $\xi'$ be the square sum of a finite number of global generators of the sheaf of $\Cal N$-ideals 
$\Cal I$ on $M_1$ defined by $Y$. 
Then $\xi'$ is a generator of $\Cal I^2$, and since $M$ is a manifold with corners there exists 
a unique Nash function $\xi$ on a semialgebraic neighborhood of $M$ in $M_1$ such that $\xi^2=
\xi'$ and $\xi>0$ on $\Int M$. 
Shrink $U$ once more. 
Then $\xi|_U$ and $\xi\circ\pi_1|_U$ are well-defined generators of $\Cal I|_U$ and $\Cal{IO}|_U$, 
respectively, and we have a positive analytic function $\gamma$ on $U$ such that $\xi\circ\pi=\gamma
\xi$ on $U$. 
We shrink $O$, and using the same notation we extend $g_{1,j}$ and $\xi$ to Nash functions on $O$. 
\par
Consider the germs on $M\times O\times\R^{n'}\times\R^{n'}\times\R$ of
the following equations of Nash functions 
in the variables $(x,y,a_{i,j},b_{i,j},c)\in U\times O\times\R^{n'}\times\R^{n'}\times\R$ for some $n'$: 
$$
h_1(y)=0,\ f_{1,i}(x)-\sum_ja_{i,j}g_{1,j}(y)=0,\ g_{1,j}(y)-\sum_ib_{i,j}f_{1,i}(x)=0,\ \xi(y)-c\xi(
x)=0.
$$
Then, since Nash Approximation Theorem II holds in the case of germs, we have Nash germ solutions on $M$ of 
the equations $y=\pi'_1(x),\ a_{i,j}=\alpha'_{i,j}(x),\ b_{i,j}=\beta'_{i,j}(x)$ and $c=\gamma'(x)$. 
The equation $\xi\circ\pi'_1=\gamma'\xi$ means $\pi'_1(M)=M$. 
Thus $\pi'_1|_M$ is the required Nash diffeomorphism of $M$. 
\qed
\enddemo
The following lemma shows countable cardinality of the normal crossing Nash $(C^\omega)$ subsets of a 
compact Nash ($C^\omega$, respectively,) manifold possibly with corners. 
\proclaim{Lemma 4.4}
Let $M$ be a compact Nash manifold possibly with corners of strictly positive dimension. 
Consider Nash immersions $\phi$ from compact Nash manifolds possibly with corners of dimension equal to 
$\dim M-1$ into $M$ such that $\Im\phi$ are normal crossing Nash
subsets of $M$, the restrictions $\phi|_{\phi^{-1}(\Reg
\Im\phi)}$ are injective and $\phi$ carry the interior and the corners into the interior and the corners, 
respectively. 
Then the cardinality of Nash R-L equivalence classes of all the $\phi$'s is countable. \par
The analytic case also holds. 
\endproclaim
\demo{Proof of lemma 4.4}
Note that the cardinality is infinite because for any $k\in\N$ we can embed $k$ copies of a sphere of 
dimension $\dim M-1$ in $M$. 
It suffices to treat only the Nash case for the following reason. \par
Let $\phi:M'\to M$ be an analytic immersion as in lemma 4.4 for analytic $M'$ and $M$. 
Assume that $M$ has no corners. 
Since a compact analytic manifold is $C^\omega$ diffeomorphic to a Nash manifold, we suppose that $M'$ 
and $M$ are Nash manifolds. 
Approximate $\phi$ by a Nash map. 
Then $\phi$ is $C^\omega$ R-L equivalent to the approximation by lemma 4.3. 
Hence we can replace $\phi$ by a Nash map. \par
Assume that $M$ has corners. 
Let $M_1\subset\R^n$ be an analytic manifold extension of $M$ such that $M$ is the closure of the union of 
some connected components of $M_1-Y$ for a normal crossing analytic subset $Y$ of $M_1$. 
We can assume that $M_1$ is compact as follows. 
Let $\alpha$ denote the function on $M_1$ which measures distance from $M$. 
Approximate $\alpha|_{M_1-M}$ by a $C^\omega$ function $\alpha'$ in the Whitney $C^\infty$ topology, 
and let $\epsilon<\epsilon'$ be positive numbers so small that $M\cup\alpha^{\prime-1}((0,\,\epsilon'
])$ is compact and such that the restrictions of $\alpha'$ to $\alpha^{\prime-1}((\epsilon,\,\epsilon'
))$ and to its intersections with strata of the canonical stratification of $Y$ are regular. 
Then $(M_1\cap\alpha^{\prime-1}((\epsilon,\,\epsilon')),Y\cap\alpha^{\prime-1}((\epsilon,\,\epsilon')))$ 
is $C^\omega$ diffeomorphic to $(M_1\cap\alpha^{\prime-1}((\epsilon+\epsilon')/2)),Y\cap\alpha^{\prime
-1}((\epsilon+\epsilon')/2)))\times(\epsilon,\,\epsilon')$. 
Hence, replacing $M_1$ with the double of $M\cup\alpha^{\prime-1}((0,\,(\epsilon+\epsilon')/2])$, we 
assume that $M_1$ is compact. \par
Next we reduce the problem to the case where $M_1$ and $Y$ are of class Nash. 
Define, as in the proof of proposition 4.1, a $C^\omega$ immersion $g:\tilde Y\to M_1$ of a compact $C^\omega$ 
manifold so that $\Im g=Y$, so that $g|_{g^{-1}(\Reg\Im g)}$ is injective and $g_y(
\tilde Y_y)$ is an analytic subset germ of $M_{1g(y)}$ for each
$y\in\tilde Y$. Give Nash structures on $M_1$ and $\tilde Y$, 
and approximate $g$ by a Nash map $g'$. 
Then by lemma 4.3 there exists a $C^\omega$ diffeomorphism $\pi$ of $M_1$ which carries $\Im g$ to $\Im 
g'$. 
Hence we can replace $Y$ with $\Im g'$ and we assume from the beginning that $M_1,\, Y$ and $M$ are of class 
Nash. 
By the same reason, we suppose that $M'$ is a Nash manifold possibly with corners and the closure of the union of 
some connected components of $M'_1-Y'$ for a compact Nash manifold extension $M'_1$ of $M'$ and a normal 
crossing Nash subset $Y'$ of $M'_1$. 
Extend $\phi$ to a $C^\omega$ immersion $\phi_1$ of a compact semialgebraic neighborhood $U$ of $M'$ in 
$M'_1$ into $M_1$, choose $U$ so small that $\phi_1(U\cap Y')\subset Y$, and approximate, as in the proof 
of step 1 in theorem 3.1,(1), $\phi_1$ by a Nash map $\tilde\phi_1$ so that $\tilde\phi_1(U\cap Y')
\subset Y$ (here we use theorems 2.7 and 2.8 in place of corollaries 2.2 and 2.4 in the proof in theorem 3.1,(1)). 
Then $\tilde\phi_1|_{M'}$ is a Nash immersion into $M$, and $\Im\tilde\phi_1|_{M'}$ is a normal crossing Nash 
subset of $M$, moreover
$\tilde\phi_1|_{(\tilde\phi_1|_{M'})^{-1}(\Reg\Im\tilde\phi_1|_{M'})}$
is injective and $\tilde\phi_1(\partial M')\subset\partial M$, and
finally $\tilde\phi_1|_{M'}$ is $C^\omega$ R-L equivalent to 
$\phi$ by lemma 4.3. 
Thus we reduce the analytic case to the Nash one. \par
Consider the Nash case. 
Let $M\subset\R^n$ and $\phi:M'\to M$ be a Nash immersion as in lemma 4.4. 
Let $M_1,\,M'_1,\,Y$ and $\phi_1:M'_1\to M_1$ be Nash manifold extensions of $M$ and $M'$, a normal 
crossing Nash subset of $M_1$ and a Nash immersion, respectively, such that $M_1\subset\R^n,\ \phi_1=
\phi$ on $M'$, $M$ and $M'$ are the closures of the unions of some connected components of $M_1-Y$ and $M'
_1-\phi_1^{-1}(Y)$, respectively, $U\cap\Im\phi_1$ is a normal crossing Nash subset of an open semialgebraic neighborhood $U$ of $M$ in $M_1$ and $\phi_1|_{\phi_1^{-1}(\Reg (U\cap\Im\phi_1))}$ is injective. 
By Artin-Mazur Theorem (see the proof of theorem 2.9) we can regard $M'_1$ as an open semialgebraic 
subset of the regular point set of an algebraic variety in $\R^n\times\R^{n'}$ for some $n'$ and 
$\phi_1$ as the restriction to $M'_1$ of the projection $\R^n\times\R^{n'}\to\R^n$. 
We will describe all such $\phi_1:M'_1\to M_1$ with fixed complexity as follows. 
Any algebraic set in $\R^n\times\R^{n'}$, and its subset of regular points where the projection to 
$\R^n$ is regular, are, respectively, described by the common zero set of polynomial functions $f_1,...,
f_l$ on $\R^n\times\R^{n'}$ for some $l\in\N$ and 
$$
\gather
\bigcup
\Sb I=\{i_1,...,i_k\}\subset\{1,...,l\}\\
I'=\{i'_1,...,i'_{k+n'}\}\subset\{1,...,l\}\endSb
\{x=(x_1,...,x_{n+n'})\in\R^{n+n'}:\qquad\qquad\qquad\qquad\qquad\qquad\\
\ f_1(x)=\cdots=f_l(x)=0,\quad |\frac{\partial(f_{i'_1},...,f_{i'_{k+n'}})}{\partial(x_{i_1},...,x_{i_k},
x_{n+1},...,x_{n+n'})}(x)|\not=0,\\
\qquad\qquad\qquad g_{I',i''}f_{i''}=\sum_{j=1}^{k+n'}g_{I',i'',j}f_{i'_j},\quad g_{I',i''}(x)\not=0,\ 
\ i''\in\{1,...,l\}-I'\}
\endgather
$$
for some polynomial functions $g_{I',i''}$ and $g_{I',i'',j}$ on $\R^n\times\R^{n'}$, where $k=n+1-\dim 
M$, and $\frac{\partial(\ )}{\partial(\ )}$ denotes the Jacobian matrix. 
Moreover, its open semialgebraic subset is its intersection with 
$$
\cup_{j'=1}^l\cap_{i=1}^l\{x\in\R^{n+n'}:h_{i,j'}(x)>0\} 
$$
for some polynomial functions $h_{i,j'}$ on $\R^n\times\R^{n'}$ (here we enlarge $l$ if necessary). 
Thus $\phi_1:M'_1\to\R^n$ is described by the family $f_i,\,g_{I',i''},\,g_{I',i'',j}$ and $h_{i,j'}$ and 
conversely, any polynomial functions $f_i,\,g_{I',i''},\,g_{I',i'',j}$ and $h_{i,j'}$ define in the above 
way a Nash manifold $M'_1$ in $\R^{n+n'}$ of dimension $\dim M-1$ such that the projection $\phi_1:M'_1
\to\R^n$ is an immersion. 
If the degree of the polynomials are less than or equal to $d\in\N$,
we say that $\phi_1:M'_1\to\R^n$ is of 
complexity $l,d,n'$. \par 
Furthermore, since a polynomial function on $\R^{n+n'}$ of degree less than or 
equal to $d$ is of the form $\sum_{\alpha\in\N_d^{n+n'}}a_\alpha x^\alpha,\ a_\alpha\in\R$, where 
$\N_d^{n+n'}=\{\alpha\in\N^{n+n'}:\ |\alpha|\le d\}$, we regard the family of $f_i,..,h_{i,j'}$ of degree less than or 
equal to $d$ as an element $a=(a_\alpha)$ of $\R^N$ for some $N\in\N$. 
We write $\phi_1:M'_1\to\R^n$ as $\phi_{1a}:M'_{1a}\to\R^n$. 
Then the set $X=\cup_{a\in\R^N}\{a\}\times M'_{1a}\subset\R^N\times\R^n\times\R^{n'}$ is semialgebraic, 
and we can identify $\phi_{1a}:M'_{1a}\to\R^n$ with $p|_{(q\circ p)^{-1}(a)}:(q\circ p)^{-1}(a)\to\{a\}
\times\R^n$, where $p:X\to\R^N\times\R^n$ and $q:\R^N\times\R^n\to\R^N$ are the projections. \par
Consider the condition $\Im\phi_{1a}\subset M_1$. 
The subset of $\R^{N}$ consisting of $a$ such that $p|_{(q\circ p)^{-1}(a)}$ fails to satisfy this condition 
is $q\circ p(X\cap\R^{N}\times(\R^n-M_1)\times\R^{n'})$ and hence is semialgebraic. 
Let $A$ denote its complement in $\R^{N}$. 
Thus $\Im\phi_{1a}\subset M_1$ if and only if $a\in A$. \par
Next consider when $U\cap\Im\phi_{1a}$ is normal crossing and $\phi_{1a}|_{\phi_{1a}^{-1}(\Reg(U\cap\Im
\phi_{1a}))}$ is injective. 
For that, remember that the tangent space $T_x M'_{1a}$ of $M'_{1a}$ at $x\in M'_{1a}$, for $M'_{1a}$ 
described by $f_i,\,g_{I',i''},...$ as above, is given by 
$$
T_x M'_{1a}=\{y\in\R^{n+n'}:d f_{1x}y=\cdots=d f_{lx}y=0\},
$$
and hence the set $TX$ defined to be $\{(a,x,y)\in X\times\R^{n+n'}:y\in T_xM'_{1a}\}$ is 
semialgebraic. 
Assume that $a\in A$. 
Set 
$$
\gather
M''_{1a}=\{(x,x')\in M'_{1a}\times M'_{1a}:x\not= x',\ \phi_{1a}(x)=\phi_{1a}(x')\in U,\qquad\qquad\qquad\\
\qquad\qquad\qquad\dim(d\phi_{1ax}(T_x M'_{1a})+d\phi_{1ax'}(T_{x'}M'_{1a}))=\dim M-1\}.
\endgather
$$
Then $M''_{1a}$ and $\cup_{a\in A}\{a\}\times M''_{1a}$ are semialgebraic, and $a\in A-q'(\cup_{a\in A}\{a
\}\times M''_{1a})$ if and only if for any $x\not=x'\in M'_{1a}$ with $\phi_{1a}(x)\in U$, the germs of 
$\phi_{1a}$ at $x$ and $x'$ intersect transversally, where $q':\R^N\times\R^n\times\R^{n'}\times\R^n\times
\R^{n'}\to\R^N$ denotes the projection. 
Repeating the same arguments on $m$-tuple of $M'_{1a}$ for any $m\le\dim M$ we obtain a semialgebraic subset 
$B$ of $A$ such that for $a\in A$, then $a\in B$ if and only if $U\cap\Im\phi_{1a}$ is normal crossing and 
$\phi_{1a}|_{\phi_{1a}^{-1}(\Reg(U\cap\Im\phi_{1a}))}$ is injective. \par
Let $\{B_i\}$ be a finite stratification of $B$ into connected Nash manifolds such that $q\circ p$ is 
Nash trivial over each $B_i$ [C-S$_1$], i.e., for each $i$ there exists a Nash diffeomorphism $\pi_i:(q
\circ p)^{-1}(B_i)\to(q\circ p)^{-1}(b_i)\times B_i$ of the form $\pi_i=(\pi'_i,q\circ p)$ for some $b_i
\in B_i$. 
For $a\in B$, set $M'_a=\phi_{1a}^{-1}(M)$ and $\phi_a=\phi_{1a}|_{M'_a}$. 
Then $\{\phi_a:M'_a\to\R^n\}_{a\in B}$ is the family of all $\phi:M'\to M$ as in lemma 4.4 which are 
extensible to $\phi_1:M'_1\to M_1$ with fixed $U$ and complexity $l,\,d,\,n'$, and if $a$ and $a'$ are in 
the same $B_i,\ i\in I$, then $\phi_a:M'_a\to M$ and $\phi_{a'}:M'_{a'}\to M$ are Nash R-L equivalent by 
lemma 4.3 for the following reason. 
As there exists a $C^0$ curve in $B_i$ joining $a$ and $a'$, considering a finite sequence of points on 
the curve we can assume that $a'$ is close to $a$ as elements of $\R^N$. 
We can replace $\phi_a$ and $\phi_{a'}$ with $\phi_a\circ(\pi'_i|_{\{a\}\times M'_a})^{-1}=p_n\circ(\pi'
_i|_{\{a\}\times M'_a})^{-1}:\{b_i\}\times M'_{b_i}\to\R^n$ and $p_n\circ(\pi'_i|_{\{a'\}\times M'_{a'}}
)^{-1}:\{b_i\}\times M'_{b_i}\to\R^n$, where $p_n$ denotes the projection $\R^N\times\R^n\times\R^{n'}
\to\R^n$. 
Hence in order to apply lemma 4.3 it suffices to see that $(\pi'_i|_{\{a'\}\times M'_{a'}})^{-1}$ is close 
$(\pi'_i|_{\{a\}\times M'_a})^{-1}$ in the $C^1$ topology. 
That is true because we can regard $(\pi'_i|_{\{a'\}\times M'_{a'}})^{-1}$ and $(\pi'_i|_{\{a\}\times 
M'_a})^{-1}$ as $\pi_i^{-1}|_{(q\circ p)^{-1}(b_i)\times\{a'\}}$ and $\pi_i^{-1}|_{(q\circ p)^{-1}(b_i)
\times\{a\}}$, respectively, $(q\circ p)^{-1}(b_i)$ is compact and because of the following fact. 
For compact $C^1$ manifolds $M_2$ and $M_3$ and for a $C^1$ function $g:M_2\times M_3\to\R$ if two points 
$u$ and $v$ in $M_3$ are close each other then the functions $M_2\ni x\to g(x,u)\in\R$ and $M_2\ni x\to 
g(x,v)\in\R$ are close in the $C^1$ topology. 
Hence the cardinality of equivalence classes of $\phi_a:M'_a\to M,\ a\in B$, is finite. 
Until now we have fixed $U$. 
We need argue for all semialgebraic neighborhood $U$ of $M$ in $M_1$. 
However, it is sufficient to treat a countable number of $U$'s since $M$ is compact. 
Thus the cardinality of all equivalence classes is countable. 
\qed
\enddemo

\subhead
4.2. Compactification of a Nash function with only normal crossing singularities
\endsubhead

The following lemmas 4.5, 4.6, 4.7 and proposition 4.8 are
preparations for proposition 4.9 that states the compactification of a
Nash function with only normal crossing singularities. The main tools
are Nash sheaf theory and the Nash version of Hironaka desingularization
theorems.\par
Lemmas 4.5 and 4.6 show that a normal crossing Nash subset of a non-compact Nash manifold is trivial at 
infinity.

\proclaim{Lemma 4.5}
Let $X$ be a normal crossing Nash subset of a Nash manifold $M$ and $f:M\to\R^m$ a proper Nash 
map whose restrictions to $M-X$ and to strata of the canonical stratification of $X$ are submersions onto 
$\R^m$. 
Then $f$ is Nash trivial, i.e., there exists a Nash diffeomorphism $\pi:M\to f^{-1}(0)\times\R^m$ of 
the form $\pi=(\pi',f)$, and $\pi'$ can be chosen so that $\pi'(X)=X\cap f^{-1}(0)$ and $\pi'=\id$ on $f^{-1
}(0)$. \par
The analytic case also holds. 
\endproclaim
This is shown in [C-S$_{1,2}$] in the case of empty $X$. 
We prove here the nonempty case. 

\demo{Proof of lemma 4.5}
Consider the Nash case. 
Let $n=\dim M$, take $k$ an integer smaller than $n$, and let $X_k$ denote the union of strata of the 
canonical stratification of $X$ of dimension less than or equal to $k$. 
We define $\pi'$ on $X_k$ by induction on $k$, and then on $M$. 
To this aim, we can assume that $\pi'$ is already given on $X$, say $\pi_X=(\pi'_X,f_X)$, by the following fact, 
where $f_A=f|_A$ for a subset $A$ of $M$. \par
{\it Fact 1.}
There exist a Nash manifold $\tilde X_k$ of dimension $k$ and a Nash immersion $p_{\tilde X_k}:\tilde X
_k\to M$ such that $\Im p_{\tilde X_k}=X_k$ and $p_{\tilde X_k}|_{p^{-1}_{\tilde X_k}(X_k-X_{k-1})}$ 
is injective. \par
{\it Proof of fact 1.} (Artin-Mazur Theorem. See the proof of theorem 2.9.)
Let $M$ be contained and closed in $\R^N$, and let $X^Z_k$ denote the Zariski closure of $X_k$ in $\R^N$. 
Then there exist an algebraic variety $\widetilde{X^Z_k}$ (the normalization of $X^Z_k$) in $\R^N\times
\R^{N'}$ for some $N'\in\N$ and the union of some connected components $\tilde X_k$ of $\widetilde{X^Z_k}$ 
such that $\widetilde{X_k^Z}$ is non-singular at $\tilde X_k$. Hence $\tilde X_k$ is a Nash manifold and 
the restriction $p_{\tilde X_k}$ to $\tilde X_k$ of the projection $p:\R^N\times\R^{N'}\to\R^N$ 
satisfies the requirements in fact 1. \par
Let $\phi_i$ be a finite number of global generators of the sheaf of $\Cal N$-ideals $\Cal I$ on $M$ 
defined by $X$, and set $\phi=\sum\phi_i^2$. 
Then $\phi>0$ on $M-X$ and $\phi$ is a global generator of $\Cal I^2$. 
For a subset $A$ of $M$ and $x\in\R^m$, set $A(x)=A\cap f^{-1}(x)$. 
We will extend $\pi'_X$ to $\pi'$. 
For that it suffices to find $\pi'$ of class semialgebraic $C^l$ for a large integer $l$ for the 
following reason. \par
{\it Fact 2.}
Let $g$ be a semialgebraic $C^l$ function on $M$ whose restriction to $X$ is of class Nash. 
Then fixing $g$ on $X$ we can approximate $g$ by a Nash function in the semialgebraic $C^{l-n}$ 
topology. \par
{\it Proof of fact 2.} 
As in the proof of theorem 3.1,(1), step 1, $g|_X$ is extensible to a Nash function $G$ on $M$ by theorem 
2.8. 
Replace $g$ with $g-G$. 
Then we can assume that $g=0$ on $X$ and $g$ is of the form $\sum g_i\phi_i$ for some semialgebraic $C^{l-n}$ 
functions $g_i$ on $M$ for the following reason. 
Reduce the problem to the case where $(M,X)=(\R^n,\{x_1\cdots x_{n'}=0\})$ and $\{\phi_i\}=\{x_1\cdots x_{n'}\}$ for some $n'\in\N$ by a partition of unity of class semialgebraic $C^l$ (remark 2.11,(5)$'$). 
Then $g$ is divisible by $x_1\cdots x_{n'}$ as a semialgebraic $C^{l-n}$ function on $M$ by elementary calculations. 
Hence $g$ is of the form $g_1x_1\cdots x_{n'}$ for some semialgebraic $C^{l-n}$ function $g_1$ on $M$. 
As usual, we approximate $g_i$ by Nash functions $\tilde g_i$ in the semialgebraic $C^{l-n}$ topology we 
obtain the required approximation $\sum\tilde g_i\phi_i$ of $g$ in fact 2. \par
We will see that there exists a finite semialgebraic $C^l$ stratification $\{B_i\}$ of 
$\R^m$ such that for each $i$, the map $\pi_X|_{X\cap f^{-1}(B_i)}$ is extensible to a semialgebraic $C^l$ 
diffeomorphism $\pi_i=(\pi'_i,f_{f^{-1}(B_i)}):f^{-1}(B_i)\to M(b_i)\times B_i$ for some point $b_i\in 
B_i$. 
In the following arguments we need to stratify $\R^m$ into $\{B_i\}$, each $B_i$ into $\{B_{i,j}:j=1,2,...\}$ 
and once more. 
However, we always use notation $\R^m$ for all strata for simplicity of notation, which does not cause problems because we can choose 
stratifications so that strata are semialgebraically $C^l$ diffeomorphic to Euclidean spaces. \par
  We recall the construction of $\pi$ as in the proof of Theorem II.6.7 in [S$_3$]. 
Without loss of generality we assume that $\pi'_X|_{X(0)}=\id$. 
First we can modify in order $\phi$ to be a semialgebraic $C^l$ function so that for each $x\in\R^m$, 
$\phi|_{M(x)-X}$ has only singularities of Morse type (Claim 2, ibid.) (here we need to stratify $\R^m$ and consider the restriction of $\phi$ to each stratum in place of $\phi$, 
and the main method of proof is a semialgebraic version of Thom's transversality theorem), $\phi$ is 
constant on each connected component of $Z\overset\text{def}\to=\cup_{x\in\R^m}\Sing(\phi|_{M(x)-X})$ 
and the values are distinct from each other (Claim 4, ibid.) after the second stratification. 
Next, let $Y$ be a connected component of $Z$ and set $\tilde Y=\phi^{-1}(\phi(Y))$. 
Then there exist an open semialgebraic neighborhood $U$ of $\tilde Y$ in $M$ and a semialgebraic 
$C^l$ embedding $u=(u',f_U):U\to U(0)\times\R^m$ such that $u'|_{U(0)}=\id$ and $\phi\circ u'=\phi|_U$ 
(Claim 5, ibid.) after the third stratification. 
Thirdly, applying lemma 4.5 without $X$ to the semialgebraic $C^l$ map $(f,\phi)|_{\phi^{-1}(I)}:\phi^
{-1}(I)\to\R^m\times I$ for each connected component $I$ of $(0,\,\infty)-\phi(Z)$, we obtain a 
semialgebraic $C^l$ diffeomorphism $\lambda=(\lambda',f_{\phi^{-1}(I)},\phi|_{\phi^{-1}(I)}):\phi^{-1}
(I)\to\phi^{-1}(I)(0,0)\times\R^m\times I$ such that $\lambda'|_{\phi^{-1}(I)(0,0)}=\id$, where $\phi
^{-1}(I)(0,0)=\phi^{-1}(t_0)\cap M(0)$ for some $t_0\in I$ (Claim 7, ibid.). 
Fourthly, we paste $u$ and $\lambda$ for all $I$ and construct a semialgebraic $C^l$ diffeomorphism 
$v=(v',f_{M-X}):M-X\to(M(0)-X)\times\R^m$ such that $v'|_{M(0)-X}=\id$ and $\phi\circ v'=\phi|_{M-X}
$, ibid. 
Hence it suffices to prove the following fact by the same idea of pasting. \par
{\it Fact 3.}
There exist an open semialgebraic neighborhood $W$ of $X$ in $M$ and a semialgebraic $C^l$ embedding 
$w=(w',f_W):W\to M(0)\times\R^m$ such that $w'=\pi'_X$ on $X$, so that $w'|_{W(0)}=\id$ and $\phi\circ w'=\phi
|_W$. \par
{\it Proof of fact 3.} 
Here the condition $\phi\circ w'=\phi|_W$ is not necessary. 
If there exists a semialgebraic $C^l$ embedding $w$ without this condition, we replace $\phi$ on $W$ 
with $\phi\circ w'$, extend it to a semialgebraic $C^l$ function on $M$ positive on $M-X$, and 
repeat the above arguments from the beginning. 
Then fact 3 is satisfied by this $w'$. \par
If $X$ is smooth, the problem becomes easier. 
Hence we reduce to the smooth case. 
Let $\tilde X\subset\R^N\times\R^{N'}$ and $p_{\tilde X}:\tilde X\to M$ be a Nash manifold and the 
restriction to $\tilde X$ of the projection $p:\R^N\times\R^{N'}\to\R^N$ defined in the proof of fact 1 
for $k=n-1$. 
For a small positive semialgebraic $C^0$ function $\epsilon$ on $\tilde X$, let $\tilde Q$ denote the 
$\epsilon$-neighborhood of $\tilde X$ in $M\times\R^{N'}$, i.e., 
$$
\tilde Q=\bigcup_{z\in\tilde X}\{z'\in M\times\R^{N'}:\dis(z,z')<\epsilon(z)\},
$$
and let $\tilde q:\tilde Q\to\tilde X$ denote the orthogonal projection, which is a Nash submersion. 
Set 
$$
\tilde M=\{(x,y)\in\tilde Q\subset M\times\R^{N'}:\tilde q(x,y)=(x',y)\ \text{for some }x'\in X\}. 
$$
Then $\tilde M$ is a Nash manifold of dimension $n$ containing $\tilde X$, and $p_{\tilde M}:\tilde M
\to M$ is a (not necessarily proper) Nash immersion, where $p_A$ denotes $p|_A$ for a subset $A$ of $M\times\R^
{N'}$. 
Set $A(0)=A\cap M(0)\times\R^{N'},$ set $f_A=f\circ p_A$ for the same $A$, and $\Tilde{\Tilde X}=p_{\tilde 
M}^{-1}(X)$. 
Then $\Tilde{\Tilde X}$ is a normal crossing Nash subset of $\tilde
M$, and $p_{\Tilde{\Tilde X}}:\Tilde
{\Tilde X}\to X$ is a (not necessarily proper) local Nash diffeomorphism at each point of $\Tilde{
\Tilde X}$. Moreover $\pi_X=(\pi'_X,f_X)$ is lifted to $\pi_{\Tilde{\Tilde X}}=(\pi'_{\Tilde{\Tilde X}},f_{
\Tilde{\Tilde X}}):\Tilde{\Tilde X}\to\Tilde{\Tilde X}(0)\times\R^m$, and there exists a Nash 
function $\tilde\phi$ on $\tilde M$ with zero set $\tilde X$ which is, locally at each point of 
$\tilde X$, the square of a regular function and such that 
$$
\tilde\phi=\tilde\phi\circ\pi'_{\Tilde{\Tilde X}}\quad\text{on}\ \Tilde{\Tilde X}. \tag$1$
$$
To be precise, we construct $\tilde\phi$ first on $\tilde M(0)$, and extend it to $\Tilde{\Tilde X}$ so that $(1)$ 
is satisfied and then to $\tilde M$ as usual. 
Moreover $\pi'_{\Tilde{\tilde X}}=\id$ on $\tilde{\tilde X}(0)$. \par
Note that $X_{m-1}=\emptyset$ since $f_{X_k-X_{k-1}}$ is a submersion onto $\R^m$ if $X_k\not=\emptyset$. 
Let $m\le k<n$. 
Then by the definition of $\tilde X$, the map $p_{\tilde X\cap p^{-1}(X_k-X_{k-1})}:\tilde X\cap p^{-1}(X_k-X_
{k-1})\to X_k-X_{k-1}$ is a Nash $(n-k)$-fold covering. 
Hence considering a semialgebraic triangulation of $X_k(0)$ compatible with $X_{k-1}$---a 
semialgebraic homeomorphism from the underlying polyhedron of some simplicial complex to $X_k(0)$ 
such that $X_{k-1}(0)$ is the image of the union of some simplexes--- and small open 
semialgebraic neighborhoods of inverse image of open simplexes by $\pi^{\prime-1}_X$ in $M-X_{k-1}$, we obtain finite 
open semialgebraic coverings $\{Q_{k,i}:i\}$ of $X_k-X_{k-1}$ in $M-X_{k-1}$ and $\{\tilde Q_{k,i,j}
:i,\,1\le j\le n-k\}$ of $\tilde X\cap p^{-1}(X_k-X_{k-1})$ in $\tilde M-\tilde X\cap p^{-1}(X_{k-1})$ 
such that $\pi_X^{\prime-1}(X(0)\cap Q_{k,i})=X\cap Q_{k,i}$, $(Q_{k,i},X_k\cap Q_{k,i})$ are Nash 
diffeomorphic to $(\R^n,\{0\}\times\R^k)$, such that $p_{\tilde Q_{k,i,j}}:(\tilde Q_{k,i,j},\tilde X\cap p^{-1}
(X_k)\cap\tilde Q_{k,i,j})\to(Q_{k,i},X_k\cap Q_{k,i})$ are Nash diffeomorphisms, and $\tilde Q_{k,i,j
}\cap\tilde Q_{k,i,j'}=\emptyset$ if $j\not=j'$. 
Define Nash functions $\phi_{k,i,j}$ on $Q_{k,i}$ by $\tilde\phi\circ p_{\tilde Q_{k,i,j}}^{-1}$. 
Then $\phi_{k,i,j}$ are the squares of Nash functions, say $\phi^{1/2}_{k,i,j}$, and we can choose $Q_
{k,i}$ so small that the maps $(f,\phi^{1/2}_{k,i,1},...,\phi^{1/2}_{k,i,n-k}):Q_{k,i}\to\R^{m+n-k}$ 
are submersions, that if $Q_{k,i}\cap Q_{k,i'}\not=\emptyset$ then 
$$
\{\phi_{k,i,j}|_{Q_{k,i}\cap Q_{k,i'}}:j=1,...,n-k\}=\{\phi_{k,i',j}|_{Q_{k,i}\cap Q_{k,i'}}:j=1,...,n-
k\},\tag 2
$$
and that if $Q_{k,i}\cap Q_{k',i'}\not=\emptyset$ for $k<k'$ then
$$
\{\phi_{k,i,j}|_{Q_{k,i}\cap Q_{k',i'}}:j=1,...,n-k\}\supset\{\phi_{k',i',j}|_{Q_{k,i}\cap Q_{k',i'}}:
j=1,...,n-k'\}.\tag 3
$$
Let $\Phi_{k,k',i,i'}$ denote the $k'-k$ Nash functions on $Q_{k,i}\cap Q_{k',i'}$ in the complement 
in (3). 
Note that (1) implies 
$$
\phi^{1/2}_{k,i,j}\circ\pi'_X=\phi^{1/2}_{k,i,j}\quad\text{on}\ X\cap Q_{k,i}. \tag"$(1)'$"
$$\par
We work from now in the semialgebraic $C^l$ category. 
Shrink again $Q_{k,i}$ (fixing always $X_k\cap Q_{k,i}$), and set $Q_k=\cup_i Q_{k,i}$. 
Then there exist semialgebraic $C^l$ submersive retractions $q_k:Q_k\to X_k-X_{k-1}$ such that 
$$
\gather
f\circ q_k=f\quad\text{on}\ Q_k,\tag 4\\
q_k\circ\pi'_X=\pi'_X\circ q_k\quad\text{on}\ X\cap Q_k,\tag 5
\endgather
$$
and the maps $(q_k|_{Q_{k,i}},\phi^{1/2}_{k,i,1},...,\phi^{1/2}_{k,i,n-k}):Q_{k,i}\to(X_k-X_{k-1})
\times\R^{n-k}$ are semialgebraic $C^l$ embeddings as follows. \par\noindent
For a while, assume that $q_k$ on $Q_k(0)$ are already given so that the following controlled conditions are satisfied. 
$$
\gather
q_k\circ q_{k'}=q_k\quad\text{on}\ Q_k(0)\cap Q_{k'}(0)\ \text{for}\ k<k',\tag"$(6)_0$"\\
\phi_{k,i,j}^{1/2}\!\circ\! q_{k'}=\phi_{k,i,j}^{1/2}\ \text{on}\ Q_{k,i}(0)\cap Q_{k',i'}(0)\ 
\text{for}\ k<k'\ \text{and}\ \phi_{k,i,j}|_{Q_{k,i}\cap Q_{k',i'}}\in\Phi_{k,k',i,i'}.\tag"$(7)_0$"
\endgather
$$
Extend each $q_k$ to $q_k:X\cap Q_k\to X_k-X_{k-1}$ so that (4) and (5) are satisfied as follows, which 
is uniquely determined, though we need to choose $Q_k$ so that $\pi'_X(X\cap Q_k)\subset Q_k(0)$. 
For $(x,y)\in Q^2_k$ with small $\dis(x,y)$, let $r(x,y)$ denote the orthogonal projection image of $x$ 
to $X_k(f(y))-X_{k-1}$. 
Let $q'_k:Q_k\to X_k-X_{k-1}$ be any semialgebraic $C^l$ extension of $q_k$, shrink $Q_k$ and define $q
_k(x)$ for $x\in Q_k$ to be $r(q'_k(x),x)$. 
Then $q_k$ is a semialgebraic $C^l$ submersive retraction of $Q_k$ to $X_k-X_{k-1}$ and satisfies (4) 
and (5). \par\noindent
Hence $(\phi^{1/2}_{k,i,1},...,\phi^{1/2}_{k,i,n-k})$ is a local 
coordinate system of $q_k^{-1}(x)\cap Q_{k,i}$ at $x$, for each $x\in X_k-X_{k-1}$. 
Therefore, by (3), for each $Q_{k,i}$ and $Q_{k',i'}$ with $k<k'$ there exists a unique semialgebraic 
$C^l$ submersion $q_{k,k',i,i'}:Q_{k,i}\cap Q_{k',i'}\to X_{k'}\cap Q_{k,i}$ such that 
$$
\gather
q_k\circ q_{k,k',i,i'}=q_k\quad\text{on}\ Q_{k,i}\cap Q_{k',i'}\ \text{and}\tag 8\\
\phi_{k,i,j}^{1/2}\circ q_{k,k',i,i'}=\phi_{k,i,j}^{1/2}\quad\text{on}\ Q_{k,i}\cap Q_{k',i'}\ \text{for}
\ \phi_{k,i,j}|_{Q_{k,i}\cap Q_{k',i'}}\in\Phi_{k,k',i,i'}.\tag 9
\endgather
$$
To be precise, the domain of definition of $q_{k,k',i,i'}$ is $q_{k'}^{-1}(Q_{k,i})\cap Q_{k,i}\cap Q
_{k',i'}$. 
However, we omit $q_{k'}^{-1}(Q_{k,i})$ for simplicity of notation. 
In the following arguments also we simplify the domains of definition of many maps. 
By the above equalities we have the following equalities $(4)',\ (5)'$ and (10). 
$$
\gather
f\circ q_{k,k',i,i'}\overset(4)\to=f\circ q_k\circ q_{k,k',i,i'}\overset(8)\to=f\circ q_k\overset(4)
\to=f\quad\text{on}\ Q_{k,i}\cap Q_{k',i'}.\tag"$(4)'$"\\
q_k\circ q_{k,k',i,i'}\circ\pi'_X\overset(8)\to=q_k\circ\pi'_X\overset(5)\to=\pi'_X\circ q_k\overset
(8)\to=\pi'_X\circ q_k\circ q_{k,k'i,i'}\overset(5)\to=q_k\circ \pi'_X\circ q_{k,k',i,i'}\\
\qquad\qquad\qquad\qquad\qquad\qquad\qquad\text{on}\ X\cap Q_{k,i}\cap Q_{k',i'},\\
\phi^{1/2}_{k,i,j}\!\circ\! q_{k,k',i,i'}\!\circ\!\pi'_X\overset(9)\to=\phi^{1/2}_{k,i,j}\!\circ\!
\pi'_X\overset(1)'\to=\phi^{1/2}_{k,i,j}\overset(9)\to=\phi^{1/2}_{k,i,j}\!\circ\! q_{k,k'i,i'}
\overset(1)'\to=\phi^{1/2}_{k,i,j}\!\circ\!\pi'_X\circ q_{k,k',i,i'}\\
\qquad\qquad\qquad\qquad\qquad\text{on}\ X\cap Q_{k,i}\cap Q_{k',i'}\ \text{for}\ \phi_{k,i,j}|_{Q_
{k,i}\cap Q_{k',i'}}\in\Phi_{k,k'i,i'},
\endgather
$$
hence by embeddingness of $(q_k|_{Q_{k,i}},\phi^{1/2}_{k,i,1},...,\phi^{1/2}_{k,i,n-k})$
$$
q_{k,k',i,i'}\circ\pi'_X=\pi'_X\circ q_{k,k',i,i'}\quad\text{on}\ X\cap Q_{k,i}\cap Q_{k',i'}.\tag"$(5
)'$"
$$
By assumption, $(6)_0$ and $(7)_0$ hold. 
Then by (4) and (5) 
$$\gather
q_k\circ q_{k'}=q_k\quad\text{on}\ X\cap Q_k\cap Q_{k'}\ \text{for}\ k<k',\tag"$(6)_X$"\\
\phi_{k,i,j}^{1/2}\circ q_{k'}=\phi_{k,i,j}^{1/2}\ \text{on}\ X\cap Q_{k,i}\cap Q_{k',i'}\ 
\text{for}\ k<k'\ \text{and}\ \phi_{k,i,j}|_{Q_{k,i}\cap Q_{k',i'}}\in\Phi_{k,k',i,i'}.\tag"$(7)_X$"
\endgather
$$
Hence by the same embeddingness as above 
$$
q_{k,k',i,i'}=q_{k'}\quad\text{on}\ X\cap Q_{k,i}\cap Q_{k',i'}.\tag 10
$$\par
Compare $q_{k,k',i_1,i'_1}$ and $q_{k,k',i_2,i'_2}$. 
By (2) and (3)
$$
\gather
q_{k,k',i_1,i'_1}=q_{k,k',i_1,i'_2}\quad\text{on}\ Q_{k,i_1}\cap Q_{k',i'_1}\cap Q_{k',i'_2},\\
q_{k,k',i_1,i'_2}=q_{k,k',i_2,i'_2}\quad\text{on}\ Q_{k,i_1}\cap Q_{k,i_2}\cap Q_{k',i'_2},\\
q_{k,k',i_1,i'_1}=q_{k,k',i_2,i'_2}\quad\text{on}\ Q_{k,i_1}\cap Q_{k,i_2}\cap Q_{k',i'_1}\cap Q_{k',i'
_2}.
\tag"hence"
\endgather
$$
Therefore, we have semialgebraic $C^l$ submersions $q_{k,k'}:Q_k\cap Q_{k'}\to X_{k'}-X_{k'-1},\ 
\allowmathbreak k<k'$, such that 
$$
\gather
f\circ q_{k,k'}=f\quad\text{on}\ Q_k\cap Q_{k'},\tag"$(4)'$"\\
q_k\circ q_{k,k'}=q_k\quad\text{on}\ Q_k\cap Q_{k'},\tag 8\\
\phi_{k,i,j}^{1/2}\circ q_{k,k'}=\phi_{k,i,j}^{1/2}\quad\text{on}\ Q_{k,i}\cap Q_{k',i'}\ \text{for}\ 
\phi_{k,i,j}|_{Q_{k,i}\cap Q_{k',i'}}\in\Phi_{k,k',i,i'},\tag 9\\
q_{k,k'}=q_{k'}\quad\text{on}\ X\cap Q_k\cap Q_{k'}.\tag 10
\endgather
$$\par
We want to shrink the $Q_k$'s and modify the $q_k$'s keeping (4) and (5) so that 
$$
q_{k,k'}=q_{k'}\quad\text{on}\ Q_k\cap Q_{k'}\ \text{for}\ k<k'. \tag 11
$$
We proceed by double induction. 
Let $m\le k_1<k_2<n\in\N$, and assume that (11) holds for $k<k'<k_2$ and for $k_1<k<k'=k_2$. 
Fix $q_k$ and $ k<k_2$. 
Then we need to modify $q_{k_2}$ so that (11) holds for $k=k_1$ and $k'=k_2$. 
Let $\xi$ be a semialgebraic $C^l$ function on $M-X_{k_1-1}$ such that
$0\le\xi\le 1$, and $\xi=1$ outside 
of a small open semialgebraic neighborhood $Q'_{k_1}\,(\subset Q_{k_1})$ of $X_{k_1}-X_{k_1-1}$ in $M-X
_{k_1-1}$ and moreover $\xi=0$ on a smaller one $Q''_{k_1}$. 
Shrink $Q_{k_2}$ and define a semialgebraic $C^l$ submersive retraction $q'_{k_2}:Q_{k_2}\to X_{k_2}-X_
{k_2-1}$ by $q'_{k_2}=q_{k_2}$ on $Q_{k_2}-Q'_{k_1}$ and for $x\in Q_{k_2}\cap Q'_{k_1}$, let $q'_{k_2}(x)$ 
be the orthogonal projection image of $\xi(x)q_{k_2}(x)+(1-\xi(x))q_{k_1,k_2}(x)\in\R^N$ to the Nash 
manifold $X_{k_2}(f(x))-X_{k_2-1}$. 
Then $q'_{k_2}$ satisfies (4) and (11) for $k=k_1$, $k'=k_2$ and for $Q_k$ replaced by $Q''_{k_1}$, the 
map $(q'_{k_2}|_{Q_{k_2,i}},\phi^{1/2}_{k_2,i,1},...,\phi^{1/2}_{k_2,i,n-k_2}):Q_{k_2,i}\to(X_{k_2-X_{k
_2-1}})\times\R^{n-k_2}$ continues to be a semialgebraic $C^l$ embedding if we shrink $Q_{k_2,i}$ (of 
course, fixing $Q_{k_2,i}\cap X_{k_2})$, $q'_{k_2}=q_{k_2}$ on $X\cap Q_{k_2}$ by (10), hence (5) for 
$q'_{k_2}$ holds, and $q'_{k_2}=q_{k_2}$ on $Q_{k_2}\cap\cup_{k_1<k<k_2}Q_k$ for the following reason. 
Let $k_1<k_3<k_2$. 
It suffices to see that $q_{k_1,k_2}=q_{k_2}$ on $Q_{k_1}\cap Q_{k_2}\cap Q_{k_3}$, which is equivalent, 
by uniqueness of $q_{k_1,k_2}$, to
$$
\gather
q_{k_1}\circ q_{k_2}=q_{k_1}\quad\text{on}\ Q_{k_1}\cap Q_{k_2}\cap Q_{k_3},\tag 12\\
\phi_{k_1,i_1,j}^{1/2}\circ q_{k_2}=\phi_{k_1,i_1,j}^{1/2}\quad\text{on}\ Q_{k_1,i_1}\cap Q_{k_2,i_2}
\cap Q_{k_3}\qquad\qquad\tag 13\\
\qquad\qquad\qquad\qquad\text{for}\ \phi_{k_1,i_1,j}|_{Q_{k_1,i_1}\cap Q_{k_2,i_2}}\in\Phi_{k_1,k_2,i_
1,i_2}.
\endgather
$$
By (8) for $k=k_1$ and $k'=k_3$ and for $k=k_3$ and $k'=k_2$ 
$$
\gather
q_{k_1}\circ q_{k_1,k_3}=q_{k_1}\quad\text{on}\ Q_{k_1}\cap Q_{k_3},\\
q_{k_3}\circ q_{k_3,k_2}=q_{k_3}\quad\text{on}\ Q_{k_2}\cap Q_{k_3}. 
\endgather
$$
By (11) for $k=k_1$ and $k'=k_3$ and for $k=k_3$ and $k'=k_2$ 
$$
\gather
q_{k_1,k_3}=q_{k_3}\quad\text{on}\ Q_{k_1}\cap Q_{k_3},\\
q_{k_3,k_2}=q_{k_2}\quad\text{on}\ Q_{k_2}\cap Q_{k_3}. 
\endgather
$$
Hence
$$
\gather
q_{k_1}\circ q_{k_3}=q_{k_1}\quad\text{on}\ Q_{k_1}\cap Q_{k_3},\tag 14\\
q_{k_3}\circ q_{k_2}=q_{k_3}\quad\text{on}\ Q_{k_2}\cap Q_{k_3}.\tag 15
\endgather
$$
Therefore, 
$$
q_{k_1}\circ q_{k_2}\overset(14)\to= q_{k_1}\circ q_{k_3}\circ q_{k_2}\overset(15)\to= q_{k_1}\circ 
q_{k_3}\overset(14)\to= q_{k_1}\quad\text{on}\ Q_{k_1}\cap Q_{k_2}\cap Q_{k_3}.\tag 12
$$
We can prove (13) in the same way because if $Q_{k_1,i_1}\cap Q_{k_2,i_2}\cap Q_{k_3,i_3}\not=
\emptyset$ then $\Phi_{k_1,k_2,i_1,i_2}|_{Q_{k_1,i_1}\cap Q_{k_2,i_2}\cap Q_{k_3,i_3}}$ is the 
disjoint union of $\Phi_{k_1,k_3,i_1,i_3}|_{Q_{k_1,i_1}\cap Q_{k_2,i_2}\cap Q_{k_3,i_3}}$ and $\Phi
_{k_3,k_2,i_3,i_2}|_{Q_{k_1,i_1}\cap Q_{k_2,i_2}\cap Q_{k_3,i_3}}$. 
Thus the induction process works, and we assume that (11) is satisfied. 
Consequently, the following {\bf controlledness} conditions are satisfied by (8), (9) and (11). 
$$\gather
q_k\circ q_{k'}=q_k\quad\text{on}\ Q_k\cap Q_{k'}\ \text{for}\ k<k',\tag 6\\
\phi_{k,i,j}^{1/2}\circ q_{k'}=\phi_{k,i,j}^{1/2}\ \text{on}\ Q_{k,i}\cap Q_{k',i'}\ 
\text{for}\ k<k'\ \text{and}\ \phi_{k,i,j}|_{Q_{k,i}\cap Q_{k',i'}}\in\Phi_{k,k',i,i'}.\tag 7
\endgather
$$\par
It remains to construct $q_k$ on $Q_k(0)$. 
First define $r$ as above, i.e., for $(x,y)\in Q^2_k(0)$ with small $\dis(x,y)$, let $r(x,y)$ denote 
the orthogonal projection image of $x$ to $X_k(0)-X_{k-1}$. 
Set $q_k(x)=r(x,x)$ for $x\in Q_k(0)$. 
Then $q_k:Q_k(0)\to X_k(0)-X_{k-1}$ are Nash submersive retractions. 
We need to modify them so that $(6)_0$ and $(7)_0$ are satisfied.
This is clearly possible by the above arguments. \par
Now we define $W$ and $w$ as in fact 3. 
Set $W=\cup_{k=m}^{n-1}Q_k$ and consider each $Q_{k,i}$. 
Shrink $Q_{k,i}$ so that
$$
(\pi'_X\circ q_k,\phi^{1/2}_{k,i,1},...,\phi^{1/2}_{k,i,n-k})(Q_{k,i})\subset(q_k,\phi^{1/2}_{k,i,1}
,...,\phi^{1/2}_{k,i,n-k})(Q_{k,i}(0)).
$$
Then for each $x\in Q_{k,i}$ there exists a unique $y\in Q_{k,i}(0)$ such that 
$$
(q_k,\phi^{1/2}_{k,i,1},...,\phi^{1/2}_{k,i,n-k})(y)=(\pi'_X\circ q_k,\phi^{1/2}_{k,i,1},...,\phi^
{1/2}_{k,i,n-k})(x).
$$
The correspondence $w'_{k,i}$ from $x$ to $y$ is a semialgebraic $C^l$ map such that $w_{k,i}=(w'_
{k,i},f_{Q_{k,i}}):Q_{k,i}\to Q_{k,i}(0)\times\R^m$ is a semialgebraic $C^l$ embedding by (4), $w'_{k,
i}=\pi'_X$ on $X\cap Q_{k,i}$ because 
$$
(q_k,\phi^{1/2}_{k,i,1},...,\phi^{1/2}_{k,i,n-k})\!\circ\!\pi'_X(x)\!\overset(1)',(5)\to=\!(\pi'_X\circ
q_k,\phi^{1/2}_{k,i,1},...,
\phi^{1/2}_{k,i,n-k})(x)\ \text{for}\ x\in\! X\!\cap Q_{k,i},
$$
and $w'_{k,i}|_{Q_{k,i}(0)}=\id$ by (4) and by the equality $\pi'_X=\id$ on $X(0)$. 
Hence it suffices to see that $w'_{k,i}=w'_{k',i'}$ on $Q_{k,i}\cap Q_{k',i,}$. 
This is clear by (2) if $k=k'$ and $Q_{k,i}\cap Q_{k',i'}\not=\emptyset$. 
Assume that $k<k'$ and $Q_{k,i}\cap Q_{k',i'}\not=\emptyset$. 
By (3) we suppose that
$$
\phi^{1/2}_{k',i',j}=\phi^{1/2}_{k,i,j+k'-k}\quad\text{on}\ Q_{k,i}\cap Q_{k',i'},\ j=1,...,n-k'. 
$$
Then by the definition of $w'_{k,i}$ and $w'_{k',i'}$ we only need to show that 
$$
q_{k'}\circ w'_{k,i}=\pi'_X\circ q_{k'}\quad\text{on}\ Q_{k,i}\cap Q_{k',i'},
$$
which is equivalent to 
$$
\gather
q_k\circ q_{k'}\circ w'_{k,i}=q_k\circ\pi'_X\circ q_{k'}\quad\text{on}\ Q_{k,i}\cap Q_{k',i'}\ \text
{and}\\
\phi^{1/2}_{k,i,j}\circ q_{k'}\circ w'_{k,i}=\phi^{1/2}_{k,i,j}\circ\pi'_X\circ q_{k'}\quad\text{on}\ 
Q_{k,i}\cap Q_{k',i'},\ j=1,...,k'-k.
\endgather
$$
We have 
$$
\gather
q_k\circ q_{k'}\circ w'_{k,i}\overset(6)\to=q_k\circ w'_{k,i}=\pi'_X\circ q_k\overset(6)\to=\pi'_X
\circ q_k\circ q_{k'}\overset(5)\to=q_k\circ\pi'_X\circ q_{k'},\\
\phi^{1/2}_{k,i,j}\circ q_{k'}\circ w'_{k,i}\overset(7)\to=\phi^{1/2}_{k,i,j}\circ w'_{
k,i}=\phi^{1/2}_{k,i,j}\overset(7)\to=\phi^{1/2}_{k,i,j}\circ q_{k'}\overset(1)'\to=\phi^{1
/2}_{k,i,j}\circ\pi'_X\circ q_{k'},\\
\qquad\qquad\qquad\qquad\qquad\qquad\qquad\qquad\qquad\qquad\qquad j=1,...,k'-k.
\endgather
$$
Thus we have completed the proof of fact 3 and hence of the construction of $\pi_i\!=\!(\pi',f_{f^{-1}(B_i)})\!:f^{-1}(B_i)\to M(b_i)\times B_i$. \par
  Next we will extend $\pi_i$ to a neighborhood of $f^{-1}(B_i)$ in $M$. 
Let $\eta_i:U_i\to B_i$ be a semialgebraic submersive $C^l$ retraction of a small semialgebraic open neighborhood of $B_i$ in $R^n$. 
Then we only need to lift $\eta_i$ to a semialgebraic submersive $C^l$ retraction $\tilde\eta_i:f^{-1}(U_i)\to f^{-1}(B_i)$ so that $\tilde\eta_i^{-1}(X_k)=X_k\cap f^{-1}(U_i)$ for each $k$ and $\pi'_X\circ\tilde\eta_i=\pi'_X$ on $X\cap f^{-1}(U_i)$ because if $\tilde\eta_i$ exists, the map $f^{-1}(U_i)\ni x\to(\pi'_i\circ\tilde\eta_i(x),f(x))\in M(b_i)\times U_i$ is the required extension of $\pi_i$. 
We proceed by two steps. 
First we define $\tilde\eta_i$ on $X\cap f^{-1}(U_i)$ and then extend it to $f^{-1}(U_i)$. \par
  The first step. 
We can assume $b_i=0$. 
Then $\pi'_i=\pi'_X$ on $X\cap f^{-1}(B_i)$. 
Hence there exists a unique semialgebraic $C^l$ diffeomorphism $\tilde\eta_{i,y}$ from $X\cap f^{-1}(y)$ to $X\cap f^{-1}(\eta_i(y))$ for each $y\in U_i$ such that $\pi'_i\circ\tilde\eta_{i,y}=\pi'_X$ on $X\cap f^{-1}(y)$. 
Define $\tilde\eta:X\cap f^{-1}(U_i)\to X\cap f^{-1}(B_i)$ by $\tilde\eta_i(x)=\tilde\eta_{i,f(x)}(x)$. 
Then $\tilde\eta_i$ satisfies the requirements. \par
  The second step. 
Since $B_i$ is Nash diffeomorphic to a Euclidean space we can regard $U_i$ as $B_i\times\R^{m'}$ and $\eta_i$ as the projection $\eta_i:B_i\times\R^{m'}\to B_i$, where $m'=m-\dim B_i$. 
Then we define $\tilde\eta_i$ on $f^{-1}(B_i\times\R^k\times\{0\})$ by induction on $k=0,...,m'$. 
For that it suffices to consider the case $m'=1$. 
Moreover we replace $\R$ of $B_i\times\R$ with the circle $S^1=\{x\in\R^2:|x|=1\}$ as follows. 
Let $\omega_i:S^1\to\R$ be a Nash function such that 0 is a regular value. 
Let $\check\eta_i:B_i\times\R\to\R$ be the projection, $\hat M$ be the
fiber product of $\check\eta_i\circ f:f^{-1}(U_i)\to\R$ and
$\omega_i:S^1\to\R$, $\hat X$ be the inverse image of $X\cup f^{-1}(B_i\times\{0\})$ under the induced map $\hat\omega_i:\hat M\to M$ and $\hat f:\hat M\to B_i$ be the naturally defined projection. 
Then $\hat M$ is a Nash manifold, $\hat f$ is a proper Nash map, $\hat X$ is a normal crossing Nash subset of $\hat M$, and the conditions in the lemma are satisfied for $\hat X,\ \hat M$ and $\hat f$. 
Define a map $\hat\pi_{i,\hat X}=(\hat\pi'_{i,\hat X},\hat f):\hat X\to(\hat X\cap\hat f^{-1}(0))\times B_i$ so that $\pi'_X\circ\hat\omega_i\circ\hat\pi'_{i,\hat X}=\pi'_X\circ\hat\omega_i$. 
Then $\hat\pi_{i,\hat X}$ is a uniquely determined Nash
diffeomorphism, $\hat\pi'_{i,\hat X}=\id$ on $\hat X\cap\hat
f^{-1}(0)$, and by fact 3 the map $\hat\pi_{i,\hat X}$ is extended to a semialgebraic $C^l$ embedding $\hat\pi_i=(\hat\pi'_i,\hat f):\hat W_i\to\hat f^{-1}(0)\times B_i$ for some open neighborhood $\hat W_i$ of $\hat X$ in $\hat M$. 
We can shrink $U_i$ and $\hat W_i$ so that $\hat f^{-1}(U_i)=\hat W_i$ since $\hat f$ is proper. 
Hence it remains to consider the problem of lifting $\eta_i$ only on $\eta_i|_{\eta^{-1}(0)}:\eta^{-1}(0)\to\{0\}$. 
Namely the problem is reduced to the case where $B_i=\{0\}$. 
This case also follows from fact 3. 
Thus $\pi_i$ is extended to $f^{-1}(U_i)$. 
We keep the notation $\pi_i$ for the extension. \par
  For the construction of $\pi$ we need to modify and paste $\pi_i$ together. 
This is what [C-S$_2$] proved. 
To be precise, [C-S$_1$] proved local Nash triviality and [C-S$_2$] proved that the local Nash triviality implies the global Nash triviality. 
They treat the case without $X$. 
However, the proof in [C-S$_2$] works in the case with $X$ (see also the proof of Theorem II.6.3, [S$_3$]). 
Thus we obtain $\pi$ and complete the proof of lemma 4.5 in the Nash case. \par
The analytic case follows in the same way. 
\qed
\enddemo
Note that the above proof shows that the lemma still holds if $M$ is a Nash manifold with corners and if the restrictions of $f$ to strata of the canonical 
stratification $\{M_k\}$ of $\partial M$ compatible with $X$ are also submersions onto $\R^m$. 
Here the canonical stratification $\{M_k\}$ compatible with $X$ is defined as follows. 
For a semialgebraic set $S$, let $\Reg S$ denote the subset of $X$ consisting of points $x$ such that 
$S_x$ is a Nash manifold germ of dimension $\dim S$. 
Then $M_{n-1}=\Reg(\partial M-X),\ M_{n-2}=\Reg(\partial M-M_{n-1}),\ M_{n-3}=\Reg(\partial M-M_{n-1}-M
_{n-2}),...$ 
Note that $\{M_k\}$ is a stratification of $\partial M$ into Nash
manifolds of dimension $k$, that $X\cap\partial M$ 
is the union of some connected components of $M_0,...,M_{n-1}$, and the method of construction of $\{M_k
\}$ is {\bf canonical}.
\proclaim{Lemma 4.6}
Let $M$ be a non-compact Nash manifold contained and closed in $\R^N$ and $X$ a normal crossing Nash 
subset of $M$. 
Let $B(r)$ denote the closed ball in $\R^N$ with center 0 and radius $r\in\R$. 
Then there exists a Nash diffeomorphism $\tau:M\to M\cap\Int B(r)$, for some large $r$, such that $\tau(X
)=X\cap\Int B(r)$. 
\endproclaim
This does not necessarily hold in the analytic case. 
\demo{Proof of lemma 4.6}
Assume that $M$ is of dimension $n$. 
Set $X_n=M-X$. 
Choose $r$ so large that the $p|_{X_i-B(r/2)}$ are submersions onto $(r/2,\,\infty)$, where $\{X_i:i=0,
...,n-1\}$ denotes the canonical stratification of $X$ and $p(x)=|x|$ for $x\in M$. 
Then by lemma 4.5 there exists a Nash diffeomorphism $\rho:M-B(r/2)\to(B\cap p^{-1}(r))\times(r/2,\,
\infty)$ of the form $\rho=(\rho',p)$ such that $\rho'(X-B(r/2))=X\cap p^{-1}(r)$. 
Let $\alpha:(-\infty,\,r)\to\R$ be a semialgebraic $C^l$ diffeomorphism such that $\alpha=\id$ on $(-
\infty,\,r/2)$, where $l$ is a sufficiently large integer. 
Set
$$
\tau_0(x)=
\cases
x\quad&\text{for}\ x\in M\cap B(r/2)\\
\rho^{-1}(\rho'(x),\alpha^{-1}\circ p(x))\quad&\text{for}\ x\in M-B(r/2). 
\endcases
$$
Then $\tau_0$ is a semialgebraic $C^l$ diffeomorphism from $M$ to $M\cap\Int B(r)$ such that $\tau_0(
X)=X\cap\Int B(r)$. 
We only need to approximate $\tau_0$ by a Nash diffeomorphism keeping the last property. 
Let $\pi:M\to M\cap\Int B(r)$ be a Nash approximation of $\tau_0$ in the semialgebraic $C^l$ topology. 
Replace $\tau_0$ with $\pi\circ\tau_0^{-1}$. 
Then what we prove is the following statement. \par
Let $\tilde M$ be a compact Nash manifold with boundary in $\R^N$, let
$\tilde X$ be a normal crossing Nash 
subset of $\tilde M$ with $\partial\tilde M\not\subset\tilde X$, and
let $\tilde\tau_0$ be a semialgebraic $C^l$ 
diffeomorphism of $\Int\tilde M$ arbitrarily close to $\id$ in the semialgebraic $C^l$ topology such 
that $\tilde\tau_0(\tilde X\cap\Int\tilde M)$ is a normal crossing Nash subset of $\Int\tilde M$. 
Then we can approximate $\tilde\tau_0$ by a Nash diffeomorphism $\tilde\tau$ of $\Int\tilde M$ in the 
semialgebraic $C^1$ topology so that $\tilde\tau(\tilde X\cap\Int\tilde M)=\tilde\tau_0(\tilde X\cap
\Int\tilde M)$. \par
We proceed as in the proof of step 1, theorem 3.1,(1). 
Let $\{\tilde X_j:j=0,...,n-1\}$ denote the canonical stratification of $\tilde X$ and set $\tilde X_
n=\tilde M-X$. 
By induction, for some $i\in\N$, assume that $\tilde\tau_0|_{\cup_{j=0}^{i-1}\tilde X_j\cap\Int\tilde M}$ is 
of class Nash. 
Let $l'\in\N$. 
Then it suffices to choose $l$ large enough and to approximate $\tilde\tau_0$ by a semialgebraic $C^l$ 
diffeomorphism $\tilde\tau$ of $\Int\tilde M$ in the semialgebraic $C^{l'}$ topology so that $\tilde
\tau(\tilde X\cap\Int\tilde M)=\tilde\tau_0(\tilde X\cap\Int\tilde M)$ and $\tilde\tau|_{\cup_{j=0}^
i\tilde X_j\cap\Int\tilde M}$ is of class Nash. 
Let $\Cal I$ denote the sheaf of $\Cal N$-ideals on $\Int\tilde M$ defined by $\cup_{j=0}^{i-1}\tilde 
X_j\cap\Int\tilde M$. 
By theorem 2.7, the sheaf $\Cal I$ is generated by a finite number of global cross-sections $\xi_1,...,\xi_k$ of 
$\Cal I$. 
Then $\tilde\tau_0|_{\cup_{j=0}^{i-1}\tilde X_j\cap\Int\tilde M}$ is an element of $H^0(\Int\tilde M,
\Cal N/\Cal I)^N$ by the same reason as in the proof in step 1. 
Hence by theorem 2.8 we have a Nash map $h:\Int\tilde M\to\R^N$ such that $h=\tilde\tau_0$ on $\cup_{j
=0}^{i-1}\tilde X_j\cap\Int\tilde M$. 
Here we can choose $h$ to be sufficiently close to $\tilde\tau_0$ in the semialgebraic $C^{l'}$ 
topology for the following reason. 
It suffices to see that $\tilde\tau_0-h$ is of the form $\sum_{j=1}^k\xi_j\beta_j$ for some 
semialgebraic $C^{l'}$ maps $\beta_j:\Int\tilde M\to\R^N$ because $h+\sum_{j=1}^k\xi_j\tilde\beta_j$ fulfills the requirements, where $\tilde\beta_j$ denote Nash approximations of $\beta_j$ in the semialgebraic $C^{l'}$ topology. 
Hence we will prove the following statement. \par
Let $\beta$ be a semialgebraic $C^l$ function on $\Int\tilde M$ vanishing on $\cup_{j=0}^{i-1}
\tilde X_j\cap\Int\tilde M$. 
Then $\beta$ is of the form $\sum_{j=1}^k\beta_j\xi_j$ for some semialgebraic $C^{l'}$ functions 
$\beta_j$ on $\Int\tilde M$. \par
By the second induction, assume that the statement holds for manifolds of dimension strictly less than $n$. 
The problem is reduced to the Euclidean case as follows. 
There exists a finite open semialgebraic covering $\{O_s\}$ of $\Int\tilde M$ such that each $(O_s,
O_s\cap\tilde X)$ is Nash diffeomorphic to $(\R^n,\{(x_1,...,x_n)\in\R^n:x_1\cdots x_{n_s}=0\})$ 
for some $n_s\in\N$. 
Let $\{\eta_s\}$ and $\{\eta'_s\}$ be a partition of unity of class semialgebraic $C^l$ subordinate 
to $\{O_s\}$, and semialgebraic $C^l$ functions on $\Int\tilde M$, respectively, such that $\eta'_s=1$ 
on $\supp\eta_s$ and $\supp\eta'_s\subset O_s$. 
If each $(\beta\eta_s)|_{O_s}$ is described to be of the form $\sum_j\beta_{j,s}\xi_j|_{O_s}$ for some 
semialgebraic $C^{l'}$ functions $\beta_{j,s}$ on $O_s$ then the
naturally defined functions $\sum_s\beta_{j,s}
\eta'_s,$ for $ j=1,...,k$, are semialgebraic $C^{l'}$ functions on $\Int\tilde M$ and $\beta=\sum_j(\sum_s
\beta_{j,s}\eta'_s)\xi_j$. 
Hence we can assume that $(\Int\tilde M,\Int\tilde M\cap X)=(\R^n,\{x_1\cdots x_{n'}=0\})$ for some $n'\in
\N$, and then $n'>0$. 
Apply the induction hypothesis to $\beta|_{\{x_1=0\}}$. 
Then there exist semialgebraic $C^{l_1}$ functions $\beta'_j$ on $\R^{n-1}$ such that 
$$
\beta(0,x_2,...,x_n)=\sum_{j=1}^k\beta'_j(x_2,...,x_n)\xi_j(0,x_2,...,x_n)
$$
because $\Cal I|_{\{x_1=0\}}$ is the sheaf of $\Cal N$-ideals on $\{x_1=0\}$ defined by $\cup_{j=0}^
{i-1}\tilde X_j\cap\{x_1=0\}$ (here $l_1>0$ is arbitrarily given and $l$ depends on $l_1$). 
Regard naturally $\beta'_j$ as semialgebraic $C^{l_1}$ functions on $\R^n$ and replace $\beta$ with 
$\beta-\sum\beta'_j\xi_j$. 
Then we can suppose that $\beta=0$ on $\{x_1=0\}$ from the beginning. 
Under this assumption $\beta/x_1$ is a well-defined semialgebraic $C^{l_1-1}$ function. 
Consider $\beta/x_1$ and $\{x_2\cdots x_{n'}=0\}$ in place of $\beta$ and $\{x_1\cdots x_{n'}=0\}$, 
and repeat the same arguments for $\{x_2=0\}$ and so on. 
Then we finally arrive at the case $\tilde X=\emptyset$. 
Thus the statement is proved, and $h$ is chosen to be close to $\tilde\tau_0$ in the semialgebraic 
$C^{l'}$ topology. \par
Set $Y=\tilde\tau_0(\tilde X\cap\Int\tilde M)$ and $Y_j=\tilde\tau_0(\tilde X_j\cap\Int\tilde M)$. 
Then $Y$ is a normal crossing Nash subset of $\Int\tilde M$, the set $\{Y_j:j=0,...,n-1\}$ is its canonical 
stratification, and $\overline Y$ is a normal crossing semialgebraic $C^l$ subset of $\tilde M$ in 
the sense that $\tilde M$ has a semialgebraic $C^l$ local coordinate system $(x_1,...,x_n)$ at each 
point of $\partial\tilde M$ with $\overline Y=\{x_1\ge 0,\ x_2\cdots x_{n'}=0\}$ for some $n'>0\in\N$ 
by the definition of semialgebraic $C^l$ topology. 
Hence there exists a tubular neighborhood $U_i$ of $Y_i$ in $\R^N$ such that for some $\epsilon>0
\in\R$ 
$$
U_i=\cup_{y\in Y_i}\{x\in\R^N:|x-y|<\epsilon\dis(y,\cup_{j=0}^{i-1}Y_j),\ (x-y)\perp T_y Y_i\}.
$$
Let $q_i:U_i\to Y_i$ denote the orthogonal projection. 
Choose $h$ so close to $\tilde\tau_0$ that $h(\tilde X_i\cap\Int\tilde M)\subset U_i$. 
Then $q_i\circ h|_{\cup_{j=0}^i\tilde X_j\cap\Int\tilde M}$ is a Nash map to $\cup_{j=0}^iY_j$ close 
to $\tilde\tau_0|_{\cup_{j=0}^i\tilde X_j\cap\Int\tilde M}$ in the semialgebraic $C^{l'}$ topology. 
Note that the map is a diffeomorphism by Lemma II.1.7 in [S$_2$]. 
Hence it remains only to extend it to a semialgebraic $C^l$ approximation $\tilde\tau:\Int\tilde M
\to\Int\tilde M$ of $\tilde\tau_0$ in the semialgebraic $C^{l'}$ topology so that $\tilde\tau(\tilde 
X\cap\Int\tilde M)=Y$. 
However, we have already proved it without the last condition. 
Moreover, the proof shows also that the condition is furnished inductively. 
Thus we complete the construction of $\tau$. \qed
\enddemo
\proclaim{Lemma 4.7}
Let $f$ and $g$ be Nash functions on a Nash manifold $M$ which have the same sign at each point of 
$M$, only normal crossing singularities at the common zero set $X$ and the same 
multiplicity at each point of $X$. 
Let $l\in\N$. 
Then there exists a Nash diffeomorphism $\pi$ of $M$ such that $\pi(X)=X$ and $f-g\circ\pi$ is 
$l$-flat at $X$. \par
If $f$ is fixed and $g$ is chosen such that the Nash function on $M$, defined to be $g/f$ on $M-X$, is close 
to 1 in the Nash topology, then $\pi$ is chosen 
to be close to id in the Nash topology. 
\endproclaim
\demo{Proof of lemma 4.7}
Let $M\subset\R^N$, set $n=\dim M$ and let $l$ be sufficiently large. 
For each $k\,(<n)\in\N$, let $X_k$ denote the union of the strata of the canonical stratification of $X$ 
of dimension less than or equal to $k$. 
By induction, assume that $f-g$ is $l$-flat at $X_{k-1}$ for some $k$. 
Then we need only to find a Nash diffeomorphism $\pi$ of $M$ such that $\pi-\id$ is $l$-flat at $X_{
k-1}$, such that $\pi(X)=X$ and $f-g\circ\pi$ is $l$-flat at $X_k$ (to be precise, we will construct $\pi$ so that 
$\pi-\id$ and $f-g\circ\pi$ are $l^{(4)}$-flat at $X_{k-1}$ and $X_k$, respectively, for some $0\ll l^{(4
)}\ll\cdot\cdot\ll l'\ll l$). \par
We proceed as in the proof of lemma 4.5. 
Let $(\tilde M,\tilde X)$ and $(\tilde M_k,\tilde X_k)$ be pairs of Nash manifolds and Nash 
submanifolds, let $p:\tilde M\to M$ and $p_k:\tilde M_k\to M$ be Nash
immersions and let $q_k:\tilde M_k\to
\tilde X_k$ be a Nash submersive retraction such that $\dim\tilde
M=\dim\tilde M_k=n$, the equalities $ p(\tilde X)=X$ and $ p_k(\tilde
X_k)=X_k$ hold, and moreover $ p|_{\tilde X-p^{-1}(X_{n-2})}$ and $p_k|_{\tilde X_k-p^{-1}_k(X_{k-1})}$ 
are injective, and $p_k(q^{-1}_k(\tilde X_k\cap p_k^{-1}(X_{k-1})))\subset X$. 
Shrink $\tilde M_k$ if necessary. 
Then we have an open semialgebraic neighborhood $U$ of $\tilde X\cap p^{-1}(X_k)$ in $\tilde M$ and 
a Nash $(n-k)$-fold covering map $r:U\to\tilde M_k$ such that $p_k\circ r=p$ on $U$. 
Let $\tilde\phi$ be a Nash function on $\tilde M$ with zero set $\tilde X$ which is, locally at each 
point of $\tilde X$, the square of a regular function. 
Then $\tilde\phi(r^{-1}(x))$ is a family of $(n-k)$-numbers possibly 
with multiplicity, for each $x\in\tilde M_k$, and there exist Nash functions $\tilde\phi_{k,1},...,\tilde\phi_{k,n-k}$ on an 
open semialgebraic neighborhood of each point of $\tilde M_k$ such that $\tilde\phi(r^{-1}(x
))=\{\tilde\phi_{k,1}(x),...,\tilde\phi_{k,n-k}(x)\}$ for $x$ in the
given neighborhood. 
For simplicity of notation, we assume that $\tilde\phi_{k,1}(x),...,\tilde\phi_{k,n-k}(x)$ are defined 
globally, which causes no problem because the following arguments are done locally and do not 
depend on the order of $\tilde\phi_{k,1}(x),...,\tilde\phi_{k,n-k}(x)$. 
Moreover, we suppose that each $\tilde\phi_{k,i}$ is the square of a regular Nash function, say $\tilde
\phi^{1/2}_{k,i}$, by the same reason as above. 
Set $\tilde f_k=f\circ p_k$ and $\tilde g_k=g\circ p_k$. \par
We want to construct a Nash diffeomorphism $\tilde\pi_k$ between semialgebraic neighborhoods of 
$\tilde X_k$ in $\tilde M_k$ such that
$\tilde\pi_k(p^{-1}_k(X))\subset p^{-1}_k(X)$, such that $ \tilde\pi_k-
\id$ is $l''$-flat at $q^{-1}_k(\tilde X_k\cap p_k^{-1}(X_{k-1}))$ and $\tilde f_k-\tilde g_k\circ
\tilde\pi_k$ is $l''$-flat at $\tilde X_k$. 
Assume that $\tilde X_k$ is connected without loss of generality. 
Since $\tilde f_k$ and $\tilde g_k$ have only normal crossing singularities at $p_k^{-1}(X)$, the 
same sign at each point of $\tilde M_k$ and the same multiplicity at each point of $p^{-1}_k(X)$, 
and since $\tilde f_k^{-1}(0)=\tilde g^{-1}_k(0)=\cup_{i=1}^{n-k}(\tilde\phi^{1/2}_{k,i})^{-1}(0)
\cup q^{-1}_k(\tilde X_k\cap p_k^{-1}(X_{k-1}))$, we have Nash functions $F$ and $G$ on $\tilde M
_k$ and $\alpha=(\alpha_1,...,\alpha_{n-k})\in(\N-\{0\})^{n-k}$ such
that the equalities $\tilde f_k=F\tilde\phi
^{1/2\alpha}_k$ and $ \tilde g_k=G\tilde\phi^{1/2\alpha}_k$ hold, such
that $FG\ge0$ on $\tilde M_k$ and $FG>0$ on 
$\tilde M_k-q_k^{-1}(\tilde X_k\cap p_k^{-1}(X_{k-1}))$, where $\tilde\phi^{1/2\alpha}_k=\prod_{i
=1}^{n-k}\tilde\phi^{1/2\alpha_i}_{k,i}$. 
Assume that $F\ge0$ and hence $G\ge0$ (the other cases can be proved in the same way). 
Note that $F$ and $G$ have zero set $q_k^{-1}(\tilde X_k\cap
p_k^{-1}(X_{k-1}))$, which has only normal crossing 
singularities and has the same multiplicity at each point. 
Shrink $\tilde M_k$ so that the map $(q_k,\tilde\phi^{1/2}_{k,1},...,\tilde\phi^{1/2}_{k,n-k}):
\tilde M_k\to\tilde X_k\times\R^{n-k}$ is a Nash embedding and let $V$ denote its image. Identify 
$\tilde M_k$ and $\tilde X_k$ with $V$ and $\tilde X_k\times\{0\}$ through this embedding, set $\tilde p
_k=p_k|_{\tilde X_k}$, regard $p_k$ as an immersion of $V$ into $M$ and $\tilde f_k$ and $\tilde g
_k$ as functions on $V$, and let $(z,y)
=(z,y_1,...,y_{n-k})\in V\subset \tilde X_k\times\R^{n-k}$. 
Then
$$
\tilde f_k(z,y)=F(z,y)y^\alpha\quad\text{and}\quad\tilde g_k(z,y)=G(z,y)y^\alpha. 
$$
Set
$$
\gather
F'=\sum_{\beta\in\N_l^{n-k}}\frac{\partial^{|\beta|} F}{\partial y^\beta}(z,0)y^\beta/\beta!,\quad
G'=\sum_{\beta\in\N_l^{n-k}}\frac{\partial^{|\beta|} G}{\partial y^\beta}(z,0)y^\beta/\beta!,\\
\tilde f'_k=F'y^\alpha\quad\text{and}\quad\tilde g'_k=G'y^\alpha,
\endgather
$$
where $\N_l^{n-k}=\{\beta\in\N^{n-k}:|\beta|\le l\}$ and $\beta!=\prod_{i=1}^{n-k}\beta_i!$. 
Then $\tilde f'_k$ and $\tilde g'_k$ are Nash functions on $V$,
moreover $\tilde f_k-\tilde f'_k$ and $
\tilde g_k-\tilde g'_k$ are $l$-flat at $\tilde X_k\times\{0\}$, and $F'$ and $G'$ have the same 
properties as $F$ and $G$. 
Hence for the construction of $\tilde\pi_k$, we can replace $\tilde f_k$ and $\tilde g_k$ with $\tilde f'
_k$ and $\tilde g'_k$. 
An advantage of $\tilde f'_k$ and $\tilde g'_k$ is the fact that $(*)$ $F'-G'$ is $l'$-flat at $V\cap
\tilde p^{-1}_k(X_{k-1})\times\R^{n-k}$, though $F-G$ is $l'$-flat only at $\tilde p^{-1}_k(X_{k-1}
)\times\{0\}$. 
Write 
$$
\tilde f'_k=\prod_{i=1}^{n-k}(F^{\prime1/|\alpha|}y_i)^{\alpha_i}\quad\text{and}\quad\tilde g'_k=
\prod_{i=1}^{n-k}(G^{\prime1/|\alpha|}y_i)^{\alpha_i}. 
$$
Then there exists a unique Nash diffeomorphism $\tilde\pi_k$ between semialgebraic neighborhoods 
of $\tilde X_k\times\{0\}$ in $V$ of the form $\tilde\pi_k(z,y)=(z,\tilde\pi'_k(z,y)y)$, for some 
positive Nash function $\tilde\pi'_k$ on the neighborhood of source, such that $\tilde f'_k=
\tilde g'_k\circ\tilde\pi_k$ on that neighborhood. Actually, we can reduce the problem to the case where 
$\tilde g'_k=z^\beta y^\alpha$ for some $\beta\in\N^k$ and some local Nash coordinate system $z=(z_
1,...,z_k)$ of $\tilde X_k$ such that $\tilde p_k^{-1}(X_{k-1})=\{z^\beta=0\}$ (by considering two 
pairs $(\tilde f'_k,z^\beta y^\alpha)$ and $(\tilde g'_k,z^\beta y^\alpha)$) and then $\tilde\pi'_
k(z,y)=(F'/z^\beta)^{1/|\alpha|}$ is the unique solution. 
Such a $\tilde\pi_k$ fulfills the requirements. 
Indeed, $\tilde\pi_k(p^{-1}_k(X))\subset p^{-1}_k(X)$ by the form of $\tilde\pi_k$ because $p^{-1}_
k(X)$ in $V$ is of the form $\tilde X_k\times\{y_1\cdots y_{n-k}=0\}\cup\tilde p^{-1}_k(X_{k-1})
\times\R^{n-k}$, because $\tilde\pi_k-\id$ is $l''$-flat at $V\cap\tilde p^{-1}_k(X_{k-1})\times\R^{n-k}$ 
because of $(*)$, and $\tilde f_k-\tilde g_k\circ\tilde\pi_k$ is $l''$-flat at $\tilde X_k\times\{0\}$ 
because
$$
\tilde f_k-\tilde g_k\circ\tilde\pi_k=(\tilde f_k-\tilde f'_k)+(\tilde f'_k-\tilde g'_k\circ\tilde
\pi_k)+(\tilde g'_k\circ\tilde\pi_k-\tilde g_k\circ\tilde\pi_k).
$$\par
Let $W$ be an open semialgebraic neighborhood of $X_k-X_{k-1}$ in $M$ so small that there exists 
an open semialgebraic neighborhood of $(\tilde X_k-\tilde p^{-1}_k(X_{k-1}))\times\{0\}$ in the 
intersection of the domain of definition of $\tilde\pi_k$ and the range of values to which the 
restriction of $p_k$ is a diffeomorphism onto $W$. 
Then $\tilde\pi_k$ induces a Nash embedding $\pi_k:W\to M\subset\R^n$ such that $\pi_k(X\cap W)
\subset X$, such that $\pi_k=\id$ on $X_k-X_{k-1}$ and $f-g\circ\pi_k$ is $l''$-flat at $X_k-X_{k-1}$. 
Though $\pi_k$ is not necessarily extensible to a neighborhood of $X_k$, there exists a Nash map 
$\eta:M\to\R^N$ such that $\eta-\id$ is $l''$-flat at $X_{k-1}$, and $\eta-\pi_k$ is $l''$-flat at $X_k
-X_{k-1}$, hence $f-g\circ\eta$ is $l''$-flat at $X_k$ for the following reason. 
Let $\Cal I$ denote the sheaf of $\Cal N$-ideals on $M$ defined by $X_k$. 
Then by theorem 2.8, it suffices to find an element $\overline\eta$ in $H^0(M,\Cal N/\Cal I^{l''})^N$ 
such that $\overline\eta_x$ is the image of $\pi_{kx}$ under the natural map $\Cal N^n_x\to(\Cal N_x
/\Cal I_x^{l''})^N$ for $x\in W$ and $\overline\eta_x=\id$ for $x\in X_{k-1}$. 
This is possible because $\tilde\pi_k-\id$ is $l''$-flat at $V\cap\tilde p_k^{
-1}(X_k)\times\R^{n-k}$. \par
We modify $\tilde\pi_k$ to show that $\eta$ can be a diffeomorphism of $M$. 
Assume that $(**)\ \tilde\pi'_k\le1$ for simplicity of notation, which is possible if we consider 
a third function $h$ on $M$ with the same properties as $f$ and $g$,
with $h/f\ge1$ on $M-X$ and $h/g\ge1$ 
on $M-X$. 
Let $\psi$ be a non-negative small Nash function on $\tilde X_k$ with zero set $\tilde p^{-1}_k(X_{
k-1})$ such that 
$$
\gather
Z\overset{\text{def}}\to=\{(z,y)\in\tilde X_k\times\R^{n-k}:|y|\le\psi(z)\}\subset \text{domain of}
\ \tilde\pi_k,\\
\tilde\pi'_k(z,sy)>|\frac{\partial\tilde\pi'_k(z,sy)}{\partial s}s|\qquad\qquad\qquad\qquad\qquad
\qquad\tag"$(3*)$"\\
\qquad\qquad\qquad\qquad\text{for}\ (z,y,s)\in\tilde X_k\times\R^{n-k}\times\R\ \text{with}\ (z,sy)
\in Z\ \text{and}\ |y|=1
\endgather
$$
and $p_k|_Z$ is injective, which exists by the \L ojasiewicz inequality. 
Let $\rho(t)$ be a semialgebraic $C^{l''}$ function on $\R$ such that
$(4*)\ 0\le\rho\le1$, such that $(5*)\ 
\frac{d\rho}{dt}\le0$, and moreover $ \rho=1$ on $(-\infty,\,1/2]$ and $\rho=0$ on $[1,\,\infty)$. 
Set 
$$
\gather
\tilde\tau'_k=
\cases
1&\text{for}\ (z,y)\in Z\cap p_k^{-1}(X_{k-1})\\
\rho(|y|/\psi(z))\tilde\pi'_k(z,y)+1-\rho(|y|/\psi(z))\ &\text{for}\ (z,y)\in Z-p_k^{-1}(X_{k-1}),
\endcases\\
\tilde\tau_k(z,y)=(z,\tilde\tau'_k(z,y)y)\quad\text{for}\ (z,y)\in Z.
\endgather
$$
Then $\tilde\tau'_k$ and hence $\tilde\tau_k$ are of class semialgebraic $C^{l^{(3)}}$ and $\tilde
\tau_k-\id$ is $l^{(3)}$-flat at $Z\cap p_k^{-1}(X_{k-1})=\tilde p_k^{-1}(X_{k-1})\times\{0\}$ 
since $\tilde\pi'_k(z,y)-1$ is $(l''-1)$-flat at $Z\cap p_k^{-1}(X_{k-1})$. 
Clearly $\tilde\tau_k=\id$ on a semialgebraic neighborhood of $\partial Z-p_k^{-1}(X_{k-1})$ in $Z$. 
Moreover, $\tilde\tau_k$ is a diffeomorphism of $Z$. Actually, we can assume that $n-k=1$ because $\tilde\tau_k=\tilde\pi_k$ on a neighborhood of $(\tilde X_
k-\tilde p_k^{-1}(X_{k-1}))\times\{0\}$ in $Z$ and because $\tilde\pi_k$ and hence $\tilde\tau_k$ 
carry each segment $\{z\}\times\{\R y\}\cap Z$ for $(z,y)\in(\tilde X_k-\tilde p_k^{-1}(X_{k-1}))
\times\R^{n-k}$ with $|y|=1$ to itself. 
Then 
$$
\gather
\frac{\partial\tilde\tau'_k(z,y)y}{\partial y}=\tilde\tau'_k(z,y)+\frac{\partial\tilde\tau'_k}{
\partial y}(z,y)y,\\
\tilde\tau'_k(z,y)=\rho(|y|/\psi(z))\tilde\pi'_k(z,y)+1-\rho(|y|/\psi(z))\overset(**),(4*)\to\ge
\tilde\pi'_k(z,y),\\
\frac{\partial\tilde\tau'_k}{\partial y}(z,y)y=\frac{d\rho}{dt}(|y|/\psi(z))(\tilde\pi'_k(z,y)-1)|
y|/\psi(z)+\rho(|y|/\psi(z))\frac{\partial\tilde\pi'_k}{\partial y}(z,y)y\overset(**),(5*)\to\ge\\
\qquad\qquad\qquad\qquad\qquad\rho(|y|/\psi(z))\frac{\partial\tilde\pi'_k}{\partial y}(z,y)y,\quad
\text{hence}\\
\frac{\partial\tilde\tau'_k(z,y)y}{\partial y}\ge\tilde\pi'_k(z,y)+\rho(|y|/\psi(z))\frac{\partial
\tilde\pi'_k}{\partial y}(z,y)y\overset(3*),(4*)\to>0\quad\text{for}\ (z,y)\in Z. 
\endgather
$$\par
Define a semialgebraic $C^{l^{(3)}}$ diffeomorphism $\tau_k$ of $M$ so that $\tau_k\circ p_k=p_k
\circ\tilde\tau_k$ on $Z$ and $\tau_k=\id$ on $M-p_k(Z)$. 
Then $\tau_k=\pi_k$ on $Z$ if we shrink $Z$, hence $\tau_k-\eta$ is
$l^{(3)}$-flat at $X_k$ and $\tau_
k(X)=X$ and moreover $f-g\circ\tau_k$ is $l^{(3)}$-flat at $X_k$. 
Let $\omega$ be a non-negative-valued global generator of the square of $\Cal I$---the sheaf of 
$\Cal N$-ideals defined by $X_k$. 
Then there exists a semialgebraic $C^{l^{(4)}}$ map $\xi:M\to\R^N$ such that $\tau_k-\eta=\omega\xi$. 
Approximate $\xi$ by a Nash map $\xi'$ in the semialgebraic $C^{l^{(4)}}$ topology, and set $\pi=
(\eta+\omega\xi')\circ o$, where $o$ denotes the orthogonal projection to $M$ of its semialgebraic 
tubular neighborhood in $\R^N$. 
Then $\pi$ is a Nash diffeomorphism of $M$ such that $\pi-\id$ is $l^{(4)}$-flat at $X_{k-1}$ and 
$f-g\circ\pi$ is $l^{(4)}$-flat at $X_k$. 
We can modify $\pi$ so that $\pi(X)=X$ in the same way as in step 1 of the proof of theorem 3.1,(1) and 
lemma 4.6, because $\pi$ is an approximation of $\tau_k$ and $\tau_k(X)=X$. 
Thus we complete the proof of the former half of lemma 4.7. \par
The latter half automatically follows from the above proof (though $(**)$ does not necessarily hold, $\pi'_k$ is close 
to 1 in the Nash topology, which is sufficient to proceed). 
\qed
\enddemo
Note that our proof of lemma 4.7 also works when $M,\, f$ and $g$ are
of class $C^\omega$ and the multiplicities of $f$ and $g$ are bounded. \par
The following lemma is also a globalization of Chapter II, Proposition 2 in [T] and shows sufficient 
conditions for two functions to be right equivalent. 

\proclaim{Proposition 4.8}
(i) Let $f$ be a $C^\omega$ function on a $C^\omega$ manifold $M$. 
Let $v_i$, for $ i=1,...,k$, be $C^\omega$ vector fields on $M$, and $I$ denote the ideal of $C^\infty(M)$ 
or $C^\omega(M)$ generated by $v_if$, for $ i=1,...,k$. 
Let $\phi$ be a small $C^\infty$ or $C^\omega$ function on $M$ contained in $I^2$ in the strong 
Whitney $C^\infty$ topology. 
Then $f$ and $f+\phi$ are $C^\infty$ or $C^\omega$ right equivalent, respectively, and the 
diffeomorphism of equivalence can be chosen to be close to id in the same topology. \par
(ii) If $f,\ M$ and $v_i$ are of class Nash or $C^\infty$ or $C^\omega$, and $\phi$ is of the form 
$\sum_{i,j=1}^k\phi_{i,j}v_if\cdot v_jf$ for some small Nash or $C^\infty$ or $C^\omega$ functions 
$\phi_{i,j}$ in the Nash or (strong) Whitney $C^\infty$ topology, then $f$ and $f+\phi$ are Nash or 
$C^\infty$ or $C^\omega$ right equivalent, respectively, by a Nash or $C^\infty$ or $C^\omega$ 
diffeomorphism close to id in the same topology. \par
(iii) Assume that $M$ is a Nash manifold and $f$ is a Nash function on $M$ with only normal crossing 
singularities. 
Set $X=f^{-1}(f(\Sing f))$. 
Let $\phi$ be a Nash function on $M$ $r$-flat at $X$ for some large $r\in\N$. 
Then there exists a Nash diffeomorphism $\pi:V_1\to V_2$ between closed semialgebraic neighborhoods 
of $X$ in $M$ close to id in the semialgebraic $C^{r'}$ topology, for $0<r'\,(\ll r)\in\N$, such that 
$f\circ\pi=f+\phi$ on $V_1$, such that $\pi-\id$ is $r'$-flat at $X$, and $\pi$ is extensible to a semialgebraic 
$C^r$ diffeomorphism of $M$. 
\endproclaim
\demo{Proof of proposition 4.8}
Consider the analytic case. 
We want to reduce (i) to (ii). 
For a while we proceed in the strong Whitney $C^\infty$ topology. 
By lemma 1.12 for $\phi$ in (i), there exist small $\phi_{i,j}\in
C^\omega$(M), for $i,j=1,...,k$, such that 
$\phi=\sum_{i,j=1}^k\phi_{i,j}v_if\cdot v_jf$. 
Consequently, (i) is reduced to (ii). 
From now on, we work in the Whitney $C^r$ topology for any $r>0\in\N$ (though we can do in the 
strong Whitney $C^\infty$ topology). 
We can assume that $M$ is open in its ambient Euclidean. Actually, let $p:\tilde M\to M$ denote the orthogonal projection of a tubular neighborhood of $M$ in its 
ambient Euclidean space. 
Assume that proposition 4.8,(ii) in the analytic case holds for $\tilde M$. 
The map $C^\omega(\tilde M)\ni\Psi\to\Psi|_M\in C^\omega(M)$ is obviously continuous, surjective 
by corollary 2.4 and open as follows. 
Let $\xi\in C^\infty(\tilde M)$ with $\xi=1$ on $M$ and $\xi=0$ outside of a small neighborhood 
of $M$ in $\tilde M$. 
Then the map $\xi C^\omega(\tilde M)\ni\xi\Psi\to\xi\Psi|_M\in C^\omega(M)$ is open because for 
$\psi\in C^\omega(M)$ and for $\Psi_0\in C^\infty(\tilde M)$, we have $\psi\circ p\in C^\omega(
\tilde M)$ and $(\xi\cdot\psi\circ p)|_M=\psi$ and the map $C^\omega(M)\ni\psi\to\xi\cdot\psi
\circ p+\xi\Psi_0-\xi\cdot\Psi_0|_M\circ p\in\xi C^\omega(\tilde M)$ is continuous and carries 
$\Psi_0|_M$ to $\xi\Psi_0$. 
Hence for small $\psi\in C^\omega(M)$, there exists small $\xi\Psi\in\xi C^\omega(\tilde M)$ such 
that $\Psi|_M=\psi$. 
Approximate $\xi$ by an analytic function $\xi'$ on $\tilde M$ so that $\xi'=1$ on $M$. 
Then $\xi'\Psi$ is analytic on $\tilde M$, close to $\xi\Psi$ and hence small since the map $C^
\infty(\tilde M)^2\ni(\alpha,\beta)\to\alpha\beta\in C^\infty(\tilde M)$ is continuous, and $\xi'
\Psi|_M=\psi$. 
Consequently, the above restriction map $\Psi\to\Psi|_M$ is open by linearity. 
Let $\tilde v_i$, for $ i=1,...,k$, be $C^\omega$ vector field extensions of $v_i$ to $\tilde M$, and 
$\tilde\phi_{i,j}$ $C^\omega$ extensions of $\phi_{i,j}$ to $\tilde M$ so small that $f
\circ p$ and $f\circ p+\sum_{i,j}^k\tilde\phi_{i,j}\tilde v_i(f\circ p)\cdot\tilde v_j(f\circ p)$ 
satisfy the condition in proposition 4.8,(ii) and hence are $C^\omega$ right equivalent by a $C^
\omega$ diffeomorphism $\tilde\pi$ close to id, i.e.,
$$
f\circ p\circ\tilde\pi=f\circ p+\sum_{i,j=1}^k\tilde\phi_{i,j}\tilde v_i(f\circ p)\cdot\tilde v_j
(f\circ p)\quad\text{on}\ \tilde M. 
$$
Set $\pi=p\circ\tilde\pi|_M$. 
Then $\pi$ is a $C^\omega$ diffeomorphism of $M$ close to id, and 
$$
f\circ\pi=f+\sum_{i,j=1}^k\phi_{i,j}v_if\cdot v_jf. 
$$
Thus proposition 4.8,(ii) is proved for $M$. 
Hence we assume that $M$ is open in $\R^n$. \par
Next we can suppose that $k=n$ and $v_i=\frac\partial{\partial x_j}$,
for $ i=1,...,n,$ because each $v_i$ 
is written as $\sum_{j=1}^n\alpha_{i,j}\frac\partial{\partial x_j}$ for some $C^\omega$ functions 
$\alpha_{i,j}$ on $M$. \par
Let $\eta$ denote the function on $M$ which measures distance from $\partial M\overset\text{def}
\to=\overline M-M$ (if $\partial M=\emptyset$ then set $\eta\equiv+\infty$). 
Set $V=\{(x,y)\in M\times\R^n:|y|<\eta(x)\}$ and consider the $C^\omega$ function 
$$
g(x,y)=f(x+y)-f(x)-\sum_{i=1}^ny_i\frac{\partial f}{\partial x_i}(x)\quad\text{for}\ (x,y)=(x_1,...
,x_n,y_1,...,y_n)\in V. 
$$
Then $g$ is a global cross-section of the sheaf of $\Cal O$-ideals $\Cal I$ on $V$ generated by $y_
iy_j$, for $ i,j=1,...,n$. 
Hence applying theorem 2.3 to the surjective homomorphism $\Cal O^{n^2}\ni(\alpha_{i,j})\to\sum_{
i,j=1}^n\alpha_{i,j}y_iy_j\in\Cal I$ we obtain $C^\omega$ functions $g_{i,j}$ on $V,\ i,j=1,...,n$, 
such that $g(x,y)=\sum_{i,j=1}^ny_iy_jg_{i,j}(x,y)$. 
Then 
$$
f(x+y)=f(x)+\sum_{i=1}^ny_i\frac{\partial f}{\partial x_i}(x)+\sum_{i,j=1}^ny_iy_jg_{i,j}(x,y). 
\tag$*$
$$
Let $\alpha=(\alpha_{i,j})_{i,j=1,...,n}$ be new variables in $\R^{n^2}$, set 
$$
\langle\alpha,\partial f\rangle=(\sum_{i=1}^n\alpha_{i,1}\frac{\partial f}{\partial x_i}(x),...,
\sum_{i=1}^n\alpha_{i,n}\frac{\partial f}{\partial x_i}(x)), 
$$
and let $W$ be a small open neighborhood of $M\times\{0\}$ in $M\times\R^{n^2}$ such that 
$$
(x,\langle\alpha,\partial f\rangle)\in V\quad\text{for}\ (x,\alpha)\in W.
$$
Take $y$ to be $\langle\alpha,\partial f\rangle$ in $(*)$. 
Then 
$$
\multline
\quad f(x+\langle\alpha,\partial f\rangle)=\\
f(x)+\sum_{i,j}^n\alpha_{i,j}\frac{\partial f}{\partial x_i}(x)\frac{\partial f}{\partial x_j}
(x)+\sum_{i,i',j,j'=1}^n\alpha_{i,i'}\alpha_{j,j'}\frac{\partial f}{\partial x_i}(x)\frac{\partial 
f}{\partial x_j}(x)G_{i',j'}(x,\alpha)
\endmultline
$$
for $C^\omega$ functions $G_{i',j'}(x,\alpha)=g_{i',j'}(x,\langle\alpha,\partial f\rangle)$ on $W$. 
Consider the map 
$$
B:W\ni(x,\alpha)\to(x,\alpha_{i,j}+\sum_{i',j'=1}^n\alpha_{i,i'}\alpha_{j,j'}G_{i',j'}(x,\alpha))
\in M\times\R^{n^2}. 
$$
Then $B$ is id and regular at $M\times\{0\}$. 
Hence, shrinking $W$, we assume that $B$ is a diffeomorphism onto an open neighborhood $O$ of $M\times\{0
\}$ in $M\times\R^{n^2}$. 
Set $B(x,\alpha)=(x,B_{i,j}(x,\alpha))$, and $B^{-1}(x,\beta)=(x,A'(x,\beta))$ for $(x,\beta)\in O$. 
Then $A'$ is a $C^\omega$ map from $O$ to $\R^{n^2}$, 
$$
\gather
f(x+\langle\alpha,\partial f\rangle)=f(x)+\sum_{i,j}^nB_{i,j}(x,\alpha)\frac{\partial f}{\partial 
x_i}(x)\frac{\partial f}{\partial x_j}(x)\quad\text{for}\ (x,\alpha)\in W, \\
f(x+\langle A'(x,\beta),\partial f\rangle)=f(x)+\sum_{i,j}^n\beta_{i,j}\frac{\partial f}{\partial x_
i}(x)\frac{\partial f}{\partial x_j}(x)\quad\text{for}\ (x,\beta)\in O. 
\endgather
$$
Choose $\Phi=(\phi_{i,j})$ so small that its graph is contained in $O$. 
Then $\pi(x)=x+\langle A'(x,\Phi(x)),\partial f\rangle$ fulfills the requirements in (ii). 
Here if $\phi_{i,j}$ are small in the Whitney $C^r$ or the strong Whitney $C^\infty$ topology, $\pi$ 
is close to id in the respective topology. \par
If $f,\, M$ and $v_i$ are of class $C^\omega$ and if $\phi$ is of class $C^\infty$, the same 
arguments as above work and the diffeomorphism of equivalence is of class $C^\infty$. 
Thus we complete the proof of (ii) in the analytic case. 
Point (ii) in the $C^\infty$ or Nash case follows also from the same proof. 
The difference is only that the existence of $C^\infty$ or Nash $g_{i,j}$ follows from a partition of 
unity of class $C^\infty$ or theorem 2.8, respectively. \par
Consider (iii). 
Assume that $M$ is not compact. 
Let $M$ be embedded in a Euclidean space so that its closure is a compact Nash manifold with boundary. 
Now, we consider an open semialgebraic tubular neighborhood of $\overline M$ and extend $f$ to 
the neighborhood as before. 
Then we can assume that $M$ is open in $\R^n$ and $\overline M$ is a compact Nash manifold with corners, 
and for the construction of $\pi$ it suffices to see that $\phi$ is of the form $\sum_{i,j}^n\phi_{i,j}
\frac{\partial f}{\partial x_i}\frac{\partial f}{\partial x_j}$ for some Nash functions $\phi_{i,j}$ on 
$M$ $r'$-flat at $X$, where $0\ll r'\ll r\in\N$. Actually, assume that there exist such $\phi_{i,j}$. 
Then by the above proof, we only need to find small semialgebraic $C^{r''}$ functions $\phi'_{i,j}$ on 
$M$ in the semialgebraic $C^{r''}$ topology such that $\phi'_{i,j}=\phi_{i,j}$ on some semialgebraic 
neighborhood of $X$ for $0<r''\ll r'\in\N$. \par
Consider only the case $r''=1$ because the general case can be proved in the same way. 
Set $g(x)=\prod_{a\in f(X)}(f(x)-a)^2$, and let $h$ be a Nash function on $M$ extensible to a Nash 
function $\overline h$ on $\overline M$ such that $0<h\le1/2$, such
that $ (1)\ |\frac{\partial h}{\partial x_k}
|\le 1,$ for $ k=1,...,n$, and $\overline h^{-1}(0)=\overline M-M$, which exists since $\overline M$ is a 
compact Nash manifold with corners. 
Let $\psi(t)$ be a semialgebraic $C^1$ function on $\R$ such that
$0\le\psi\le1$ and $ \psi=1$ on $(-
\infty,\,1]$ whereas $\psi=0$ on $[2,\,\infty)$. 
Then $\phi'_{i,j}=\phi_{i,j}\psi(g/h^m)$ fulfill
the requirements for some $m\in\N$. Actually,
Clearly $\phi'_{i,j}=\phi_{i,j}$ on a semialgebraic neighborhood $\{x\in M:g(x)\le h^m(x)\}$ of $X$ 
in $M$, and $\phi'_{i,j}=0$ on $\{g(x)\ge2h^m(x)\}$. 
Hence we prove that each $\phi'_{i,j}$ is small on $V\overset\text{def}\to=\{g(x)\le2h^m(x)\}$ in the 
semialgebraic $C^1$ topology. 
Let $\epsilon>0\in\R$. 
Let $\xi$ denote the Nash function on $M$ defined to be $\phi_{i,j}/g^2$ on $M-X$ and 0 on $X$. 
Then $\xi,\frac{\partial g}{\partial x_k}$ and $\frac{\partial\phi_{i,j}}{\partial x_k},\ k=1,...,n$, 
vanish at $X$. 
Hence there exists a semialgebraic neighborhood $W$ of $X$ in $M$ where 
$$
(2)\ |\phi_{i,j}|\le\epsilon g^2,\quad(3)\ |\frac{\partial g}{\partial x_k}|\le1,\quad(4)\ |\frac{
\partial\phi_{i,j}}{\partial x_k}|\le\epsilon.
$$
By the \L ojasiewicz inequality, we have $V\subset W$ for large $m$. 
Note that (5) $g\le1/2^{m-1}$ on $V$ since $h\le1/2$. 
Set $c=\max|\frac{d\psi}{dt}|$. 
Then on $V$ 
$$
\gather
|\phi'_{i,j}|=|\phi_{i,j}\psi(\frac{g}{h^m})|\overset(2)\to\le\epsilon g^2\overset(5)\to<\epsilon,\\
|\frac{\partial\phi'_{i,j}}{\partial x_k}|\le|\frac{\partial\phi_{i,j}}{\partial x_k}\psi(\frac{g}{h
^m})|+|\phi_{i,j}\frac{d\psi}{dt}(\frac{g}{h^m})|(|\frac{\partial g}{\partial x_k}|/h^m+m|g\frac{
\partial h}{\partial x_k}|/h^{m+1}),\\
|\frac{\partial\phi_{i,j}}{\partial x_k}\psi(\frac{g}{h^m})|\overset(4)\to\le\epsilon,\\
|\phi_{i,j}\frac{d\psi}{dt}(\frac{g}{h^m})\frac{\partial g}{\partial x_k}|/h^m\overset(3)\to\le\frac{c|
\phi_{i,j}|}{h^m}\overset\text{by def. of $V$}\to\le\frac{2c|\phi_{i,j}|}{g}\overset(2),(5)\to\le c
\epsilon,\\
m|\phi_{i,j}\frac{d\psi}{dt}(\frac{g}{h^m})g\frac{\partial h}{\partial x_k}|/h^{m+1}\overset(1)\to\le
\frac{mc|\phi_{i,j}|g}{h^{m+1}}\overset(2)\to\le2^{\frac{m+1}{m}}mc\epsilon g^{2-\frac{1}{m}}\overset
(5)\to\le2^{4+\frac{1}{m}-2m}mc\epsilon.
\endgather
$$
Hence $\phi'_{i,j}$ is small on $V$ for large $m$. \par
It remains to find $\phi_{i,j}$. 
Let $\Cal K$ denote the sheaf of $\Cal N$-ideals on $M$ defined by $X$. 
Then $\phi$ is a cross-section of $\Cal K^r$ since $\phi$ is $r$-flat at $X$ and since $X$ is normal 
crossing. 
On the other hand, $\sum_{i=1}^n\frac{\partial f}{\partial x_i}\Cal N\supset\Cal K^{r'}$ since $f$ 
has only normal crossing singularities. 
Hence $\phi$ is a cross-section of $\sum_{i,j=1}^n\frac{\partial f}{\partial x_i}\frac{\partial f}
{\partial x_j}\Cal K^{r'}$ because of $r'\ll r$. 
Let $g_l$, for $ l=1,...,k'$, be global generators of $\Cal K^{r'}$ (theorem 2.7). 
Apply theorem 2.8 to the surjective $\Cal N$-homomorphism that assigns
to $(\alpha_{i,j,l})\in \Cal N^{n^2k'}_a \subset \Cal N^{n^2k'}$, for $a\in M$, the value
$$\sum\alpha_{i,j,l}g_{la}(\frac{\partial f}{\partial x_i})_a(\frac{\partial f}{
\partial x_j})_a\in\sum_{i,j=1}^n(\frac{\partial f}{\partial x_i})_a(\frac{\partial f}{\partial x_j}
)_a\Cal K^{r'}_a\subset\sum_{i,j=1}^n\frac{\partial f}{\partial x_i}\frac{\partial f}{\partial x_j}
\Cal K^{r'}.$$ 
Then there exist Nash functions $\phi_{i,j},$ for $ i,j=1,...,n$, in $H^0(M,\Cal K^{r'})$ such that $\phi
=\sum_{i,j=1}^n\phi_{i,j}\frac{\partial f}{\partial x_i}\frac{\partial f}{\partial x_j}$. 
It follows that $\phi_{i,j}$ are $r'$-flat at $X$. \par
The case of compact $M$ is clear by the above arguments. 
\qed
\enddemo
\proclaim{Proposition 4.9} 
(Compactification of a Nash function with only normal crossing singularities)
Let $f$ be a bounded Nash function on a non-compact Nash manifold $M$ with only normal crossing 
singularities. 
Then there exist a compact Nash manifold with corners $M'$ and a Nash diffeomorphism $\pi:M\to\Int M'$ 
such that $f\circ\pi^{-1}$ is extensible to a Nash function on $M'$ with only normal crossing 
singularities. 
\endproclaim
The analytic case does not necessarily hold. \par
We cannot necessarily choose $M'$ with smooth boundary. 
For example any compact Nash manifold with boundary whose interior is Nash diffeomorphic to $M=\R^3$ is Nash 
diffeomorphic to a closed ball in $\R^3$ (Theorem VI.2.2 in
[S$_2$]). But there does not exist Nash function on a 
2-sphere with only normal crossing singularities (see remark (v) after theorem 3.2). \par
Extensibility of a Nash function to a compact Nash manifold with corners is shown in Proposition 
VI.2.8 in [S$_2$]. 
Hence the problem is to impose to the extension to have only normal crossing singularities. 
\demo{Proof of proposition 4.9}
Set $n=\dim M$, and $X=f^{-1}(f(\Sing f))$, set $ B^N=\{x\in\R^N:|x|\le1\}$ for a positive integer $N$ and $S^{N-1
}=\partial B^N$. 
Since there exists a Nash embedding of $M$ into $\R^N$ such that the image is closed in $\R^N$, we can 
assume by lemma 4.6 that $M\subset\Int B^N,$ that $ \overline
M-M\subset S^{N-1}$, that $ \overline M$ is a compact Nash 
manifold with boundary, and moreover $\overline M$ intersects transversally with $S^{N-1}$ in the sense that some 
Nash manifold extension $\tilde M$ of $\overline M$ intersects transversally with $S^{N-1}$, 
that $\overline X$ is a normal crossing Nash subset of $\overline M$, and there exists a Nash function $g$ on 
$\tilde M$ with only normal crossing singularities such that $g(\Sing
g)=f(\Sing f)$, the equality $g^{-1}(g(\Sing g)
)\cap\overline M=\overline X$ holds, and such that $g=f$ on $X$, 
for each $a\in X$, $g(x)-g(a)$ has the same multiplicity as $f(x)-f(a)$ at $a$ and $g(b)>g(a)$ if and 
only if $f(b)>f(a)$ for $b\in M$. 
We do not know whether $g|_M$ is Nash right equivalent to $f$. 
We will modify $\overline M$ and $g$ so that this is indeed the case
and so that $g|_{\overline M}$ has only normal 
crossing singularities. \par
Let $\phi$ be a polynomial function on $\R$ such that $\phi^{-1}(0)=f(\Sing f)$ and $\phi$ is regular 
at $\phi^{-1}(0)$. 
Let $r\in\N$ be large enough. 
Apply lemma 4.7 to $\phi\circ f$ and $\phi\circ g|_M$. 
Then we have a Nash diffeomorphism $\tau_1$ of $M$ such that $\tau_1(X)=X$ and $f\circ\tau_1-g|_M$ is 
$r$-flat at $X$. 
Hence replacing $f$ with $f\circ\tau_1$, we assume that $f-g$ is $r$-flat at $X$. 
Next, by proposition 4.8,(iii) there exists a semialgebraic $C^r$ diffeomorphism $\tau_2$ of $M$ such 
that $g=f\circ\tau_2$ on a semialgebraic neighborhood $V$ of $X$ in $M$ and $\tau_2$ is of class Nash 
on $V$. 
We can choose $V$ of the form $\{x\in M:\phi^2\circ g(x)\le c(x)\xi^m(x)\}$ by the \L ojasiewicz 
inequality, where $\xi(x)=(1-|x|^2)/2$ for $x\in\tilde M$, where $c$ is a positive Nash function on $\tilde M$ 
such that $c$ depends on only $|x|$ and $m$ is a large odd integer. 
Shrink $\tilde M$ so that $\xi<0$ on $\tilde M-\overline M$. 
We can choose, moreover, $c$ and $m$ so that $\phi^2\circ g-c\xi^m$ is regular at $A-S^{N-1
}$, where $A$ denotes the zero set of $\phi^2\circ g-c\xi^m$, and
hence $V$ is a Nash manifold with boundary $\{x\in M:\phi^2\circ
g(x)=c(x)\xi^m(x)\}$. Actually, let $0<\epsilon_0\in\R$ be small. 
Then for any $0<\epsilon\in\R$ with $\epsilon<\epsilon_0$, $\xi^{-1}(\epsilon)\cup(\phi\circ g)^{-1}(0)$ 
is normal crossing in $\tilde M$, and hence for small $c$ and large
$m$, the function $\phi^2\circ g$ on $\{x\in\xi^{
-1}(\epsilon):0<\phi^2\circ g(x)<2c(x)\xi^m\}$ is regular. 
We can choose $c$ and $m$ independently of $\epsilon$. 
Therefore, $\phi^2\circ g-c\xi^m$ is regular at $A\cap\xi^{-1}((0,\,\epsilon_0))$. 
Moreover, if we choose $c$ and $m$ so that $c\xi^m$ is close to a small constant on $M-\xi^{-1}
((0,\,\epsilon_0/2])$, then $\phi^2\circ g-c\xi^m$ is regular at $A-\xi^{-1}((0,\,\epsilon_0/2])$. 
Hence $\phi^2\circ g-c\xi^m$ can be regular at $A-S^{N-1}$. 
However, we omit $c$ for simplicity of notation. 
We want first to modify $M$ so that $\overline V$ is a neighborhood of $\overline X$ in $\overline M$. 
\par
Apply theorem 2.10 to the two sheaves of $\Cal N$-ideals on $\tilde M$ defined by $(\phi\circ g)^{-1}(0
)$ and generated by $\xi\cdot(\phi^2\circ g-\xi^m)$. Note that the former sheaf is normal crossing, the stalk 
of the latter is not generated by one regular function germ at a point of $\overline X-X$ only, and at least 
one of the two stalks of both sheaves at each $x\not\in\overline X-X$ is $\Cal N_x$. 
Then we have a composition of a finite sequence of blowings-up $\tau_3:\hat M\to\tilde M$ along smooth Nash 
centers such that $\tau_3|_{\tau^{-1}_3(\tilde M-(\overline X-X))}:\tau_3^{-1}(\tilde M-(\overline X-X))\to
\tilde M-(\overline X-X)$ is a Nash diffeomorphism and $(\phi\circ g\cdot\xi\cdot(\phi^2\circ g-\xi^m))
\circ\tau_3$ has only normal crossing singularities at its zero set, say $Y$. 
It follows that $(\hat M,Y,\overline{\tau^{-1}_3(M)})$ is Nash diffeomorphic to $(\R^n,\{(x_1,...,x_n)\in\R^n:x_1\cdots x_{n'}=0\},B)$ locally at each point of $\overline{\tau^{-1}_3(M)}$ for some $n'\,(\le n)\in\N$, where $B$ denotes the closure of the union of some connected components of 
$\{x_1\cdots x_{n'}\not=0\})$, and 
$\overline{\tau_3^{-1}(M)}-\tau^{-1}_3(\overline X-X)$ is a Nash manifold with boundary. 
However, $\overline{\tau^{-1}_3(M)}$ is not necessarily a manifold with corners. 
It may happens that $\tau^{-1}_3(M)$ is locally diffeomorphic to the union of more than one 
connected components of $\{(x_1,...,x_n)\in\R^n:x_1\cdots x_{n'}\not=0\}$ at some point of $\tau^{-1}
_3(\overline X-X)$, for $0<n'\,(\le n)\in\N$. 
Then we need to separate these connected components. 
That is possible as shown in the proof of Theorem VI.2.1 in [S$_2$]. 
Namely, there exist a compact Nash manifold $L$ with corners and a Nash immersion $\tau_4:L\to
\overline{\tau^{-1}_3(M)}$ such that $\tau_4|_{L-\Sing\partial L}$ is a Nash diffeomorphism to its 
image and the image contains $\overline{\tau^{-1}_3(M)}-\tau^{-1}_3(\overline X-X)\,(\supset\tau_3^{-1
}(M))$. \par
Clearly $(\phi\circ g\cdot\xi\cdot(\phi^2\circ g-\xi^m))\circ\tau_3\circ\tau_4$ has only normal 
crossing singularities at its zero set $\tau^{-1}_4(Y)$ since $\tau_4$ is an immersion. 
Set $\tau=\tau_2\circ\tau_3\circ\tau_4|_{\Int L}$ and
$h=g\circ\tau_3\circ\tau_4$. Define $W=(\tau_3\circ\tau_4
)^{-1}(V)$ and $W'=\overline W-\overline{\partial W
}$ and set $Z=(\tau_3\circ\tau_4)^{-1}(X)$. 
Then $W$ is a non-compact Nash manifold with boundary; $\tau$ is a semialgebraic $C^r$ diffeomorphism from $\Int L$ to $M$ and of class Nash on $W$; $h$ 
is a Nash function on $L$; $h=f\circ\tau$ on $W$; $h$ is regular on $\Int L-Z$; $h|_{\Int L\cup W'}$ 
has only normal crossing singularities at $\overline Z$ though $h$ is not necessarily so globally; 
$\overline W$ is a 
neighborhood of $\overline Z$ in $L$ because if it were not, $\overline Z\cap\overline{(\tau_3\circ\tau
_4)^{-1}(\{x\in M:\phi^2\circ g(x)=\xi^m(x)\})}$ could be not empty and of dimension $n-2$ but 
contained in $(\tau_3\circ\tau_4)^{-1}(\xi^{-1}(0))$, which contradicts the normal crossing property 
of $(\phi\circ g\cdot\xi\cdot(\phi^2\circ g-\xi^m))\circ\tau_3\circ\tau_4$. 
Note that $W'$ and $\overline W$ are Nash manifolds with corners by the next fact and the normal crossing 
property of $(\xi\cdot(\phi^2\circ g-\xi^m))\circ\tau_3\circ\tau_4$. 
Thus $\overline V$ is changed to $\overline W$---a neighborhood of $\overline Z$ in $L$. 
We consider $h$ on $L$ in place of $g$ on $\overline M$. \par
We replace $\tau$ by a Nash diffeomorphism. 
Let $0\ll r\in\N$, set $\overline\psi=(\phi^r\circ h\cdot\xi^r)\circ\tau_3\circ\tau_4$ on $L$ and 
$\psi=\overline\psi|_{\Int L}$, and let $\Cal I$ denote the sheaf of $\Cal N$-ideals on $\Int L$ 
generated by $\psi$. 
Then we regard $\tau$ as an element of $H^0(\Int L,\Cal N/\Cal I)^N$ because $\supp\Cal N/\Cal I=Z$ 
and $\tau$ is of class Nash near there. 
Hence by theorem 2.8 there exists a Nash map $\tau':\Int L\to\R^N$ such that $\tau-\tau'=\psi\theta$ 
for some semialgebraic $C^r$ map $\theta:\Int L\to\R^N$ of class Nash on $W$. 
Approximate $\theta$ by a Nash map $\theta':\Int L\to\R^N$ in the semialgebraic $C^r$ topology, and 
set $\tau''=p\circ(\tau'+\psi\theta')$, where $p$ denotes the orthogonal projection of a semialgebraic 
tubular neighborhood of $M$ in $\R^N$. 
Then $\tau''$ is a well-defined Nash diffeomorphism from $\Int L$ to $M$ and close to $\tau$ in the 
semialgebraic $C^r$ topology; $f\circ\tau''-h|_{\Int L}=\psi\delta$ for some semialgebraic $C^r$ 
function $\delta$ on $\Int L$ though $f\circ\tau''-h|_{\Int L}$ does not necessarily vanish on $W$; 
moreover, $\delta$ is extensible to a semialgebraic $C^{r'}$ function $\overline\delta$ on 
$\Int L\cup W'$ for $0\ll r'\,(\ll r)\in\N$ by the definition of the semialgebraic $C^{r'}$ topology, 
by the fact that a small semialgebraic $C^r$ function on $\Int L$ is extensible to a semialgebraic 
$C^r$ function on $L$ and by 
$$
f\circ\tau''-h|_{\Int L}=f\circ p\circ(\tau+\psi\cdot(\theta'-\theta))-f\circ p\circ\tau\quad\text{
on}\ W. 
$$
The last equality implies also that $\delta$ is of class Nash on $W$, and hence on $\Int L$ since 
$f\circ\tau''$ and $h$ are Nash functions and $\psi^{-1}(0)\subset W$. \par
Next we modify $h$. 
Let $\overline\delta'$ be a Nash approximation on $\Int L\cup W'$ of $\overline\delta$ in the 
semialgebraic $C^{r'}$ topology, and set $\delta'=\overline\delta'|_{\Int L}$ and $h'=h+\psi
\overline\delta'$ on $\Int L\cup W'$. 
Then $h'$ is a Nash function on $\Int L\cup W'$ and has only normal crossing singularities at 
$\overline Z$ by the same property of $h|_{\Int L\cup W'}$ and by the definition of $h'$, and $f
\circ\tau''-h'|_{\Int L}$ is of the form $\psi\cdot(\delta-\delta')$. 
Hence $f$ and $h'\circ\tau^{\prime\prime-1}$ satisfy the conditions in proposition 4.8,(ii) because 
$\phi^r\circ h\circ\tau^{\prime\prime-1}$ is of the form $\sum_{i,j=1}^k\psi_{i,j}v_f\cdot v_jf$ 
for some Nash functions $\psi_{i,j}$ on $M$ and Nash vector fields
$v_i$, for $ i=1,...,k$, on $M$ which 
span the tangent space of $M$ at each point of $M$ and because $\xi^r|_M\cdot(\delta-\delta')\circ
\tau^{\prime\prime-1}$ is small as a semialgebraic $C^{r'}$ function on $M$. 
Consequently, $f$ and $h'\circ\tau^{\prime\prime-1}$ are Nash right equivalent, and we can replace 
$f$ with $h'|_{\Int L}$. \par
We can assume that $W'\cap\partial L$ is the union of some connected components $\sigma$ of strata of the 
canonical stratification $\{L_i\}$ of $\partial L$ such that $\overline\sigma\cap\overline Z\not=
\emptyset$. Actually, let $\psi_L$ be a non-negative Nash function on
$L$ with zero set $\overline Z$, and let $\epsilon>
0\in\R$ be such that the restriction of $\psi_L$ to $\psi_L^{-1}((0,\,2\epsilon))$ is regular. 
Then $\psi_L^{-1}(\epsilon)$ is a compact Nash manifold with corners
equal to $\partial L\cap\psi^{-1}_L(
\epsilon)$. 
Let $\{L_{\epsilon,i}\}$ denote the canonical stratification of $\partial L\cap\psi_L^{-1}(
\epsilon)$. 
We blow up $L_{\epsilon,i}$ as follows. 
Let $L'$ and $\tilde L'$ be a compact Nash submanifold possibly with corners of $L$ and some Nash manifold 
extension of $L'$ respectively. 
If $\tilde L'\cap L=L'$ and $(L,\tilde L')$ is locally diffeomorphic 
to $(\{(x_1,...,x_n)\in\R^n:x_1\ge 0,...,x_{n'}\ge0\},\{x_{n_1}=\cdots=x_{n_k}=0\})$ for some $n'\,
(\le n),1\le n_1<\cdots<n_k\le n\in\N$, then we say $L'$ has the property $(*)$. 
For $L'$ with $(*)$, consider $\gamma:\Gamma\to L$---the restriction of the blowing-up of a small 
Nash manifold extension $\tilde L$ of $L$ along center $\tilde L\cap\tilde L'$ to the closure of 
inverse image of $L-L'$, modify $\gamma:\Gamma\to L$ so that $\Gamma$ is a compact Nash manifold with 
corners by the idea in the proof of Theorem VI.2.1 in [S$_2$] as before, use the same notation $\gamma:
\Gamma\to L$, and call it the $(*)$-blowing-up of $L$ along center $L'$. 
Note that $\gamma^{-1}(L')$ is the closure of the union of some connected components of $\Reg\partial\Gamma$. 
Set $\Gamma_{-1}=L$ and let $0\le k\le n-2$. 
Inductively we define $(*)$-blowing-up $\gamma_k:\Gamma_k\to\Gamma_{k-1}$ of $\Gamma_{k-1}$ along 
center $L_{\epsilon,0}$ if $k=0$ and along center $\overline{(\gamma_0\circ\cdots\circ\gamma_{k-1})^{-1}(L
_{\epsilon,k})}$ if $k>0$, which is possible because $L_{\epsilon,0}$ and $\overline{(\gamma_0\circ
\cdots\circ\gamma_k)^{-1}(L_{\epsilon,k+1})}$ for $0\le k\le n-3$ are compact Nash submanifolds with 
corners of $\Gamma_{-1}$ and $\Gamma_k$ with $(*)$, respectively. 
Thus we assume that the above condition on $W$ holds considering $(\Gamma_{n-2},(\psi_L\circ\gamma_0
\circ\cdots\gamma_{n-2})^{-1}([0,\,\epsilon])-\partial\Gamma_{n-2})$ in place of $(L,W)$. 
Here we choose $\epsilon$ so small that $(**)$ $h'$ is extensible to a Nash function on an open 
semialgebraic neighborhood of $\Int L\cup\overline W$ in $L$ with only normal crossing singularities. 
\par
Moreover, we can assume that the closure of each connected component of $\Reg\partial L$ is a Nash 
manifold possibly with corners. 
Indeed, we obtain this situation if we repeat the same arguments as above to the canonical stratification of 
$\partial L$ compatible with $\{x\in\partial L:\dis(x,L_k)=\epsilon_k,\,\dis(x,L_i)\ge\epsilon_i,
\ i=0,...,k-1\}$, for $k=0,...,n-2$. 
Here we naturally define the canonical stratification of $\partial L$ 
compatible with the above family in the same way as in the remark after the proof of lemma 4.5. 
After this modification of $L$, the property $(**)$ continues to hold. \par
Let $M_j$, for $ j\in J,$ be the set of closures of the connected components of $\Reg\partial L$, and let $J_0$ 
denote the subset of $J$ consisting of $j$ such that $M_j\cap\overline Z=\emptyset$. 
Let $\tilde L$ and $\tilde M_j$ be Nash manifold extensions of $L$ and $M_j$, respectively, which 
are contained and closed in a small open semialgebraic neighborhood $U$ of $L$ in the ambient 
Euclidean space such that $\cup_{j\in J}\tilde M_j$ is normal crossing in $\tilde L$ and for each $j
\in J$ there is one and only one connected component of $\tilde L-\tilde M_j$ which does not 
intersect with $L$. 
Let $\tilde Z$ denote the smallest Nash subset of $\tilde L$ containing $Z$. 
Then $\tilde Z$ is normal crossing in $\tilde L$, and there exist Nash functions $\chi_j$ on 
$\tilde L$ with zero set $\tilde M_j$, regular there and with $\chi_j>0$ on $\Int L$. \par
By $(**)$ we can choose a sufficiently small $U$ so that $h'$ can be extended to a Nash function $h'_+$ on 
$L_+\overset\text{def}\to=\{x\in\tilde L:\chi_j(x)>0,\ j\in J_0\}$,
such that $h'_+(\Sing h'_+)=h'(Z)$ and $h'_+$ 
has only normal crossing singularities. 
Now we smooth $h'_+$ at $\overline{L_+}-L_+$ as in the proof of Proposition VI.2.8 in [S$_2$]. 
Let $\tilde L\subset\R^N$, set $G=\graph h'_+\subset L_+\times\R$, and let $G^Z$ be the Zariski closure of 
$G$ in $\R^N\times\R$ and $Q$ be the normalization of $G^Z$ in $\R^N\times\R\times\R^{N'}$ for some 
$N'\in\N$, and let $r:Q\to\R^N\times\R$ and $q:Q\to\R^N$ denote the restrictions to $Q$ of the 
projections $\R^N\times\R\times\R^{N'}\to\R^N\times\R$ and $\R^N\times\R\times\R^{N'}\to\R^N$, 
respectively. 
Then it is known that $r$ is a proper map to $G^Z$, and by Artin-Mazur Theorem there exists a union 
of connected components $R$ of $Q-r^{-1}(\overline G-G)$ such that $R\subset\Reg Q$ and $r|_R$ is a 
Nash diffeomorphism onto $G$. 
Here we can replace $r^{-1}(\overline G-G)$ with a Nash subset $q^{-1}((\prod_{j\in J_0}\chi_j)^{-1}
(0))$ of $Q$ because $r^{-1}(\overline G-G)\subset q^{-1}((\prod_{j\in J_0}\chi_j)^{-1}(0))$ and $R
\cap q^{-1}((\prod_{j\in J_0}\chi_j)^{-1}(0))=\emptyset$; $q|_R$ is a Nash diffeomorphism onto $L_+$; 
the map $h'_+\circ q|_R$ is the restriction of the projection $\R^N\times\R\times\R^{N'}\to\R$ and hence 
extensible to a smooth rational function on $Q$; the set $\overline{R\cap q^{-1}(L)}$ is compact because $r$ 
is proper and because $\overline{G\cap L\times\R}$ is compact by
boundedness of $f$; the function $h'_+\circ 
q|_R$ has only normal crossing singularities because the same is true for $h'_+$. 
However, $\chi_j\circ q$ are now not necessarily regular at their zero sets. 
By theorem 2.8, $R\cap q^{-1}(\tilde Z)$ is a Nash subset of $\Reg Q$ and there exists a Nash function 
$\alpha$ on $\Reg Q$ whose zero set is $R\cap q^{-1}(\tilde Z)$ and which has only normal crossing 
singularities there since $R\cap q^{-1}(\tilde Z)$ is a Nash subset of $R$ and since its closure in 
$\Reg Q$ does not intersect with $\overline R-R$. \par
Thus replacing $\tilde L,\ L_+,\ h'_+,\ \chi_j$ and $\tilde Z$ with $\Reg Q,\ R,\ h'_+\circ q|_R,
\ \chi_j\circ q|_{\Reg Q}$ and $R\cap q^{-1}(\tilde Z)$ we assume from the beginning that $M$ and $f$ 
satisfy moreover the following conditions. \par
(i) $\tilde f$ and $\chi_j$, for $ j\in J,$ are a finite number of Nash functions on a Nash manifold 
$\tilde M$, and $\tilde X$ is a normal crossing Nash subset of $\tilde M$. \newline
(ii) $M$ is the union of some connected components of $\tilde
M-(\prod_{j\in J}\chi_j)^{-1}(0)$, the set $\overline M$ is compact,
the equalities $f=\tilde f|_M$ and $X=\tilde X\cap M$ hold (we do not assume that $\overline M$ is a 
manifold with corners). \par
We make $\prod_{j\in J}\chi_j$ normal crossing at its zero set. 
Apply theorem 2.10 to the sheaf of $\Cal N$-ideals on $\tilde M$ defined by $\tilde X$ and the sheaf 
of $\Cal N$-ideals $\prod_{j\in J}\chi_j\Cal N$. 
Then via blowings-up, $\prod_{j\in J}\chi_j$ becomes to have only normal crossing singularities at its 
zero set, and the conditions (i) and (ii) do not change because the subset of $\tilde M$ where we 
modify by blowings-up is contained in $(\prod_{j\in J}\chi_j)^{-1}(0)$. \par
It remains to make $\tilde f$ together with $(\prod_{j\in J}\chi_j)^{-1}(0)$ normal crossing. 
Let $\{\tilde M_i\}$ denote the canonical stratification of $(\prod_{j\in J}\chi_j)^{-1}(0)$, set 
$\tilde M_n\!=\!\tilde M-(\prod_{j\in J}\chi_j)^{-1}(0)$, and let $\tilde\phi$ be a polynomial 
function on $\R$ such that $\tilde\phi^{-1}(0)=\cup_{i=0}^n\tilde f(\Sing\tilde f|_{\tilde M_i})$ 
and $\tilde\phi$ is regular at $\tilde\phi^{-1}(0)$. 
Once more, apply theorem 2.10 to the sheaf of $\Cal N$-ideals on $\tilde M$ defined by $\tilde X\cup(
\prod_{j\in J}\chi_j)^{-1}(0)$ and the sheaf of $\Cal N$-ideals $[\tilde\phi\circ\tilde f\Cal N:\cap_
i\Cal I_i^{\alpha_i}]\overset\text{def}\to=\cup_{x\in\tilde M}\{\rho\in\Cal N_x:\rho\cap_i\Cal I^{
\alpha_i}_{ix}\subset\tilde\phi\circ\tilde f\Cal N_x\}$, where $\cap_i\Cal I_i$ is the decomposition 
of the sheaf of $\Cal N$-ideals on $\tilde M$ defined by $\tilde X$ to irreducible finite sheaves of 
$\Cal N$-ideals and each $\alpha_i$ is the maximal integer such that $\tilde\phi\circ\tilde f\Cal N$ 
is divisible by $\Cal I_i^{\alpha_i}$. 
Then $(\tilde f-\tilde f(x_0))\prod_{j\in J}\chi_j$ becomes to have only normal crossing 
singularities at its zero set for each $x_0\in\tilde M$ and the subset of $\tilde M$ where we modify 
now by blowings-up does not intersect with $M$ because the stalk of the latter sheaf at each point 
of $M$ is generated by a regular function germ and because 
$$(\tilde X\cup(\prod_{j\in J}\chi_j)^{-1}
(0))\cap\supp\Cal N/[\tilde\phi\circ\tilde f\Cal N:\cap_i\Cal I_i^{\alpha_i}]\cap M=\emptyset.$$ 
Finally, we separate as before $M$ at the points of $\overline M$ where $M$ is not locally connected 
so that $\overline M$ is a compact Nash manifold with corners. 
Then $\tilde f|_{\overline M}$ has only normal crossing singularities, and we complete the proof. 
\qed
\enddemo

\head 
5. Proofs of theorem 3.2 and theorems 3.1,(2) and 3.1,(3)
\endhead

\subhead
5.1. Proof of theorem 3.2
\endsubhead\par

By proposition 4.1 it suffices to prove the Nash case and, moreover, that the cardinality of Nash R-L 
equivalence classes of Nash functions with only normal crossing singularities on a compact Nash 
manifold possibly with corners is zero or countable. 
The reasons are that first we can restrict functions to being bounded by the fact that $\R$ is Nash diffeomorphic to $(0,\,1)$ and secondly by proposition 4.9 
we can regard a non-compact Nash manifold $M$ and a bounded Nash function $f$ with only normal crossing singularities on $M$ as the interior of a compact Nash manifold with corners $M'$ and the 
restriction to $M$ of a Nash function on $M'$ with only normal crossing singularities. 
Assume that there is at least one Nash function $f$ on $M$ with only normal crossing singularities. 
Then the cardinality is infinite because we can increase arbitrarily the cardinality of the critical 
value set, which is finite, by replacing $f$ with $\pi\circ f$ for some Nash function $\pi$ on $\R$. 
Let $\{X_\alpha\}_{\alpha\in A}$ denote all normal crossing Nash subsets of $M$. 
We define $\alpha$ and $\alpha'$ in $A$ to be equivalent if there exists a Nash diffeomorphism of 
$M$ which carries $X_\alpha$ to $X_{\alpha'}$. 
Then by lemma 4.4 the cardinality of equivalence classes of $A$ is countable. 
Hence it suffices to see that for each $X_\alpha$ there exist at most a countable number of Nash R-L 
equivalence classes of Nash functions $f$ on $M$ with only normal crossing singularities such that 
$f^{-1}(f(\Sing f))=X_\alpha$. 
Let $F_\alpha$ denote all such Nash functions. 
Clearly there are a finite number of equivalence classes of $\{f|_{X_\alpha}:X_\alpha\to\R:f\in F_
\alpha\}$ under the Nash left equivalence relation since the value sets are finite. 
Moreover, there are at most a countable number of choices of multiplicity of $f-f(a)$ at $a$ for $f
\in F_\alpha$ and $a\in X_\alpha$. 
Hence we reduce the problem to the following one. 
Fix $f\in F_\alpha$, and let $F_f$ denote the family of $g\in F_\alpha$ such that $g=f$ on $X_\alpha$ 
and $g-g(a)$ has the same multiplicity as $f-f(a)$ at each point $a$ of $X_\alpha$. 
Then the cardinality of Nash right equivalence classes of functions in $F_f$ is finite. 
Moreover, it suffices to prove that each element of $F_f$, say $f$, is stable in $F_f$ in the sense that 
any $g\in F_f$ near $f$ in the $C^\infty$ topology is Nash right equivalent to $f$ because there are 
only a finite number of connected components in $F_f$. \par
Set $n=\dim M$, embed $M$ in $\R^N$, and let $\{M_i\}$ denote the canonical stratification of $M$. 
There exist Nash vector fields $v_1,...,v_k$ on $M$ such that $v_{1x},...,v_{kx
}$ span the tangent space $T_xM_i$ of $M_i$ at each $x\in M_i$. If we regard $M$ as $\{(x_1,...,x_n)\in\R^n:x
_1\ge0,...,x_{n'}\ge0\}$ by its local coordinate system, then $x_i\frac
\partial{\partial x_i}$ is contained in the linear space over $N(M)$
spanned by $v_1,...,v_k$ for each $1\le i\le n'$. Actually, set
$L_i=\cup_{j=0}^iM_j$, for $ i=0,...,n-1$, and choose a Nash manifold extension $\tilde M$ of $M$ 
and Nash subset extensions $\tilde L_i$ of $L_i$ in $\tilde M$ so that $\tilde L_{n-1}$ is normal 
crossing in $\tilde M$ and $\{\tilde L_i-\tilde L_{i-1}\}$ is the canonical stratification of 
$\tilde L_{n-1}$. 
Set $L_n=M$ and $\tilde L_n=\tilde M$ also. 
Then when we describe $(\tilde M,\tilde L_{n-1})$ by a local coordinate system as $(*)$ $(\R^n,\{(x_
1,...,x_n)\in\R^n:x_1\cdots x_{n'}=0\})$, 
$$
\tilde L_i=\bigcup_{1\le j_1<\cdots<j_{n-i}\le n'}\{(x_1,...,x_n)\in\R^n:x_{j_1}=\cdots=x_{j_{n-i}}
=0\},\ n-n'\le i\le n.
$$
We consider the situation on $\tilde M$ rather than on $M$ because the
existence of $v_1,...,v_k$ follows from the
existence of Nash vector fields on $\tilde M$ with the same properties. \par
First, let $w_{n,1},...,w_{n,k_n}$ be Nash vector fields on $\tilde M$ which span the tangent 
space of $\tilde M$ at each point, and $\alpha_n$ a global generator of the sheaf of $\Cal N
$-ideals on $\tilde M$ defined by $\tilde L_{n-1}$---we can choose $\tilde M$ so that $\alpha_n$ 
exists because $M$ is a manifold with corners. 
Then $v_{n,1}=\alpha_nw_{n,1},...,v_{n,k_n}=\alpha_nw_{n,k_n}$ are Nash vector fields on $\tilde M$, 
span the tangent space of $\tilde M$ at each point of $\tilde M-\tilde L_{n-1}$ and vanish at 
$\tilde L_{n-1}$, and in the case $(*)$, for each $1\le i\le n'$, $x_i\frac\partial{\partial x_i}$ 
on $\{(x_1,...,x_n)\in\R^n:x_j\not=0$ for $1\le j\le n'$ with $j\not=i\}$ is contained in the linear 
space over the Nash function ring on the set spanned by $v_{n,1},...,v_{n,k_n}$. \par
Next fix $i<n$ and consider on $\tilde L_i$. 
Then it suffices to prove the following two statements.

(i) There exist Nash vector fields $v_{i,1},...,v_{i,k_i}$ on $\tilde L_i$---Nash cross-sections 
of the restrictions to $\tilde L_i$ of the tangent bundle of $\R^N$, i.e. the restrictions to 
$\tilde L_i$ of Nash vector fields on $\R^N$ by theorem 2.8---which span the tangent space of 
$\tilde L_i-\tilde L_{i-1}$ at its each point and vanish at $\tilde L_{i-1}$ and such that in the 
case of $(*)$ the condition on each irreducible component $\{(x_1,...,x_n)\in\R^n:x_{j_1}=\cdots=x_{
j_{n-i}}=0\}$, same as on $\tilde M$, is satisfied for $1\le j_1<\cdots<j_{n-i}\le n'$; to be precise, 
for any $1\le j\le n'$ other than $j_1,...,j_{n-i}$, then
$x_j\frac\partial{\partial x_j}$ on 
$$\{(x_1,...
,x_n)\in\R^n:x_{j_1}=\cdots=x_{j_{n-i}}=0,\ x_l\not=0
\ \text{if}\ l\in \{1,...,n'\}\setminus\{j_1,...,j_{n-i}
,j\}\}$$ 
is contained in the linear space over the Nash function ring on the set spanned by $v_{i,1},...,v_
{i,k_i}$.

(ii) Any Nash vector field on $\tilde L_i$ tangent to $\tilde L_j-\tilde L_{j-1}$ at its each point 
for $j\le i$ is extensible to a Nash vector field on $\tilde L_{i+1}$ tangent to $\tilde L_{i+1}-
\tilde L_i$ at each its point. 
\par
Proof of (i). 
By considering the Zariski closure of $\tilde L_i$ and its normalization and by Artin-Mazur Theorem, 
we have a Nash manifold $P_i$ and a Nash immersion $\xi_i:P_i\to\tilde L_i$ such that $\xi_i|_{P_i-
\xi_i^{-1}(\tilde L_{i-1})}$ is a Nash diffeomorphism onto $\tilde L_i-\tilde L_{i-1}$. 
Note that $\xi_i^{-1}(\tilde L_{i-1})$ is normal crossing in $P_i$. 
Apply the same arguments to $(P_i,\xi_i^{-1}(\tilde L_{i-1}))$ as on $(\tilde M,\tilde L_{n-1})$. 
Here the difference is only that we need a finite number of global generators $\alpha_{i,1},\alpha_{i,2},...$ 
of the sheaf of $\Cal N$-ideals on $P_i$ defined by $\xi_i^{-1}(\tilde L_{i-1})$. 
Then there exist Nash vector fields $w_{i,1},...,w_{i,k_i}$ on $P_i$ with the corresponding 
properties, and they induce semialgebraic $C^0$ vector fields $v_{i,1},...,v_{i,k_i}$ on $\tilde L
_i$ through $\xi_i$ because $w_{i,1},...,w_{i,k_i}$ vanish on $\xi_i^{-1}(\tilde L_{i-1})$. 
Such $v_{i,1},...,v_{i,k_i}$ are of class Nash by the normal crossing property of 
$\tilde L_{n-1}$ in $\tilde M$ and satisfy the conditions in (i). \par
Proof of (ii). 
Let $v$ be a Nash vector field on $\tilde L_i$ in (ii), and $\xi_{i+1}:P_{i+1}\to\tilde L_{i+1}$ 
the same as above. 
Then since $\xi_{i+1}$ is an immersion, $v$ pulls back a Nash cross-section $w$ of the restriction 
to $\xi_{i+1}^{-1}(\tilde L_i)$ of the tangent bundle of the Nash manifold $P_{i+1}$, and by theorem 
2.8 we obtain a Nash vector field on $P_{i+1}$ whose restriction to $\xi_{i+1}^{-1}(\tilde L_i)$ 
is $w$. 
This vector field induces a Nash vector field of $\tilde L_{i+1}$ through $\xi_{i+1}$, which is an 
extension of $v$, by the same reason as in the proof of (i). \par
Let $g\in F_f$ near $f$. 
It suffices to see that $f$ and $g$ are $C^\omega$ right
equivalent. Actually, assume that there exists a $C^\omega$ diffeomorphism $\pi$ of $M$ such that $f=g\circ\pi$. 
Let $\tilde M$ and $\tilde L_{n-1}$ be the same as above and so small that $f$ and $g$ are 
extensible to Nash functions $\tilde f$ and $\tilde g$ on $\tilde M$, respectively. 
Extend $\pi$ to a $C^\omega$ diffeomorphism $\tilde\pi:U_1\to U_2$ between open neighborhoods of 
$M$ in $\tilde M$ so that $\tilde\pi(U_1\cap\tilde L_{n-1})\subset\tilde L_{n-1}$ and $\tilde f=
\tilde g\circ\pi$. 
As above, let $\alpha_n$ be a global generator of the sheaf of $\Cal N$-ideals on $\tilde M$ 
defined by $\tilde L_{n-1}$. 
Then $\alpha_n\circ\tilde\pi=\beta\alpha_n$ on $U_1$ for some positive $C^\omega$ function $\beta$ 
on $U_1$. 
Consider the following equations in variables $(x,y,z)\in\tilde M^2\times\R$. 
$$
f(x)-g(y)=0\quad\text{and}\quad\alpha_n(y)-z\alpha_n(x)=0 
$$
Here the second equation means that if $x\in\tilde L_{n-1}$ then $y\in\tilde L_{n-1}$. 
Then $y=\tilde\pi(x)$ and $z=\beta(x)$ are $C^\omega$ solutions. 
Hence by Nash Approximation Theorem II, there exist Nash germ $M$ $y=\pi'(x)$ and $z=\beta'(x)$ solutions on $M$,
which are approximations of the germs of $\tilde\pi$ and $\beta$ on $M$. 
Thus $\pi'|_M$ is a Nash diffeomorphism of $M$ and $f=g\circ\pi'$ on $M$. \par
Now we show the $C^\omega$ right equivalence of $f$ and $g$. 
Set $G(x,t)=(1-t)f(x)+tg(x)$ for $(x,t)\in M\times[0,\,1]$. 
Then $G(x,0)=f(x)$ and $G(x,1)=g(x)$. 
Hence by the same reason as in the proof of theorem 3.1,(1) it suffices to find a $C^\omega$ vector field 
$v$ on $M\times[0,\,1]$ of the form $\frac\partial{\partial t}+\sum_{i=1}^ka_iv_i$ for some $C^
\omega$ functions $a_i$ on $M\times[0,\,1]$ such that $vG=0$ on $M\times[0,\,1]$, i.e.,
$$
f-g=\sum_{i=1}^ka_i(v_if+tv_i(g-f)).\tag **
$$
Moreover, as shown there, we only need to solve this equations locally at each point $(x_0,t_0)$ of 
$M\times[0,\,1]$ since $M$ is compact. \par
If $x_0\not\in X_\alpha$, then $(v_if)(x_0)\not=0$ for some $i$ and hence we have solutions of $(**)$ $a_j=0$ for $j\not=i$ and 
$a_i=(f-g)/(v_if+tv_i(g-f))$ around $(x_0,t_0)$ because $g-f$ and hence $tv_i(g-f)$ are small in the $C^\infty$ topology. \par
Let $x_0\in X_\alpha$. 
Then we can assume that $M=\{x=(x_1,...,x_n)\in\R^n:|x|\le 1,\,x_1\ge0,...,x_{n'}\ge0\}$ for some $n'\,
(\le n)\in\N$, that $x_0=0$ and $ f(x)=x^\beta$ for some $\beta=(\beta_1,...,\beta_n)\in\N^n$ with $|\beta|>
0$, that $ k=n$ and $ v_1=x_1\frac\partial{\partial x_1},...,v_{n'}=x_{n'}\frac\partial{\partial x_{n'}},v_
{n'+1}=\frac\partial{\partial x_{n'+1}},...,v_n=\frac\partial{\partial
x_n}$ and that $f-g=bx^\beta$ 
for some small $C^\omega$ function $b$ on $M$ by lemma 2.12. 
Let $i$ be such that $\beta_i\not=0$. 
Then $v_if=\beta_ix^\beta/x_i$ and $v_i(f-g)=b\beta_ix^\beta/x_i+\frac{\partial b}{\partial x_i}x
^\beta$ if $i>n'$, and $v_if=\beta_ix^\beta$ and $v_i(f-g)=b\beta_ix^\beta+x_i\frac{\partial b}{
\partial x_i}x^\beta$ if $i\le n'$. 
In any case $(**)$ is solved as before. 
Thus theorem 3.2 is proved. 
%%%%%%%%%%%%%%%%%%%%%%%%%%%%%%%%%%%%%%%%%%%%%%%%%%%%%%%%%%%%%%%%%%%%%%%%%%%%%%%%%%%
\subhead
5.2. Proof of theorems 3.1,(2) and 3.1,(3)
\endsubhead\par
Let us consider the case where $M$ is a manifold without corners. \par
{\it Proof of (2).} 
Set $X=f^{-1}(f(\Sing f))$ and $Y=g^{-1}(g(\Sing g))$, and let $\pi$ be a $C^2$ diffeomorphism of 
$M$ such that $f\circ\pi=g$. 
Then $X$ and $Y$ are normal crossing, $\pi(Y)=X$, and we assume that $\pi$ is close to id in the Whitney 
$C^2$ topology by replacing $f$ and $\pi$ with $f\circ\pi'$ and $\pi^{\prime-1}\circ\pi$ for a $C^
\infty$ approximation $\pi'$ of $\pi$ in the Whitney $C^2$ topology. 
Hence by lemma 4.2 and properness of $f$ and $g$, there exists a $C^\infty$ 
diffeomorphism $\pi''$ of $M$ close to id in the Whitney $C^2$ topology such that $\pi''(Y)=X$. 
Replace $f$ and $\pi$, once more, by $f\circ\pi''$ and $\pi^{\prime\prime-1}\circ\pi$. 
Then we can assume that, moreover, $X=Y$. 
We want to modify $\pi$ to be of class $C^\infty$ on a neighborhood of $X$. 
Set $B(\epsilon)=\{x\in\R^n:|x|\le\epsilon\}$ for $\epsilon>0\in\R$. 
Let $\{U_i\}$ and $\{U'_i\}$ be locally finite open coverings of $X$ in $M$ such that $\overline{U'
_i}\subset U_i$, such that $ \pi(\overline{U'_i})\subset U_i$, each $f|_{U_i}$ is $C^\infty$ right equivalent to 
the function $\prod_{j=1}^nx_j^{\alpha_j}+\,$constant, for
$x=(x_1,...,x_n)\in\Int B(\epsilon_i)$ and for 
some $\epsilon_i>0\in\R$ and some $\alpha=(\alpha_1,...,\alpha_n)\in\N^n$ depending on $i$ with 
$\alpha_1>0,...,\alpha_{n'}>0,\alpha_{n'+1}=\cdots\alpha_n=0$ and that $U_i\cap X$ and $U'_i$ are 
carried to $\Int B(\epsilon_i)\cap\{x_1\cdots x_{n'}=0\}$ and $B(\epsilon_i/2)$ by the diffeomorphism 
of equivalence. 
Then by induction on $i$ it suffices to prove the following statement (for simplicity of notation we 
assume that $\epsilon_i=3$ and $\overline{U'_i}$ is carried to $B(1)$). \par
Let $C$ be a closed subset of $B(3)$. 
Let $f$ and $g$ be $C^\infty$ functions on $\R^n$ such that $f$ is of the form $x^\alpha=\prod_{j=1}
^nx_j^{\alpha_j}$ for the above $\alpha$ and $g$ is of the form $x^\alpha g'$ for some positive $C^
\infty$ function $g'$ on $\R^n$. 
Let $\pi$ be a $C^2$ embedding of $B(3)$ into $\R^n$ such that $f\circ\pi=g$ on $B(3)$ and $\pi(X\cap 
B(3))\subset X$ where $X=\{x^\alpha=0\}$. 
Let $\tau:B(3)\to\R^n$ be a $C^2$ approximation of $\pi$ in the $C^1$ topology such that $\tau(X\cap 
B(3))\subset X$, such that $f\circ\tau=g$ on a neighborhood of $C$ in $B(3)$ and $\tau$ is of class $C^\infty$ 
there. 
Then, fixing on $(B(3)-B(2))\cup C$, we can approximate $\tau$ by a $C^2$ embedding $\tilde\tau:B(3)\to
\R^n$ in the $C^1$ topology so that $\tilde\tau(X\cap B(3))\subset X$,
so that $f\circ\tilde\tau=g$ on $B(1)$ 
and $\tilde\tau$ is of class $C^\infty$ on $B(1)$. \par
We prove the statement. 
Set $\tau(x)=(\tau_1(x),...,\tau_n(x))$. 
Then $\tau_j(x)$ for each $1\le j\le n'$ is divisible by $x_j$, to be precise, there exists a 
positive $C^1$ function $F_j$ on $B(3)$ such that $\tau_j(x)=x_jF_j(x)$ since $\pi(X\cap B(3))\subset 
X$ and $ X=\{0\}\times\R^{n-1}\cup\cdots\cup\R^{n'-1}\times\{0\}\times\R^{n-n'}$ and $\pi$ is close to id. 
The required approximation $\tilde\tau=(\tilde\tau_1,...,\tilde\tau_n)$ also has to have the form $(x_
1\tilde F_1,...,x_{n'}\tilde F_{n'},\tilde\tau_{n'+1},...,\tilde\tau_n)$ for some positive $C^1$ 
functions $\tilde F_j$ and $C^2$ functions $\tilde\tau_{n'+1},...,\tilde\tau_n$. 
Set $F=(F_1,...,F_{n'})$ and $\tilde F=(\tilde F_1,...,\tilde F_{n'})$. 
Then $F$ is of class $C^\infty$ on a neighborhood of $C$, the condition $f\circ\tilde\tau=g$ on $B(1)$ 
coincides with the one $\tilde F^\alpha=g'$ on $B(1)$, and the other conditions which $\tilde F,
\tilde\tau_{n'+1},...,\tilde\tau_n$ satisfy are that $\tilde F=F$ on
$(B(3)-B(2))\cup C$, that $(\tilde F,
\tilde\tau_{n'+1},...,\tilde\tau_n)$ is an approximation of $(F,\tau_{n'+1},...,\tau_n)$ in the $C^1$ 
topology and that $\tilde\tau$ is of class $C^2$ on $B(3)$ and of class $C^\infty$ on $B(1)$. \par
Set $Z=\{(x,y)\in B(3)\times\R^{n'}:y^\alpha=g'(x)\}$, which is a $C^\infty$ submanifold with 
boundary of $B(3)\times\R^{n'}$ by the implicit function theorem since $g'$ is positive. 
Note that $\tilde F^\alpha=g'$ on $B(1)$ if and only if $\graph\tilde F|_{B(1)}\subset Z$ and that 
$\graph F|_C\subset Z$. 
We can construct a $C^\infty$ projection $p:W\to Z$ of a tubular neighborhood of $Z$ in $B(3)\times\R^{n'}$ 
such that $p(x,y)$ for $(x,y)\in W$ is of the form $(x,p_2(x,y))$ as follows. 
Since $g'$ is positive, $Z\cap\{x\}\times\R^{n'}$ for each $x\in B(3)$ is smooth and, moreover, the 
restriction to $Z$ of the projection $B(3)\times\R^{n'}\to B(3)$ is submersive. 
Hence if we define $p(x,y)$ for each $(x,y)\in B(3)\times\R^{n'}$ near $Z$ to be the orthogonal 
projection image of $(x,y)$ to $Z\cap\{x\}\times\R^{n'}$, then $p$ satisfies the requirements. 
Let $(\hat F,\hat\tau_{n'+1},...,\hat\tau_n)$ be a $C^\infty$ approximation of $(F,\tau_{n'+1},...,
\tau_n)$ in the $C^1$ topology, fixed on a neighborhood of $C$, and $\phi$ a $C^\infty$ 
function on $B(3)$ such that $0\le\phi\le1$, $\phi=1$ on $B(1)$ and $\phi=0$ on $B(3)-B(2)$. 
Define a $C^2$ map $\tilde F=(\tilde F_1,...,\tilde F_{n'}):B(3)\to\R^{n'}$ by
$$
\tilde F(x)=\phi(x)p_2(x,\hat F(x))+(1-\phi(x))F(x)\quad\text{for}\ x\in B(3),
$$
and set $\tilde\tau=(x_1\tilde F_1,...,x_{n'}\tilde F_{n'},\hat\tau_{n'+1},...,\hat\tau_n)$ on $B(3)$. 
Then $\graph\tilde F|_{B(1)}$ is included in $Z$ because $\tilde F|_{B(1)}$ coincides with the map $:B(1)\ni x
\to p_2(x,\tilde F(x))\in\R^{n'}$ whose graph is contained in $Z$;
then $\tilde F=F$ on $B(3)-B(2)$ since 
$\phi=0$ there; then $\tilde F=F$ on $C$ since $\hat F=F$ there and
since $p(x,F(x))=(x,F(x))$ there; then $(\tilde F,\hat\tau_{n'+1},...,\hat\tau_n)$ is an approximation of $(F,\tau_{n'+1},...,\tau_n)$ in 
the $C^1$ topology since so is $(\hat
F,,\hat\tau_{n'+1},...,\hat\tau_n)$; then $\tilde\tau$ is of class 
$C^2$ because if we set $p_2(x,y)=(p_{2,1}(x,y),...,p_{2,n'}(x,y))$ then 
$$
\tilde\tau_j(x)=\phi(x)x_jp_{2,j}(x,\hat F(x))+(1-\phi(x))\tau_j(x),\ 1\le j\le n'; 
$$
finally $\tilde\tau$ is of class $C^\infty$ on $B(1)$ since $\tilde F(x)=p_2(x,\hat F(x))$ on $B(1)$. 
Thus the statement is proved. \par
In conclusion, for some closed neighborhood $V$ of $f(\Sing f)$ in $\R$ each of whose connected 
components contains one point of $f(\Sing f)$, there exists a $C^2$ diffeomorphism $\tau$ of $M$ 
sufficiently close to $\pi$ in the Whitney $C^1$ topology such that $\tau$ is of class $C^\infty$ on 
$f^{-1}(V)$ and $f\circ\tau=g$ on $f^{-1}(V)$. 
Then the restrictions of $f$ and $g$ to $f^{-1}(\overline{\R-V})$ are proper and locally trivial maps 
onto $\overline{\R-V}$, moreover $f\circ\pi=g$ on $f^{-1}(\overline{\R-V})$ and $\tau|_{f^{-1}(\overline{\R-V}
)}$ is an approximation of $\pi|_{f^{-1}(\overline{\R-V})}$ in the Whitney $C^1$ topology. 
Hence we can modify $\tau$ so that $f\circ\tau=g$ and $\tau$ is of class $C^\infty$ everywhere 
fixing on $f^{-1}(V)$. 
Therefore, $f$ and $g$ are $C^\infty$ right equivalent, which proves (2). \par

{\it Proof of (3).} 
Let $0\ll l\in\N$. 
We prove first that $f$ and $g$ are semialgebraically $C^l$ right equivalent and later that semialgebraic 
$C^l$ right equivalence implies Nash right equivalence. 
We proceed with the former step as in the above proof of (2). 
Let $\pi$ be a semialgebraic $C^2$ diffeomorphism of $M$ such that $f\circ\pi=g$, and set $X=f^{-1}(f
(\Sing f))$ and $g^{-1}(g(\Sing g))$. 
Let $\pi'$ be a Nash approximation of $\pi$ in the semialgebraic $C^2$ topology (Approximation 
Theorem I). 
Then $\pi'$ is a diffeomorphism of $M$ and $\pi^{\prime-1}\circ\pi$ is a semialgebraic $C^2$ 
approximation of id in the semialgebraic $C^2$ topology. 
Hence by replacing $f$ and $\pi$ with $f\circ\pi'$ and $\pi^{\prime-1}\circ\pi$, we assume that $\pi$ is 
close to id in the semialgebraic $C^2$ topology. 
Moreover, we suppose that $X=Y$ as in the proof of (2) by using lemma 4.3 and its remark in place of lemma 
4.2. Furthermore, by using lemma 4.6 we can reduce the problem to the case where $M$ is the interior of a 
compact Nash manifold possibly with boundary $M_1$ and for each
$x\in\partial M_1$, the germ $(M_{1x},X_x)$ is 
Nash diffeomorphic to the germ at 0 of $(\R^{n-1}\times[0,\,\infty),\{(x_1,...,x_{n-1})\in\R^{n-1}:x_
1\cdots x_{n'}=0\}\times(0,\,\infty))$ for some $n'\,(<n)\in\N$. \par
We modify $\pi$ on a semialgebraic neighborhood of $X$. 
By lemma 4.7 and proposition 4.8,(iii) there exist {\bf finite} open semialgebraic coverings $\{U_i\}$ 
and $\{U'_i\}$ of $X$ in $M$ such that the closure $\overline{U'_i}$
in $M$ is contained in $U_i$, such that $\pi(
\overline{U'_i})$ is contained in $U_i$, such that $f|_{U_i}$ is Nash right equivalent to $x^\alpha+\,$constant on 
$\Int B_{\xi_i}(\epsilon_i)$ where $\alpha=(\alpha_1,...,\alpha_n)\in\N^n$ depending on $i$ with $\alpha
_1>0,...,\alpha_{n'}>0,\alpha_{n'+1}=\cdots\alpha_n=0$, for $n'\,(<n)\in\N-\{0\}$ and $B_{\xi_i}(\epsilon_i)=
\{x=(x_1,...,x_n)\in\R^n:x_n>0,\,|x^\alpha|\le\xi_i(x_n),\,|x|\le\epsilon_i\}$ for some $\epsilon_i>0
\in\R$ and some positive Nash function $\xi_i$ on $(0,\,\infty)$ and that $U_i\cap X$ and $U'_i$ are 
carried to $\Int B_{\xi_i}(\epsilon_i)\cap\{x_1\cdots x_{n'}=0\}$ and $\Int B_{\xi_i/2}(\epsilon_i/2)$ 
by the diffeomorphism of equivalence. 
For modification of $\pi$ on a semialgebraic neighborhood of $X$ we need the following statement. 
Let $l'\in\N$ such that $l\le l'\le l+\#\{i\}$. \par
Let $\xi$ be a small positive Nash function on $(0,\,\infty)$, and $C$ a closed semialgebraic 
subset of $B_{3\xi}(3)$. 
Let $f$ and $g$ be Nash functions on $B_{4\xi}(4)$ such that $f$ is of the form $x^\alpha$ for the 
above $\alpha$ and $g$ is of the form $x^\alpha g'$ for some positive Nash function $g'$ on $B_{4\xi
}(4)$. 
Let $\pi$ be a semialgebraic $C^2$ embedding of $B_{3\xi}(3)$ into $B_{4\xi}(4)$ close to id in the 
semialgebraic $C^2$ topology such that $f\circ\pi=g$ on $B_{3\xi}(3)$ and $\pi(X\cap B_{3\xi}(3))
\subset X$ where $X=\{(x_1,...,x_n)\in\R^n:x_1\cdots x_{n'}=0\}$. 
Let $\tau:B_{3\xi}(3)\to B_{4\xi}(4)$ be a semialgebraic $C^2$ approximation of $\pi$ in the 
semialgebraic $C^1$ topology such that $\tau(X\cap B_{3\xi}(3))\subset
X$, such that $ f\circ\tau=g$ on a closed 
semialgebraic neighborhood $V$ of $C$ in $B_{3\xi}(3)$ and $\tau$ is of class $C^{l'}$ there. 
Then, fixing on $(B_{3\xi}(3)-B_{2\xi}(2))\cup C$ we can approximate $\tau$ by a semialgebraic $C^2$ 
embedding $\tilde\tau:B_{3\xi}(3)\to B_{4\xi}(4)$ in the semialgebraic $C^1$ topology so that 
$\tilde\tau(X\cap B_{3\xi}(3))\subset X$, such that $ f\circ\tilde\tau=g$ on $B_\xi(1)$ and $\tilde\tau$ is of 
class $C^{l'-1}$ on $B_\xi(1)$. \par
We prove the statement. 
As before, set $\tau=(\tau_1,...,\tau_n)$, $\tilde\tau=(\tilde\tau_1,...,\tilde\tau_n)$, let $F_j$ 
and $\tilde F_j,\,1\le j\le n',$ be positive semialgebraic $C^1$ functions on $B_{3\xi}(3)$ such that 
$\tau_j=x_jF_j(x)$ and $\tilde\tau_j=x_j\tilde F_j(x)$ on $B_{3\xi}(3)$, and set $F=(F_1,...,F_{n'})$ 
and $\tilde F=(\tilde F_1,...,\tilde F_{n'})$. 
Note that $F_j$ are of class $C^{l'-1}$ on a semialgebraic neighborhood of $C$ in $B_{3\xi}(3)$, which 
is different to $F_j$ in the proof of (2) where they are of class $C^\infty$. 
Then the required conditions are that $\tilde F^\alpha=g'$ on
$B_\xi(1)$, that $\tilde F=F$ on $(B_{3\xi}(3)-
B_{2\xi}(2))\cup C$, that $(\tilde F,\tilde\tau_{n'+1},...,\tilde\tau_n)$ is a semialgebraic $C^1$ 
approximation of $(F,\tau_{n'+1},...,\tau_n)$ in the semialgebraic $C^1$ topology, and $\tilde\tau$ 
is of class $C^2$ on $B_{3\xi}(3)$ and of class $C^{l'-1}$ on $B_\xi(1)$. \par
Set $Z=\{(x,y)\in B_{3\xi}(3)\times\R^{n'}:y^\alpha=g'(x)\}$, which is a Nash submanifold with 
boundary of $B_{3\xi}(3)\times\R^{n'}$, and let $p:W\to Z$ be a Nash projection of a semialgebraic 
tubular neighborhood of $Z$ in $B_{3\xi}(3)\times\R^{n'}$ such that $p(x,y)$ for $(x,y)\in W$ is of 
the form $(x,p_2(x,y))$, which is constructed as before. 
Let $(\hat F,\hat\tau_{n'+1},...,\hat\tau_n)$ be a Nash approximation of $(F,\tau_{n'+1},...,\tau_n)$ 
in the semialgebraic $C^1$ topology, and $\phi$ and $\psi$ semialgebraic $C^{l'}$ functions on 
$B_{3\xi}(3)$ such that $0\le\phi\le1$, such that $\phi=1$ on
$B_\xi(1)$ and $\phi=0$ on $B_{3\xi}(3)-B_{2\xi}(2)$, such that 
$0\le\psi\le1$ and $\psi=1$ on $B_{3\xi}(3)-V$ whereas $\psi=0$ on a semialgebraic neighborhood of $C$ in $B
_{3\xi}(3)$ smaller than $\Int V$. 
Set 
$$
\gather
\tilde F(x)=\phi(x)p_2\big(x,\psi(x)\hat F(x)+(1-\psi(x))F(x)\big)+(1-\phi(x))F(x)\quad\text{for}\ x\in B_{
3\xi}(3),\\
\tilde\tau=(x_1\tilde F_1,...,x_{n'}\tilde F_{n'},\hat\tau_{n'+1},...,\hat\tau_n)\quad\text{on}\ B_
{3\xi}(3).\tag"and"
\endgather
$$
Then we see as before that the required conditions are satisfied. 
Hence the statement is proved. \par
By the statement, a partition of unity of class semialgebraic $C^l$ and by remark 2.11,(5)$'$ we obtain an open 
semialgebraic neighborhood $U$ of $X$ and a semialgebraic $C^2$ diffeomorphism $\tau$ of $M$ close 
to $\pi$ in the semialgebraic $C^1$ topology such that $\tau$ is of class $C^l$ on $U$ and $f\circ
\tau=g$ on $U$ (the point is that after fixing $U$ we can choose $\tau$ so as to be arbitrarily close 
to id). 
Then we modify $\tau$ so that $\tau$ is of class semialgebraic $C^l$ and $f\circ\tau
=g$, i.e., $f$ and $g$ are semialgebraically $C^l$ right equivalent as follows. \par
Let $\eta$ be a semialgebraic $C^l$ function on $M$ such that
$0\le\eta\le1$, such that $\eta=0$ outside of 
$U$ and $\eta=1$ on a smaller semialgebraic neighborhood of $X$, and set $A=\{(x,y)\in(M-X)^2:f(y)=
g(x)\}$. 
Then $A$ is a Nash manifold and there exists a Nash projection $q:Q\to A$ of a small semialgebraic 
tubular neighborhood of $A$ in the square of the ambient Euclidean space of $M$ of the form $q(x,y)=(x,
q_2(x,y))$ for $x\in M-X$. 
Let $\check\tau$ be a Nash approximation of $\tau$ in the semialgebraic $C^1$ topology, and set 
$$
\Check{\Check\tau}=q_2\big(x,\eta(x)\tau(x)+(1-\eta(x))\check\tau(x)\big)\quad\text{for}\ x\in M. 
$$
Then $\Check{\Check\tau}$ is well-defined because the graph of the map from $M$ to the ambient Euclidean 
space of $M$ $:x\to\eta(x)\tau(x)+(1-\eta(x))\check\tau(x)$ is contained in $Q$, hence $\Check{
\Check\tau}$ is a semialgebraic $C^l$ diffeomorphism of $M$ and $f\circ\Check{\Check\tau}=g$. 
Thus the former step of the proof is achieved. \par
Let $0\ll l^{(3)}\ll\cdot\cdot\ll l\in\N$. 
For the latter step also we can assume that $X=Y$ and that there exists a semialgebraic $C^l$ 
diffeomorphism $\pi$ of $M$ close to id in the semialgebraic $C^l$ topology such that $f\circ\pi=g$. 
Let $\mu$ be a Nash function on $\R$ such that $\mu^{-1}(0)=f(\Sing f)$ and $\mu$ is regular at 
$\mu^{-1}(0)$. 
Consider $\mu\circ f$ and $\mu\circ g$. 
Their zero sets are $X$, they have only normal crossing singularities at $X$, the same sign at each 
point of $M$ and the same multiplicity at each point of $X$, and we see easily that the Nash 
function on $M$, defined to be $\mu\circ g/\mu\circ f$ on $M-X$, is close to 1 in the semialgebraic 
$C^l$ topology. 
Hence the conditions in lemma 4.7 are satisfied and there exists a Nash diffeomorphism $\pi'$ of 
$M$ close to id in the semialgebraic $C^{l'}$ topology such that $\pi'(X)=X$ and $f\circ\pi'-g$ is 
$l'$-flat at $X$. 
Thus, replacing $f$ and $\pi$ with $f\circ\pi'$ and $\pi^{\prime-1}\circ\pi$, we assume that $f-g$ is 
$l'$-flat at $X$ and $\pi$ is close to id in the semialgebraic $C^{l'}$ topology. \par
By proposition 4.9 we can assume that $M$ is the interior of a compact Nash manifold possibly with corners 
$M_1$ and $f$ is the restriction to $M$ of a Nash function $f_1$ on $M_1$ with only normal crossing 
singularities. 
Then by the definition of semialgebraic $C^l$ topology, $\pi$ is extensible to a semialgebraic $C^{
l'}$ diffeomorphism $\pi_1$ of $M_1$ such that $\pi_1-\id$ is $l'$-flat at $\partial M_1$. 
Hence $g$ also is extensible to a semialgebraic $C^{l'}$ function $g_1$ on $M_1$, and $f_1-g_1$ is 
close to 0 in the $C^{l'}$ topology and $l'$-flat at $\partial M_1$. 
Let $v_i$, for $ i=1,...,N$, be Nash vector fields on $M_1$ spanning the tangent space of $M_1$ at each 
point, $\nu_1$ a non-negative Nash function on $M_1$ with zero set $\partial M_1$ and regular 
there, and set $\nu_2=\sum_{i=1}^N(v_if_1)^2$ and $\nu=\nu_1^{l''}\nu_2$. 
Then the radical of $\nu_2\Cal N$ is the sheaf of $\Cal N$-ideals defined by $X\cup\partial M_1$, and 
$f_1-g_1$ is divisible by $\nu$; to be precise, there exists a semialgebraic $C^{l''}$ function 
$\beta$ on $M_1$ such that $f_1-g_1=\nu\beta$. 
Moreover, $\beta$ is close to 0 in the $C^{l^{(3)}}$
topology. Actually, by lemma 2.12 the map $C^\infty(M_1)\ni h\to\nu h\in\nu C^\infty(M_1)$ is open. 
Hence for $h\in C^\infty(M_1)$, if $\nu h$ is close to 0 in the $C^{l''}$ topology then $h$ is close 
to 0 in the $C^{l^{(3)}}$ topology. 
This holds for $h\in C^{l''}(M_1)$ also because $h$ of class $C^{l''}$ is approximated by a $C^\infty$ 
function $h'$ in the $C^{l''}$ topology and $\nu\cdot(h-h'$) and hence $\nu h'$ are close to 0 in the 
$C^{l''}$ topology. Therefore, $\beta$ is close to 0 in the $C^{l^{(3)}}$ topology. \par
It follows from the definition of semialgebraic $C^l$ topology that $\nu_1^{l''}\beta|_M$ is close 
to 0 in the semialgebraic $C^{l^{(3)}}$ topology. 
Then the conditions in proposition 4.8,(ii) for $f$ and $g\,(=f-\nu_1^{l''}\beta\sum_{i=1}^N(v_if_1)^2
|_M)$ are satisfied. 
Hence $f$ and $g$ are Nash right equivalent. \par
We can prove the case with corners in the same way. 
\qed

\enddocument